\documentclass{gtart_h}  

\def\ifplaintex{\expandafter\ifx\csname documentclass\endcsname\relax}

\ifplaintex 
\else
\fi

\expandafter\ifx\csname beginpicture\endcsname\relax
\expandafter\ifx\csname documentclass\endcsname\relax
\input pictex \else
\input prepictex \input pictex \input postpictex \fi\fi

\def\gt{{\mathsurround=0pt\it $\cal G\mskip-2mu$eometry \&\ 
$\cal T\!\!$opology}}        

\def\gtp{{\mathsurround=0pt\it $\cal G\mskip-2mu$eometry \&\ 
$\cal T\!\!$opology $\cal P\!$ublications}}  

\def\lognumber#1{\def\thelognumber{#1}}
\def\volumenumber#1{\def\thevolumenumber{#1}}
\def\papernumber#1{\def\thepapernumber{#1}}
\def\volumeyear#1{\def\thevolumeyear{#1}}

\def\pagenumbers#1#2{\def\startpage{#1}\def\finishpage{#2}}
\def\published#1{\def\publishdate{#1}}
\def\proposed#1{\def\theproposer{#1}}
\def\seconded#1{\def\theseconders{#1}}
\def\received#1{\def\receiveddate{#1}}
\def\revised#1{\def\reviseddate{#1}}
\def\accepted#1{\def\accepteddate{#1}}
\def\asciititle#1{\def\theasciititle{#1}}
\def\covertitle#1{\def\thecovertitle{#1}}
\def\coverauthors#1{\def\thecoverauthors{#1}}
\def\asciiauthors#1{\def\theasciiauthors{#1}}
\def\asciiaddress#1{\def\theasciiaddress{#1}}

\long\def\asciiabstract#1{\long\def\theasciiabstract{#1}}
\def\asciikeywords#1{\def\theasciikeywords{#1}}

\let\\\par\let\thelognumber\relax
\let\thevolumenumber\relax\let\thepapernumber\relax
\let\thevolumeyear\relax\let\thesamplenumber\relax\let\startpage\relax
\let\finishpage\relax\let\publishdate\relax\let\receiveddate\relax
\let\reviseddate\relax\let\accepteddate\relax\let\theasciititle\relax
\let\thecovertitle\relax\let\theasciiauthors\relax\let\theasciiaddress\relax
\let\theasciiabstract\relax\let\theasciikeywords\relax
\let\theasciiemail\relax\let\theshortauthors\relax\let\theshorttitle\relax
\let\thecoverauthors\relax

\long\def\maketitlep{   

\count0=\startpage

\gt\hfill      
\beginpicture
\setcoordinatesystem units <0.33truein, 0.33truein> point at 2.2 0.9
\setplotsymbol ({$\cal G$})
\plotsymbolspacing=9truept
\circulararc 315 degrees from 0 1 center at 0 0
\setplotsymbol ({$\cal T$})
\circulararc 315 degrees from 1 -1 center at 1 0
\endpicture
\break
{\small\ifx\thesamplenumber\relax 
Volume \else Sample
\fi\thevolumenumber\ (\thevolumeyear)
\startpage--\finishpage\nl
Published: \publishdate}
\vglue 0.5truein plus 0.4fil minus 0.1truein

{\parskip=0pt\leftskip 0pt plus 1fil\def\\{\par\smallskip}{\ifplaintex\large
\else\Large\fi\bf\thetitle}\par\medskip}   

\vglue 0pt plus 0.1fil 

{\parskip=0pt\leftskip 0pt plus 1fil\def\\{\par}{\sc\theauthors}
\par\medskip}

\vglue 0pt plus 0.1fil 

{\small\parskip=0pt\let\newline\\
{\leftskip 0pt plus 1fil\def\\{\par}{\sl\theaddress}\par}
\expandafter\ifx\theemail\relax    
\relax\else\vglue 5pt plus 0.02fil minus 2pt\def\\{\stdspace{\rm 
and}\stdspace} 
\cl{Email:\stdspace\tt\theemail}\fi
\ifx\theurl\relax                  
\relax\else\vglue 5pt plus 0.02fil minus 2pt\def\\{\stdspace{\rm 
and}\stdspace}
\cl{URL:\stdspace\tt\theurl}\fi\par}

\vglue 7pt plus 0.3fil minus 3pt

{\bf Abstract}
\vglue 5pt plus 0.1fil minus 2pt

\theabstract

\vglue 7pt plus 0.3fil minus 3pt

{\bf AMS Classification numbers}\quad Primary:\quad \theprimaryclass

Secondary:\quad \thesecondaryclass

\vglue 5pt plus 0.3fil minus 2pt

{\bf Keywords:}\quad \thekeywords

\vglue 10pt plus 0.5fil minus 5pt

{\small  Proposed: \theproposer\hfill Received: \receiveddate\nl
Seconded: \theseconders\hfill 
\ifx\reviseddate\relax                         
Accepted: \accepteddate                        
\else
Revised: \reviseddate                          
\fi}
\eject
}       

\let\maketitlepage\maketitlep
\let\maketitle\maketitlepage

\font\phead=cmsl9 scaled 950
\font=cmsl9 scaled 1050
\font\pnum=cmbx10 scaled 913
\font\lnum=cmbx10 
\font\pfoot=cmsl9 scaled 950
\font=cmsl9 scaled 1050
\ifplaintex
\headline{\vbox to 0pt{\vskip -4.5mm\line{\small\phead\ifnum
\count0=\startpage ISSN 1364-0380 (on line)
1465-3060 (printed) \hfill {\pnum\folio}\else\ifodd\count0\def\\{ }%
\ifx\theshorttitle\relax\thetitle\else\theshorttitle\fi\hfill{\pnum\folio}
\else\def\\{ and }{\pnum\folio}\hfill\ifx\theshortauthors\relax\theauthors
\else\theshortauthors\fi\fi\fi}\vss}}
\footline{\vbox to 0pt{\vglue 0mm\line{\small\pfoot\ifnum\count0=\startpage
\copyright\ \gtp\hfill\else
\gt, Volume \thevolumenumber\ (\thevolumeyear)\hfill\fi}\vss
}}
\else
\makeatletter
\def\@oddhead{{\small\lhead\ifnum\count0=\startpage ISSN 1364-0380 (on line)
1465-3060 (printed) \hfill {\lnum\number\count0}\else\ifodd\count0
\def\\{ }\ifx\theshorttitle\relax \thetitle \else\theshorttitle\fi\hfill
{\lnum\number\count0}\else\def\\{ and }{\lnum\number\count0}
\hfill\ifx\theshortauthors\relax 
\theauthors\else\theshortauthors\fi\fi\fi}}\def\@evenhead{\@oddhead}
\def\@oddfoot{\small\lfoot\ifnum\count0=\startpage\copyright\ \gtp\hfill\else
\gt, Volume \thevolumenumber\ (\thevolumeyear)\hfill\fi}
\def\@evenfoot{\@oddfoot}
\makeatother
\fi

\newwrite\gtoutfile
\long\gdef\makeheadfile{  
{\def\\{, }\def\s{ }
\immediate\openout\gtoutfile head.xxx
\immediate\write\gtoutfile{Proxy-for: \ifx\theasciiauthors\relax
\theauthors\else\theasciiauthors\fi\s<\ifx\theasciiemail\relax\theemail\else\theasciiemail\fi>}
\immediate\write\gtoutfile{\noexpand\\}
\immediate\write\gtoutfile{Authors: \ifx\theasciiauthors\relax
\theauthors\else\theasciiauthors\fi}
{\def\\{ }\immediate\write\gtoutfile{Title: \ifx\theasciititle\relax
\thetitle\else\theasciititle\fi}}
\immediate\write\gtoutfile{Subj-class: GT or SG or MG etc}
\immediate\write\gtoutfile{MSC-class: \theprimaryclass\ifx\thesecondaryclass\relax\else, \thesecondaryclass\fi}
\immediate\write\gtoutfile{Journal-ref: Geom. Topol. \thevolumenumber
(\thevolumeyear) \startpage-\finishpage}
\immediate\write\gtoutfile{Comments: Published by Geometry and Topology at}
\immediate\write\gtoutfile{\s\s http://www.maths.warwick.ac.uk/gt/GTVol\thevolumenumber/paper\thepapernumber.abs.html}
\immediate\write\gtoutfile{\noexpand\\}
\immediate\write\gtoutfile{}
\ifx\theasciiabstract\relax
\immediate\write\gtoutfile{\theabstract}\else
\immediate\write\gtoutfile{\theasciiabstract}\fi
\immediate\write\gtoutfile{}
\immediate\write\gtoutfile{\noexpand\\}
\immediate\write\gtoutfile{}
\immediate\closeout\gtoutfile}}  

\def\maketitlepage{\maketitlep\makeheadfile}
\let\maketitle\maketitlepage

\lognumber{369}
\received{9 October 2003}
\volumenumber{9}\papernumber{1}\volumeyear{2005}
\pagenumbers{1}{93}   
\revised{1 December 2004}
\published{17 December 2004}
\accepted{29 September 2004}
\proposed{Dieter Kotschick}
\seconded{Simon Donaldson, Tomasz Mrowka}

\usepackage{amsmath}
\usepackage{amsfonts}
\usepackage{amssymb}

\newtheorem{lemma}{Lemma}[section]

\newtheorem{prop}{Proposition}[section]
\newtheorem{thm}{Theorem}[section]
\newtheorem{cor}{Corollary}[section]
\newtheorem{defn}{Definition}[section]
\newtheorem{claim}{Claim}[section]

\newcommand{\col}{\colon\thinspace}

\newcommand{\Ref}[1]{(\ref{#1})}

\newcommand{\cb}{\mathcal B}
\def\cc{\mathcal C}
\def\cD{\mathcal D}
\def\ce{\mathcal E}
\newcommand{\cg}{\mathcal G}
\def\harm{\mathcal H}
\def\ci{\mathcal I}
\def\cll{\mathcal L}
\newcommand{\cm}{\mathcal M}
\def\cO{\mathcal O}

\def\fl{\mathcal R}
\def\ct{\mathcal T}

\def\cw{\mathcal W}

\newcommand{\al}{\alpha}
\newcommand{\del}{\delta}
\def\Del{\Delta}
\newcommand{\eps}{\epsilon}
\def\ga{\gamma}
\def\Ga{\Gamma}
\def\ka{\kappa}
\def\lla{\lambda}
\def\La{\Lambda}
\def\om{\omega}
\def\Om{\Omega}
\newcommand{\si}{\sigma}
\def\Si{\Sigma}

\def\BB{\mathbf B}

\renewcommand{\co}{\mathbb C}

\def\bl{\mathbb L}
\newcommand{\n}{\mathbb N}
\def\q{{\mathbb Q}}
\newcommand{\R}{\mathbb R}
\def\bs{\mathbb S}

\def\bv{\mathbb V}
\newcommand{\z}{\mathbb Z}

\newcommand{\uset}{\underset}
\newcommand{\oset}{\overset}
\newcommand{\uline}{\underline}
\newcommand{\oline}{\overline}

\newcommand{\la}{\langle}
\newcommand{\ra}{\rangle}
\newcommand{\st}{\,:\,}
\newcommand{\ti}{\widetilde}

\newcommand{\prtl}{\partial}
\def\square{\kern20pt{\vbox{\hrule height.4pt
        \hbox{\vrule width.4pt height 6pt\kern6pt
                \vrule width.4pt}
        \hrule height.4pt}}}
\newcommand{\ry}{\R\times Y}
\newcommand{\rpy}{\R_+\times Y}

\def\orp{\overline{\mathbb R}_+}
\def\orpy{\overline{\mathbb R}_+\times Y}

\def\im{\text{im}}
\def\Im{\text{Im}}

\def\inv{^{-1}}

\def\aut{\text{Aut}}
\def\const{\text{const}\cdot}
\def\rf{{\R}^4}
\def\loc{{\text{loc}}}

\def\U#1{\text{U}(#1)}

\def\SO#1{\text{SO}(#1)}

\def\spc{\text{spin}^c}
\def\Spc{\text{Spin}^c}

\def\Spin#1{\text{Spin}(#1)}
\def\Spinc#1{\text{Spin}^c(#1)}

\def\GL{\text{GL}}
\def\glp#1{\text{GL}^+(#1)}
\def\glc#1{\text{GL}^c(#1)}
\def\tglp#1{\widetilde{\text{GL}}\vphantom{L}^+(#1)}
\def\pso{P_{\text{SO}}}
\def\pglp{P_{\text{GL}^+}}
\def\pglc{P_{\text{GL}^c}}
\def\pspc{P_{\text{Spin}^c}}

\def\ind{\text{\rm ind}}
\def\index{\text{index}}

\def\cutoff{c}
\def\g{\mathfrak g}

\def\fq{\mathfrak q}
\def\s{\mathfrak s}

\def\cd{{\vartheta}}
\def\csd#1{\cd_{#1}}

\def\cocl{\mathbb C\ell}

\def\loc{\text{loc}}

\def\llw#1#2#3{{L^{#1,#2}_#3}}

\def\lw#1#2{{L^{#1,#2}}}

\def\id{\text{Id}}

\def\sdot{\,\cdot\,}
\def\pp#1{\frac\partial{\partial#1}}

\def\lpk{{L^{p,\kappa}}}
\def\lpkn{{L^{p,\kappa_n}}}
\def\llpk{{L^{p,\kappa}_1}}
\def\llpkn{{L^{p,\kappa_n}_1}}
\def\lpkj{{L^{p,\kappa_j}}}
\def\llpke{{L^{p,\kappa_e}_1}}

\def\bbe{\bar\beta_e}
\def\tu{^{(T)}}
\def\xt{{X\tu}}
\def\mt{M\tu}

\def\rr{_o}

\def\ual{{\uline\al}}
\def\adt{a'}
\def\uh{\uline h}

\def\av{_{\text{av}}}

\def\hx{{\mathfrak b}}
\def\hatf{\hat F}

\def\tcd{\ti\cd}
\def\sj{_{(j)}}
\def\Sj{^{(j)}}
\def\nj{_{n,j}}
\def\tno{t_{n,1}}
\def\vin{v_{i,n}}
\def\vjn{v_{j,n}}
\def\wkn{w_{k,n}}
\def\prt{\mathfrak p}
\def\Prt{\mathfrak P}
\def\tPrt{\ti{\mathfrak P}}
\def\scal{\mathbf s}
\def\quadr{Q}
\def\he#1{^{#1}}
\def\hg#1{^{[#1}}
\def\gl{^\#}

\def\endt#1{_{:#1}}
\def\pxt{X^{\{T\}}}
\def\xtn{X^{(T(n))}}
\def\ortn{o}
\def\katn{\ka_n}
\def\tmin{\check T}
\def\wt{W^{(T)}}
\def\wtn{W^{T(n)}}
\def\hatfpl{\hat F\vphantom{F}^+}

\def\Pkd{\Pi_\del}
\def\Pid{\Pi_\del}
\def\SW{\text{SW}}
\def\BF{\widetilde{\text{SW}}}
\def\sw{\Theta}
\def\tsw{\ti\Theta}
\def\itr{\text{int}}

\def\uk{\Upsilon_k}
\def\mkd{M_{k,\del}}
\def\wkd{W_{k,\del}}
\def\skd{\si_{k,\del}}
\def\cbr{\cb_\rho}
\def\mkdr{M_{k,\del,\rho}}
\def\gkdr{G_{k,\del,\rho}}

\def\pb{_{\prtl\BB}}
\def\Ttn{\Theta_{\tau,n}}

\begin{document}

\title{Monopoles over $4$--manifolds containing long necks, I}
\covertitle{Monopoles over $4$--manifolds\\containing long necks, I}
\asciititle{Monopoles over 4-manifolds containing long necks, I}

\author{Kim A Fr\o yshov}
\coverauthors{Kim A Fr\noexpand\o yshov}
\asciiauthors{Kim A Froyshov}

\address{Fakult\"at f\"ur Mathematik, Universit\"at Bielefeld\\Postfach
100131, D-33501 Bielefeld, Germany}
\asciiaddress{Fakultaeat fuer Mathematik, Universitaet Bielefeld\\Postfach
100131, D-33501 Bielefeld, Germany}

\email{froyshov@mathematik.uni-bielefeld.de}

\begin{abstract}
We study moduli spaces of Seiberg--Witten monopoles over ${\rm
spin}^c$ Riemannian 4--manifolds with long necks and/or tubular
ends. This first part discusses compactness, exponential decay, and
transversality. As applications we prove two vanishing theorems for
Seiberg--Witten invariants.
\end{abstract}

\asciiabstract{%
We study moduli spaces of Seiberg-Witten monopoles over spin^c
Riemannian 4-manifolds with long necks and/or tubular ends. This first
part discusses compactness, exponential decay, and transversality. As
applications we prove two vanishing theorems for Seiberg-Witten
invariants.}

\primaryclass{57R58}

\secondaryclass{57R57}

\keywords{Floer homology, Seiberg--Witten, Bauer--Furuta, compactness,
monopoles}
\asciikeywords{Floer homology, Seiberg-Witten, Bauer-Furuta, compactness,
monopoles}

\maketitlepage

\section{Introduction}

This is the first of two papers devoted to the study of moduli spaces of
Seiberg--Witten
monopoles over $\spc$ Riemannian $4$--manifolds with long necks and/or 
tubular ends. Our principal motivation is to provide the 
analytical foundations for subsequent work on Floer homology.
Such homology groups should appear naturally when one attempts
to express the Seiberg--Witten invariants of a closed $\spc$
$4$--manifold $Z$ cut along a hypersurface $Y$, say
\[Z=Z_1\cup_YZ_2,\]
in terms of relative invariants
of the two pieces $Z_1,Z_2$. The standard approach, familiar from
instanton Floer theory (see \cite{F1,D5}), is to construct a $1$--parameter
family $\{g_T\}$ of Riemannian metrics on $Z$ by stretching along $Y$
so as to obtain a neck $[-T,T]\times Y$, and  study the monopole
moduli space $\mt$ over $(Z,g_T)$ for large $T$. There are
different aspects of this problem: compactness, transversality, and gluing. 
The present paper will focus particularly on compactness, and also establish
transversality results sufficient for the construction of Floer homology
groups of rational homology $3$--spheres. The second paper in this series
will be devoted to gluing theory.

Let the monopole equations over the neck $[-T,T]\times Y$
be perturbed by a closed $2$--form $\eta$
on $Y$, so that temporal gauge solutions to these equations correspond to
downward gradient flow lines of the correspondingly perturbed
Chern--Simons--Dirac functional $\csd\eta$ over $Y$.
Suppose all critical points of $\csd\eta$ are non-degenerate.
Because each moduli space $\mt$ is compact,
one might expect, by analogy with Morse theory,
that a sequence $\om_n\in M^{(T_n)}$
where $T_n\to\infty$ has a subsequence which
converges in a suitable sense to a pair of monopoles over the cylindrical-end
manifolds associated to $Z_1,Z_2$ together with a broken gradient line of
$\csd\eta$ over $\ry$.
The first results in this direction were obtained
by Kronheimer--Mrowka \cite{KM4} (with $\eta=0$) and
Morgan--Szab\'o--Taubes \cite{MST} (in a particular case, with $\eta$ non-exact).
Nicolaescu's book \cite{nico1}
contains some foundational results in the case $\eta=0$.
Marcolli--Wang \cite{Marcolli-Wang1} proved a general compactness theorem
for $\eta$ exact.
In this paper we will consider the general case when $\eta$ may be non-exact.
Unfortunately, compactness as stated above may then fail.
(A simple class of
counter-examples is described after Theorem~\ref{thm:neck-comp} below.)
It is then natural to seek
topological conditions which ensure that compactness does hold.
We will consider
two approaches which provide different sufficient conditions.
In the first approach, which is a refinement of well-known 
techniques (see \cite{D4,nico1}), one first establishes global bounds
on what is morally
a version of the energy functional (although the energy concept is ambiguous
in the presence of perturbations) and then derives local $L^2$ bounds on
the curvature forms. In the second approach, which appears to be new,
one begins by placing the
connections in Coulomb gauge with respect to a given reference connection
and then obtains global bounds on the corresponding connection forms 
in suitably weighted Sobolev norms, utilizing the \textit{a priori} pointwise
bounds on the spinors.

The paper also contains expository sections on configuration spaces and
exponential decay, borrowing some ideas from
Donaldson \cite{D5}, to which we also refer for the Fredholm theory.
In the transversality theory of moduli spaces
we mostly restrict ourselves, for the time being,
to the case when all ends of the $4$--manifold in question
are modelled on rational homology spheres.
The perturbations used here are minor modifications of the ones introduced
in \cite{Fr1}. 
It is not clear to us that these perturbations 
immediately carry over to the case of more general ends,
as has apparently been assumed
by some authors, although we expect that a modified version
may be shown to work with the aid of gluing theory.

In most of this paper we make an assumption on the cohomology class of
$\eta$ which rules out the hardest case 
in the construction of Floer homology (see Subsection~\ref{intro-csd}).
This has the advantage that we can use relatively simple perturbations.
A comprehensive monopole Floer theory including the hardest case
is expected to appear in a forthcoming book by Kronheimer--Mrowka \cite{KM5}.
An outline of their construction (and much more) can be found in
\cite{KMOS1}.

A large part of this work was carried out during a one-year stay at
the Institut des Hautes \'Etudes Scientifiques, and the author is grateful
for the hospitality and excellent research environment which he enjoyed there.
This work was also supported by a grant from the National Science Foundation.

\subsection{Vanishing results for Seiberg--Witten invariants}
\label{intro:vanishing}

Before describing our compactness results in more detail we will mention two
applications to Seiberg--Witten invariants of closed $4$--manifolds.

By a {\em $\spc$ manifold} we shall mean an oriented smooth manifold with a 
$\spc$ structure. If $Z$ is a $\spc$ manifold then $-Z$ will refer to 
the same smooth manifold equipped with the opposite orientation
and the corresponding $\spc$ structure.

If $Z$ is a closed, oriented $4$--manifold then by an
{\em homology orientation} of $Z$ we
mean an orientation of the real vector space
$H^0(Z)^*\oplus H^1(Z)\oplus H^+(Z)^*$, where $H^+(Z)$ is any maximal
positive subspace for the intersection form on $H^2(Z)$. The dimension of
$H^+$ is denoted $b^+_2$.

In \cite{Bauer-Furuta} Bauer and Furuta introduced a refined Seiberg--Witten
invariant for closed $\spc$ $4$--manifolds $Z$. This invariant $\BF(Z)$
lives in a certain equivariant stable cohomotopy group. 
If $Z$ is connected and $b^+_2(Z)>1$, and given
an homology orientation of $Z$, then according to \cite{Bauer3}
there is a natural homomorphism from this stable cohomotopy group to $\z$
which maps $\BF(Z)$ to the ordinary Seiberg--Witten invariant $\SW(Z)$
defined by the homology orientation. In \cite{bauer2} Bauer showed
that, unlike the ordinary
Seiberg--Witten invariant, the refined invariant does not in general
vanish for connected sums where both summands have $b^+_2>0$.
However, $\BF(Z)=0$ provided there exists a
metric and perturbation $2$--form on $Z$ for which the Seiberg--Witten moduli
space $M_Z$ is empty (see \cite[Remark 2.2]{Bauer3} and
\cite[Proposition~6]{Ishida-LeBrun1}).

\begin{thm}\label{sw-thm1}
Let $Z$ be a closed $\spc$ $4$--manifold
and $Y\subset Z$ a closed, orientable $3$--dimensional submanifold.
Suppose
\begin{enumerate}
\item[\rm(i)]$Y$ admits a Riemannian metric with positive scalar curvature, and
\item[\rm(ii)]$H^2(Z;\q)\to H^2(Y;\q)$ is non-zero.
\end{enumerate}
Then there exists a metric and perturbation $2$--form on $Z$ for which
$M_Z$ is empty, hence $\BF(Z)=0$.
\end{thm}

This generalizes a result of Fintushel--Stern \cite{FS7} and
Morgan--Szab\'o--Taubes \cite{MST} which concerns the special case when
$Y\approx S^1\times S^2$ is the link of an embedded $2$--sphere of
self-intersection $0$. One can derive Theorem~\ref{sw-thm1} from Nicolaescu's
proof \cite{nico1} of their result and the classification of closed orientable
$3$--manifolds admitting positive scalar curvature metrics (see
\cite[p\,325]{LM}). However, we shall give a direct (and much simpler)
proof where the main idea is to perturb the monopole equations on $Z$
by a suitable $2$--form such that the corresponding perturbed
Chern--Simons--Dirac functional on $Y$ has no critical points. One then
introduces a long neck $[-T,T]\times Y$. See Section~\ref{proofs12}
for details.

We now turn to another application, for which we need a little preparation.
For any compact $\spc$ $4$--manifold $Z$ whose boundary is a disjoint union
of rational homology spheres set
\[d(Z)=\frac14\left(c_1(\cll_{Z})^2-\si(Z)\right)+b_1(Z)-b_2^+(Z).\]
Here $\cll_Z$ is the determinant line bundle of the $\spc$ structure, and
$\si(Z)$ the signature of $Z$.
If $Z$ is closed then the moduli space $M_Z$
has expected dimension $d(Z)-b_0(Z)$. 

In \cite{Fr4} we will assign to every $\spc$ rational homology $3$--sphere $Y$
a rational number $h(Y)$. (A preliminary version of this invariant
was introduced in \cite{Fr1}.)
In Section~\ref{proofs12} of the present paper this invariant will be defined
in the case when $Y$ admits a metric with positive scalar curvature.
It satisfies $h(-Y)=-h(Y)$. In particular, $h(S^3)=0$.

\begin{thm}\label{sw-thm2}
Let $Z$ be a closed, connected $\spc$ $4$--manifold, and let
$W\subset Z$ be a compact, connected, codimension~$0$ submanifold whose
boundary is a disjoint union of rational homology spheres $Y_1,\dots,Y_r$, $r\ge1$,
each of which admits a metric of positive scalar curvature. Suppose
$b_2^+(W)>0$ and set $W^c=Z\setminus\text{int}\,W$. Let each $Y_j$
have the orientation and $\spc$ structure inherited from $W$. Then the
following hold:
\begin{enumerate}
\item[\rm(i)]If $2\sum_jh(Y_j)\le-d(W)$ then there exists a metric and
perturbation $2$--form on $Z$ for which $M_Z$ is empty, hence $\BF(Z)=0$.
\item[\rm(ii)]If $b^+_2(Z)>1$ and $2\sum_jh(Y_j)<d(W^c)$ then $\SW(Z)=0$.
\end{enumerate}
\end{thm}

Note that (ii) generalizes the classical theorem (see \cite{das1,nico1})
which says that $\SW(Z)=0$ if
$Z$ is a connected sum where both sides have $b^+_2>0$.

\subsection{The Chern--Simons--Dirac functional}
\label{intro-csd}

Let $Y$ be a closed, connected Riemannian $\spc$ $3$--manifold. We consider
the Seiberg--Witten monopole equations over $\ry$, perturbed by adding
a $2$--form to the curvature part of these equations. This
$2$--form should be the pull-back of a closed form $\eta$ on $Y$. 
Recall from \cite{KM4,MST} that in temporal gauge
these perturbed monopole equations can be 
described as the downward gradient flow equation for a perturbed
Chern--Simons--Dirac functional, which we will denote by $\csd\eta$, or just
$\cd$ when no confusion can arise.

For transversality reasons we will add a further small perturbation
to the monopole equations over $\ry$, similar to those
introduced in \cite[Section~2]{Fr1}. This perturbation depends on a parameter
$\prt$ (see Subsection~\ref{subsec:pert}). When $\prt\neq0$ then the perturbed
monopole equations are no longer of gradient flow type.
Therefore,
$\prt$ has to be kept small in order for
the perturbed equations to retain certain properties
(see Subsection~\ref{subsec:curv-bounds}).

If $S$ is a configuration over $Y$ (ie a spin connection together with a 
section of the spin bundle) and $u\col Y\to\U1$ then
\begin{equation}\label{cduS}
\cd(u(S))-\cd(S)=2\pi\int_Y\ti\eta\wedge[u],
\end{equation}
where $[u]\in H^1(Y)$ is the pull-back by $u$ of the fundamental class of
$\U1$, and
\begin{equation}\label{ti-eta}
\ti\eta=\pi c_1(\cll_Y)-[\eta]\in H^2(Y).
\end{equation}
Here $\cll_Y$ is the determinant line bundle of the $\spc$ structure of $Y$.

Let $\fl_Y$ be the space of (smooth)
monopoles over $Y$ (ie critical points of $\cd$)
modulo all gauge transformations $Y\to\U1$,
and $\ti\fl_Y$ the space of monopoles over $Y$
modulo null-homotopic gauge transformations.

When no statement is made to the contrary,
we will always make the following two assumptions:
\begin{enumerate}
\item[(O1)]$\ti\eta$ is a real multiple of some rational cohomology class.
\item[(O2)]All critical points of $\cd$ are non-degenerate.
\end{enumerate}
The second assumption implies that $\fl_Y$ is a finite set.
This rules out the case when
$\ti\eta=0$ and $b_1(Y)>0$, because if $\ti\eta=0$ then the subspace of
reducible points in $\fl_Y$ is homeomorphic to a $b_1(Y)$--dimensional
torus. If $\ti\eta\neq0$ or $b_1(Y)=0$ then the non-degeneracy condition
can be achieved by perturbing $\eta$ by an exact form
(see Proposition~\ref{frechet}).

For any $\al,\beta\in\ti\fl_Y$ let $M(\al,\beta)$ denote the moduli space
of monopoles over $\ry$ that are asymptotic to $\al$ and $\beta$
at $-\infty$ and $\infty$, respectively.
Set $\check M=M/\R$.
By a {\em broken gradient line} from $\al$ to $\beta$ we mean
a sequence $(\om_1,\dots,\om_k)$ where $k\ge0$ and
$\om_j\in\check M(\al_{j-1},\al_j)$
for some $\al_0,\dots,\al_k\in\ti\fl_Y$ with $\al_0=\al$, $\al_k=\beta$, and 
$\al_{j-1}\neq\al_j$ for each $j$. If $\al=\beta$ then we allow the empty
broken gradient line (with $k=0$).

\subsection{Compactness}
\label{intro:comp-single}

Let $X$ be a $\spc$ Riemannian $4$--manifold with tubular ends $\orp\times Y_j$,
$j=1,\dots,r$, where $r\ge0$ and each $Y_j$ is a closed, connected
Riemannian $\spc$ $3$--manifold. Setting $Y=\cup_jY_j$ this
means that we are given
\begin{itemize}
\item an orientation preserving isometric embedding
$\iota\col\orp\times Y\to X$ such that
\begin{equation}\label{endt-def}
X\endt t=X\setminus\iota((t,\infty)\times Y)
\end{equation}
is compact for any $t\ge0$,
\item an isomorphism between the $\spc$ structure
on $\orp\times Y$
induced from $Y$ and the one inherited from $X$ via the embedding $\iota$.
\end{itemize}
Here $\R_+$ is the set of positive real numbers and $\orp=\R_+\cup\{0\}$.
Usually we will just regard $\orp\times Y$ as
a (closed) submanifold of $X$.

Let $\eta_j$ be a closed $2$--form on $Y_j$ and define $\ti\eta_j\in H^2(Y_j)$
in terms of $\eta_j$ as in \Ref{ti-eta}. We write
$\cd$ instead of $\cd_{\eta_j}$ when no confusion is likely to arise.
We perturb the
curvature part of the monopole equations over $X$ by adding a $2$--form
$\mu$ whose restriction to $\rpy_j$ agrees with the pull-back of $\eta_j$.
In addition we perturb the equations over $\ry_j$ and the corresponding
end of $X$ using a
perturbation parameter $\prt_j$.
If $\vec\al=(\al_1,\dots,\al_r)$ with $\al_j\in\ti\fl_{Y_j}$ let
$M(X;\vec\al)$ denote the moduli space of monopoles over $X$ that
are asymptotic to $\al_j$ over $\rpy_j$.

Let $\lla_1,\dots,\lla_r$ be positive constants.
We consider the following two equivalent conditions
on the $\spc$ manifold $X$ and $\ti\eta_j,\lla_j$:
\begin{enumerate}
\item[(A)]There exists a class $\ti z\in H^2(X;\R)$ such that
$\ti z|_{Y_j}=\lla_j\ti\eta_j$ for $j=1,\dots,r$.
\item[{\mathsurround0pt(${\rm A}'$)}]For configurations $S$ over $X\endt 0$ the sum
$\sum_j\lla_j\cd(S|_{\{0\}\times Y_j})$ depends only on the gauge equivalence
class of $S$.
\end{enumerate}
Note that if $\lla_j=1$ for all $j$ then (A) holds precisely when
there exists a class $z\in H^2(X;\R)$ such that $z|_{Y_j}=[\eta_j]$
for $j=1,\dots,r$.

\begin{thm}\label{thm:sing-comp}
If Condition~(A) is satisfied and each $\prt_j$ has sufficiently small
$C^1$ norm then the following holds.
For $n=1,2,\dots$ let $\om_n\in M(X;\vec\al_n)$, where
$\vec\al_n=(\al_{n,1},\dots,\al_{n,r})$. If
\begin{equation}\label{eqn:energy-bound1}
\inf_n\sum_{j=1}^r\lla_j\cd(\al\nj)>-\infty
\end{equation}
then there exists a subsequence of $\om_n$ which chain-converges
to an $(r+1)$--tuple $(\om,\vec v_1,\dots,\vec v_r)$ where $\om$ is an element
of some moduli space $M(X;\vec\beta)$ and $\vec v_j$ is a broken gradient line
over $\ry_j$ from $\beta_j$ to some $\ga_j\in\ti\fl_{Y_j}$.
Moreover, if
$\om_n$ chain-converges to $(\om,\vec v_1,\dots,\vec v_r)$ then
for sufficiently large $n$ there is a gauge transformation $u_n\col X\to\U1$
which is translationary invariant over the ends and
maps $M(X;\vec\al_n)$ to $M(X;\vec\ga)$
\end{thm}

The assumption~\Ref{eqn:energy-bound1} imposes an ``energy bound'' over
the ends of $X$, as we will show in Subsection~\ref{subsec:single}.
The notion of chain-convergence is defined in
Subsection~\ref{subsec:chain-conv}. The limit, if it exists,
is unique up to gauge equivalence
(see Proposition~\ref{prop:chain-conv} below).

\subsection{Compactness and neck-stretching}
\label{intro:comp-neck}

In this subsection cohomology groups will have real coefficients.

We consider again a $\spc$ Riemannian $4$--manifold $X$ as in the previous
subsection, but we now assume that the ends of $X$ are given by orientation
preserving isometric embeddings
\begin{align*}
\iota'_j\col \orpy'_j\to X&,\qquad j=1,\dots,r',\\
\iota^\pm_j\col \orp\times(\pm Y_j)\to X&,\qquad j=1,\dots,r,
\end{align*}
where $r,r'\ge0$. Here each $Y'_j,Y_j$ should be a closed, connected $\spc$ 
Riemannian $3$--manifold, and as before there should be the appropriate 
identifications of $\spc$ structures.
For every $T=(T_1,\dots,T_r)$ with $T_j>0$ for each $j$, let $\xt$ denote
the manifold obtained from $X$ by  gluing, for $j=1,\dots,r$, the two ends
$\iota^\pm_j(\orpy_j)$ to form a neck $[-T_j,T_j]\times Y_j$.
To be precise, let
$\pxt\subset X$ be the result of deleting from $X$ the sets
$\iota^\pm_j([2T_j,\infty)\times Y_j)$, $j=1,\dots,r$. Set
\[\xt=\pxt/\sim,\]
where we identify 
\[\iota^+_j(t,y)\sim\iota^-_j(2T_j-t,y)\]
for all $(t,y)\in(0,2T_j)\times Y_j$ and
$j=1,\dots,r$. We regard $[-T_j,T_j]\times Y_j$ as a submanifold of $\xt$
by means of the isometric embedding
$(t,y)\mapsto\pi_T\iota^+_j(t+T_j,y)$, where
$\pi_T\col \pxt\to\xt$. Also, we write $\R_+\times(\pm Y_j)$ instead of
$\iota^\pm_j(\rpy_j)$, and similarly for $\R_+\times Y'_j$,
if this is not likely to cause any confusion.

Set $X\gl=\xt$ with $T_j=1$ for all $j$.
The process of constructing $X\gl$ from $X$ (as smooth manifolds)
can be described by the unoriented
graph $\ga$ which has one node for every connected component of $X$ and,
for each $j=1,\dots,r$,
one edge representing the pair of embeddings $\iota^\pm_j$.

A node in an oriented graph is called a {\em source} if it has no
incoming edges. 
If $e$ is any node in $\ga$ let $X_e$ denote the corresponding
component of $X$. Let $Z_e=(X_e)\endt1$ be the corresponding truncated manifold
as in \Ref{endt-def}. Let $\ga\setminus e$ be the graph
obtained from $\ga$ by deleting the node $e$ and all edges of which $e$
is a boundary point. Given an orientation $\ortn$ of $\ga$ let 
$\prtl^-Z_e$ denote the union of all boundary components of $Z_e$ corresponding
to incoming edges of $(\ga,\ortn)$. Let $F_e$ be the kernel of
$H^1(Z_e)\to H^1(\prtl^-Z_e)$, and set
\[\Si(X,\ga,\ortn)=\dim\,H^1(X\gl)-\sum_e\dim\,F_e.\]
It will follow from Lemma~\ref{gk-inj} below that
$\Si(X,\ga,\ortn)\le0$ if each connected component of $\ga$ is simply-connected.

We will now state a condition on $(X,\ga)$ which is recursive with respect to
the number of nodes of $\ga$. 
\begin{enumerate}
\item[(C)]If $\ga$ has more than one node then it should
admit an orientation $\ortn$ such that the following two conditions hold:
\begin{itemize}
\item $\Si(X,\ga,\ortn)=0$,
\item Condition~(C) holds for $(X\setminus X_e,\ga\setminus e)$ for all
sources $e$ of $(\ga,\ortn)$.
\end{itemize}
\end{enumerate}
We are only interested in this condition
when each component of $\ga$ is simply-connected.
If $\ga$ is connected and has exactly two nodes $e_1,e_2$ then (C) holds
if and only if
$H^1(X\gl)\to H^1(Z_{e_j})$ is surjective for at least one value of $j$,
as is easily seen from the Mayer--Vietoris sequence.
See Subsection~\ref{subsec:onC} and the proof of
Proposition~\ref{prop:local-energy-bound2} for more information
about Condition~(C).

Let the Chern--Simons--Dirac functionals on $Y_j,Y'_j$ be defined in terms
of closed $2$--forms $\eta_j,\eta'_j$ respectively.
Let $\ti\eta_j$ and $\ti\eta'_j$ be the corresponding classes
as in \Ref{ti-eta}. Let $\lla_1,\dots,\lla_r$ and $\lla'_1,\dots,\lla'_{r'}$
be positive constants. The following conditions on
$X,\ti\eta_j,\ti\eta'_j,\lla_j,\lla'_j$
will appear in Theorem~\ref{thm:neck-comp} below.

\begin{enumerate}
\item[(B1)]There exists a class in $H^2(X\gl)$ whose 
restrictions to $Y_j$ and $Y'_j$ are $[\eta_j]$ and $[\eta'_j]$, respectively,
and all the constants $\lla_j,\lla'_j$ are equal to $1$.
\item[(B2)]There exists a class in $H^2(X\gl)$ whose 
restrictions to $Y_j$ and $Y'_j$ are $\lla_j\ti\eta_j$ and $\lla'_j\ti\eta'_j$,
respectively. 
Moreover, the graph $\ga$
is simply-connected, and Condition~(C) holds for $(X,\ga)$.
\end{enumerate}

Choose a $2$--form $\mu$ on $X$ whose restriction to each end
$\R_+\times(\pm Y_j)$
is the pull-back of $\eta_j$, and whose restriction to
$\rpy'_j$ is the pull-back of $\eta'_j$.
Such a form $\mu$ gives rise,
in a canonical way, to a form $\mu^{(T)}$ on $\xt$. We use the forms
$\mu,\mu^{(T)}$ to perturb the curvature part of the monopole equations
over $X$, $\xt$, respectively.
We use the perturbation parameter $\prt'_j$ over $\ry'_j$ and the corresponding
ends, and $\prt_j$ over $\ry_j$ and the corresponding ends and necks.

Moduli spaces over $X$ will be denoted $M(X;\vec\al_+,\vec\al_-,\vec\al')$,
where the $j$'th component of $\vec\al_\pm$ specifies the limit over the
end $\R_+\times(\pm Y_j)$ and the
$j$'th component of $\vec\al'$ specifies the limit over
$\rpy'_j$. 

\begin{thm}\label{thm:neck-comp}
Suppose at least one of the conditions (B1), (B2) holds, and
for $n=1,2,\dots$ let $\om_n\in M(\xtn;\vec\al'_n)$, where
$\vec\al'_n=(\al'_{n,1},\dots,\al'_{n,r'})$ and
$T_j(n)\to\infty$ for $j=1,\dots,r$. Suppose also that
the perturbation parameters $\prt_j,\prt'_j$ are admissible
for each $\vec\al'_n$, and that
\[\inf_n\sum_{j=1}^{r'}\lla'_j\cd(\al'\nj)>-\infty.\]
Then there exists a subsequence of $\om_n$ which
chain-converges to an $(r+r'+1)$--tuple
$\bv=(\om,\vec v_1,\dots,\vec v_r,
\vec v'_1,\dots,\vec v'_{r'})$, where
\begin{itemize}
\item $\om$ is an element of some moduli space
$M(X;\vec\al_1,\vec\al_2,\vec\beta')$,
\item $\vec v_j$ is a broken gradient line over $\ry_j$ from $\al_{1j}$
to $\al_{2j}$,
\item $\vec v'_j$ is a broken gradient line over $\ry'_j$ from $\beta'_j$
to some $\ga'_j\in\ti\fl_{Y'_j}$.
\end{itemize}
Moreover, if $\om_n$ chain-converges to $\bv$ then
for sufficiently large $n$ there is a gauge transformation
$u_n\col \xtn\to\U1$
which is translationary invariant over the ends and
maps $M(\xtn;\vec\al'_n)$ to $M(\xtn;\vec\ga')$.
\end{thm}

The notion of chain-convergence is defined in
Subsection~\ref{subsec:chain-conv}. Note that the chain-limit is 
unique only up to gauge equivalence, see Proposition~\ref{prop:chain-conv}.

What it means for the perturbation parameters
$\prt_j,\prt'_j$ to be ``admissible'' is
defined in Definition~\ref{defn:suff-small}. As in Theorem~\ref{thm:sing-comp},
if (B2) holds and the 
perturbation parameters have sufficiently small $C^1$ norm then
they are admissible
for any $\vec\al'$, see Proposition~\ref{prop:local-energy-bound2}. 
If (B1) is satisfied but perhaps not (B2)
then for any $C_1<\infty$ there is a $C_2>0$ such that if the 
perturbation parameters have $C^1$ norm $<C_2$ then they
are admissible for all $\vec\al'$ 
satisfying $\sum_{j=1}^{r'}\lla'_j\cd(\al'_j)>-C_1$,
see the remarks after Proposition~\ref{prop:facd-bound}.

The conditions (B1), (B2) in the theorem correspond to the two approaches
to compactness referred to at the beginning of this introduction:
If  (B1) is satisfied then one can take
the ``energy approach'', whereas if (B2) holds one can use the
``Hodge theory approach''

The conclusion of the theorem does not hold in general when neither
(B1) nor (B2) are satisfied. For in that case Theorem~\ref{sw-thm1} would
hold if instead of (ii) one merely assumed that $b_1(Y)>0$. Since $\R^4$
contains an embedded $S^1\times S^2$ this would contradict the fact that
there are many $\spc$ $4$--manifolds with $b^+_2>1$ and non-zero 
Seiberg--Witten invariant.

For the moment we will abuse language and say that (B2) holds
if it holds for some choice of constants $\lla_j,\lla'_j$, and similarly
for (B1). Then
a simple example where (B1) is satisfied but not (B2) is $X=\ry$, where one
glues the two ends to obtain $\xt=(\R/2T\z)\times Y$. There are also many
examples where (B2) is satisfied but not (B1). For instance, consider the case
when $X$ consists of two copies of $\ry$, say
$X=\ry\times\{1,2\}$ with $Y$ connected, and
one glues $\R_+\times Y\times\{1\}$ with $\R_-\times Y\times\{2\}$.
In this case $r=1$ and $r'=2$, so we are given closed $2$--forms
$\eta_1,\eta'_1,\eta'_2$ on $Y$. Condition~(B1) now
requires that these three $2$--forms
represent the same cohomology class, while (B2) holds as long as
there are $a_1,a_2>0$ such that $[\eta_1]=a_1[\eta'_1]=a_2[\eta'_2]$.

\subsection{Outline}

Here is an outline of the content of the remainder of this paper.
In Section~\ref{sec:orbit} we study orbit spaces of configurations over
Riemannian $n$--manifolds with tubular ends. This includes results on the
Laplacians on such manifolds, which are also applied later in
Subsection~\ref{subsec:hodge} to the study of the $d^*+d^+$ operator in the
case $n=4$. Section~\ref{sec:moduli} introduces monopoles, perturbations,
and moduli spaces. Section~\ref{sec:local-comp1}
establishes local compactness results for monopoles,
and technical results on perturbations. In Section~\ref{sec:local-comp2} 
the Hodge theory approach to local compactness is presented, which is an
alternative to the energy approach of Subsection~\ref{subsec:neckstr1}.
Section~\ref{sec:exp-decay} is devoted to exponential decay, which is needed
for the global compactness results in Section~\ref{sec:global-comp},
where Theorems~\ref{thm:sing-comp} and \ref{thm:neck-comp} are proved.
Section~\ref{section:transv} discusses non-degeneracy of critical points
and regularity of moduli spaces. Finally, in Section~\ref{proofs12} we
prove Theorems~\ref{sw-thm1} and \ref{sw-thm2}. There are also two appendices.
The first of these explains how to patch together sequences of local
gauge transformations (an improvement of Uhlenbeck's technique). The second
appendix contains a quantitative inverse function theorem which is
used in the proof of exponential decay.

\section{Configuration spaces}\label{sec:orbit}

\subsection{Configurations and gauge transformations}
\label{conf-gt}

Let $X$ be a Riemannian $n$--manifold with tubular ends $\orp\times Y_j$,
$j=1,\dots,r$, where $n\ge1$, $r\ge0$, and each $Y_j$ is a closed, connected
Riemannian $(n-1)$--manifold. 
This means that we are given for each $j$ an isometric embedding
\[\iota_j\col \orp\times Y_j\to X;\]
moreover, the images of these embeddings are disjoint and their union
have precompact complement. Usually we will just regard $\orp\times Y_j$ as
a submanifold of $X$. Set $Y=\cup_jY_j$ and, for $t\ge0$,
\[X\endt t=X\setminus(t,\infty)\times Y.\]
Let $\bs\to X$ and $\bs_j\to Y_j$ be Hermitian complex vector bundles,
and $L\to X$ and $L_j\to Y_j$ principal $\U1$--bundles.
Suppose we are given, for each $j$, isomorphisms
\[\iota_j^*\bs\oset\approx\to\orp\times\bs_j,\quad
\iota^*_jL\oset\approx\to\orp\times L_j.\]
By a {\em configuration} in $(L,\bs)$ we shall mean a pair $(A,\Phi)$
where $A$ is a connection in $L$ and $\Phi$ a section of $\bs$.
Maps $u\col X\to\U1$ are referred to as {\em gauge transformations} and these
act on configurations in the natural way:
\begin{equation*}
u(A,\Phi)=(u(A),u\Phi).
\end{equation*}
The main goal of this section is to prove a "local slice" theorem for certain
orbit spaces of configurations modulo gauge transformations.

We begin by setting up suitable function spaces. 
For $p\ge1$ and any non-negative integer $m$ let $L^p_m(X)$
be the completion of the space of
compactly supported smooth functions on $X$ with
respect to the norm
\[\|f\|_{m,p}=\|f\|_{L^p_m}=
\left(\sum_{k=0}^m\int_X|\nabla^kf|^p\right)^{1/p}.\]
Here the covariant derivative is computed using
some fixed connection in the tangent bundle $TX$ which is translationary
invariant over each
end. Define the Sobolev space $L^p_k(X;\bs)$ of sections of $\bs$ similarly.

We also need weighted Sobolev spaces. For any smooth function
$w\col X\to\R$ set $L^{p,w}_m(X)=e^{-w}L^p_m(X)$ and
\[\|f\|_{L^{p,w}_k}=\|e^wf\|_{L^p_k}.\]
In practice we require that $w$ have a specific form over the ends, namely
\[w\circ\iota_j(t,y)=\si_j t,\]
where the $\si_j$'s are real numbers. 

The following Sobolev embeddings (which hold in $\R^n$, hence over $X$)
will be used repeatedly:
\begin{gather*}
L^p_{m+1}\subset L^{2p}_m\quad\text{if $p\ge n/2$, $m\ge0$,}\\
L^p_2\subset C^0_B\quad\text{if $p>n/2$.}
\end{gather*}
Here $C^0_B$ denotes the Banach space of bounded continuous functions, with
the supremum norm. Moreover, if $pm>n$ then multiplication defines a
continuous map $L^p_m\times L^p_k\to L^p_k$ for $0\le k\le m$.

For the remainder of this section fix $p>n/2$.
Note that this implies $L^p_1\subset L^2$ over compact $n$--manifolds.

We will now define an affine space $\cc$ of $L^p_{1,\loc}$
configurations in $(L,\bs)$.
Let $A\rr$ be a smooth connection in $L$.
Choose a smooth section $\Phi\rr$ of $\bs$ whose restriction to $\rpy_j$
is the pull-back of a section $\psi_j$ of $\bs_j$.
Suppose $\psi_j=0$ for $j\le r_0$, and $\psi_j\not\equiv0$
for $j>r_0$, where $r_0$ is a non-negative integer..
Fix a weight function $w$ as above with
$\si_j\ge0$ small for all $j$, and 
$\si_j>0$ for $j>r_0$. Set
\[\cc=\{(A\rr+a,\Phi\rr+\phi)\st a,\phi\in\llw pw1\}.\]
We topologize $\cc$ using the $\llw pw1$ metric.

We wish to define a Banach Lie group $\cg$ of $L^p_{2,\loc}$ gauge 
transformations over $X$ such that $\cg$ acts smoothly on $\cc$ and such 
that if $S,S'\in\cc$ and $u(S)=S'$ for some $L^p_{2,\loc}$ gauge 
transformations $u$ then $u\in\cg$. If $u\in\cg$ then we must certainly have
\[(-u\inv du,(u-1)\Phi\rr)=u(A\rr,\Phi\rr)-(A\rr,\Phi\rr)\in\llw pw1.\]
Now
\begin{equation}\label{duu}
\|du\|_{\llw pw1}\le\const(\|u\inv du\|_{\llw pw1}+
\|u\inv du\|_{\llw pw1}^2),
\end{equation}
and vice versa, $\|du\|_{\llw pw1}$ controls $\|u\inv du\|_{\llw pw1}$,
so we try 
\[\cg=\{u\in L^p_{2,\loc}(X;\U1)\st du,(u-1)\Phi\rr\in\llw pw1\}.\]
By $L^p_{2,\loc}(X;\U1)$ we mean the set of elements of $L^p_{2,\loc}(X;\co)$
that map into $\U1$.
We will see that $\cg$ has a natural smooth structure such that the above
criteria are satisfied.

(This approach to the definition of $\cg$ was inspired by \cite{D5}.)

\subsection{The Banach algebra}

Let $\ti x$ be a finite subset of $X$ which contains at least one
point from every connected component of $X$ where $\Phi\rr$ vanishes
identically.

\begin{defn}
{\rm Set
\[\ce=\{f\in L^p_{2,\loc}(X;\co)\st df,f\Phi\rr\in\llw pw1\},\]
and let $\ce$ have the norm
\[\|f\|_{\ce}=\|df\|_{\llw pw1}+\|f\Phi\rr\|_{\llw pw1}
+\sum_{x\in\ti x}|f(x)|.\]}
\end{defn}

We will see in a moment that $\ce$ is a Banach algebra
(without unit if $r_0<r$).
The next lemma shows that the topology on $\ce$ is independent of the choice
of $\Phi\rr$ and $\ti x$. 

\begin{lemma}\label{Z-ineq}
Let $Z\subset X$ be a compact, connected codimension~$0$ submanifold.
\begin{enumerate}
\item[\rm(i)]If $\Phi\rr|_Z\not\equiv0$ then there is a constant $C$ such that
\begin{equation}\label{zineq}
\int_Z|f|^p\le C\int_Z|df|^p+|f\Phi\rr|^p
\end{equation}
for all $f\in L^p_1(Z)$.
\item[\rm(ii)]There are constants $C_1,C_2$ such that 
\[|f(z_2)-f(z_1)|\le C_1\|df\|_{L^{2p}(Z)}
\le C_2\|df\|_{L^p_1(Z)}\]
for all $f\in L^p_2(Z)$ and $z_1,z_2\in Z$.
\end{enumerate}
\end{lemma}

\begin{proof} Part~(i) follows from the compactness of the embedding
$L^p_1(Z)\to L^p(Z)$. The first inequality in (ii) can either be
deduced from the compactness of $L^{2p}_1(Z)\to C^0(Z)$, or
one can prove it directly, as a step towards proving the Rellich lemma,
by considering the integrals of $df$ along a suitable family
of paths from $z_1$ to $z_2$.
\end{proof}

\begin{lemma}\label{Y-ineq}
Let $Y$ be a closed Riemannian manifold, and $\si>0$.
\begin{enumerate}
\item[\rm(i)]If $q\ge1$ and $f\col \rpy\to\R$ is a $C^1$ function
such that $\lim_{t\to\infty}f(t,y)=0$ for all $y\in Y$ then
\[\|f\|_{L^{q,\si}(\rpy)}\le\si\inv\|\prtl_1f\|_{L^{q,\si}(\rpy)}.\]
\item[\rm(ii)]If $q>1$, $T\ge1$, and $f\col [0,T]\times Y\to\R$
is a $C^1$ function then
\[\|f\|_{L^q([T-1,T]\times Y)}\le\|f_0\|_{L^q(Y)}+
(\si r)^{-1/r}\|\prtl_1f\|_{L^{q,\si}([0,T]\times Y)},\]
where $f_0(y)=f(0,y)$ and $\frac1q+\frac1r=1$.
\end{enumerate}
\end{lemma}
Here $\prtl_1$ is the partial derivative in the first variable, ie in the
$\R_+$ coordinate.

\begin{proof} Part~(i):
\begin{align*}
\|f\|_{L^{q,\si}(\rpy)}
&=\left(\int_{\rpy}\left|\int_0^\infty e^{\si t}\prtl_1f(s+t,y)\,ds
\right|^qdt\,dy\right)^{1/q}\\
&\le\int_0^\infty\left(\int_{\rpy}|e^{\si t}\prtl_1f(s+t,y)|^qdt\,dy\right)
^{1/q}ds\\
&\le\left(\int_0^\infty e^{-\si s}ds\right)
\left(\int_{\rpy}|e^{\si(s+t)}\prtl_1f(s+t,y)|^qdt\,dy\right)^{1/q}\\
&\le\si\inv\|\prtl_1f\|_{L^{q,\si}(\rpy)}.
\end{align*}
Part~(ii) follows by a similar computation.
\end{proof}

Parts~(i)--(iv) of the following proposition are essentially due to 
Donaldson~\cite{D5}.

\begin{prop}\label{banach-alg}$\phantom{9}$
\begin{enumerate}
\item[\rm(i)]There is a constant $C_1$ such that
\[\|f\|_\infty\le C_1\|f\|_{\ce}\]
for all $f\in \ce$.
\item[\rm(ii)]For every $f\in \ce$ and $j=1,\dots,r$ the restriction
$f|_{\{t\}\times Y_j}$ converges
uniformly to a constant function $f^{(j)}$ as $t\to\infty$, and
$f^{(j)}=0$ for $j>r_0$.
\item[\rm(iii)]There is a constant $C_2$ such that if $f\in \ce$ and $f^{(j)}=0$
for all $j$ then
\[\|f\|_{\llw pw2}\le C_2\|f\|_{\ce}.\]
\item[\rm(iv)]There is an exact sequence
\[0\to\llw pw2\oset\iota\to\ce\oset e\to\co^{r_0}\to0,\]
where $\iota$ is the inclusion and $e(f)=(f^{(1)},\dots,f^{(r_0)})$.
\item[\rm(v)]$\ce$ is complete, and multiplication defines a continuous map
$\ce\times\ce\to\ce$.
\end{enumerate}
\end{prop}

\begin{proof}
First observe that for any $f\in L^p_{2,\loc}(X)$ and $\eps>0$ there exists
a $g\in C^\infty(X)$ such that $\|g-f\|_{\llw pw2}<\eps$. Therefore
it suffices to prove (i)--(iii) when $f\in \ce$ is smooth.
Part~(i) is then a consequence of Lemma~\ref{Z-ineq} and
Lemma~\ref{Y-ineq}~(ii), while Part~(ii) for $r_0<j\le r$ follows from 
Lemma~\ref{Z-ineq}.

We will now prove (ii) when $1\le j\le r_0$. Let $f\in \ce$ be smooth.
Since $\int_{\rpy_j}|df|<\infty$
by the H\"older inequality, we have
\[\int_0^\infty|\prtl_1f(t,y)|\,dt<\infty\qquad\text{for almost all\ $y\in Y_j$}.\]
For $n\in\n$ set $f_n=f|_{[n-1,n+1]\times Y_j}$,
regarded as a function on $B=[n-1,n+1]\times Y_j$. Then 
$\{f_n\}$ converges a.e., so by Egoroff's
theorem $\{f_n\}$ converges uniformly over some subset
$T\subset B$ of positive measure.
There is then a constant $C>0$, depending on $T$, such that for
every $g\in L^p_1(B)$ one has
\[\int_B|g|^p\le C(\int_B|dg|^p+\int_T|g|^p).\]
It follows that $\{f_n\}$ converges in $L^p_2$ over $B$,
hence uniformly over $B$,
to some constant function.

Part~(iii) follows from Lemma~\ref{Z-ineq} and Lemma~\ref{Y-ineq}~(i).
Part~(iv) is an immediate consequence of (ii) and (iii). 
It is clear from (i) that $\ce$ is complete. The multiplication property
follows easily from (i) and the fact that smooth functions are dense in $\ce$.
\end{proof}

\subsection{The infinitesimal action}

If $f\col X\to i\R$ and $\Phi$ is a section of $\bs$ we
define a section of $i\La^1\oplus\bs$ by
\[\ci_\Phi f=(-df,f\Phi)\]
whenever the expression on the right makes sense. Here $\La^k$ denotes the
bundle of $k$--forms (on $X$, in this case). If $S=(A,\Phi)$ is a configuration
then we will sometimes write $\ci_S$ instead of $\ci_\Phi$. Set
\[\ci=\ci_{\Phi\rr}.\]
If $\Phi$ is smooth then the formal adjoint of the operator $\ci_\Phi$ is
\[\ci_\Phi^*(a,\phi)=-d^*a+i\la i\Phi,\phi\ra_{\R},\]
where $\la\sdot,\sdot\ra_{\R}$ is the real inner product on $\bs$.
Note that
\[\ci_\Phi^*\ci_\Phi=\Delta+|\Phi|^2\]
where $\Delta$ is the positive Laplacian on $X$.

Set
\[L\cg=\{f\in\ce\st\text{$f$ maps into $i\R$}\}.\]
 From Proposition~\ref{banach-alg}~(i)
we see that the operators
\[\ci_\Phi\col L\cg\to\llw pw1,\quad
\ci_\Phi^*\col \llw pw1\to\lw pw\]
are well-defined and bounded for every
$\Phi\in\Phi\rr+\llw pw1(X;\bs)$.

\begin{lemma}\label{laplace-inj}
For every $\Phi\in\Phi\rr+\llw pw1(X;\bs)$,
the operators $\ci_\Phi^*\ci_\Phi$ and $\ci_\Phi$ have the same kernel
in $L\cg$.
\end{lemma}

\begin{proof} Choose a smooth function $\beta\col \R\to\R$ such that
$\beta(t)=1$ for $t\le1$, $\beta(t)=0$ for $t\ge2$.
For $r>0$ define a compactly supported function
\[\beta_r\col X\to\R\]
by $\beta_r|_{X\endt 0}=1$, and $\beta_r(t,y)=\beta(t/r)$ for
$(t,y)\in\rpy_j$.

Now suppose $f\in L\cg$ and $\ci^*_\Phi\ci_\Phi f=0$. 
Lemma~\ref{banach-alg}~(i) and elliptic regularity gives
$f\in L^p_{2,\loc}$, so we certainly have
\[\ci_\Phi f\in L^p_{1,\loc}\subset L^2_{\loc}.\]
Clearly,
\[\|\ci_\Phi f\|_2\le\liminf_{r\to\infty}\|\ci_\Phi(\beta_rf)\|_2.\]
Over $\rpy_j$ we have
\[\ci^*_\Phi\ci_\Phi=-\prtl_1^2+\Delta_{Y_j}+|\Phi|^2,\]
where $\prtl_1=\frac\prtl{\prtl t}$ and
$\Delta_{Y_j}$ is the positive Laplacian on $Y_j$, so
\begin{align*}
\|\ci_\Phi(\beta_rf)\|^2_2
&=\int_X\ci^*_\Phi\ci_\Phi(\beta_rf)\cdot\beta_r\bar f\\
&=-\sum_j\int_{\rpy_j}((\prtl_1^2\beta_r)f+2(\prtl_1\beta_r)(\prtl_1 f))\cdot
\beta_r\bar f\\
&\le C_1\|f\|^2_\infty\int_0^\infty|r^{-2}\beta''(t/r)|\,dt \\
&\qquad+ C_1\|f\|_\infty\|df\|_p
\left(\int_0^\infty|r\inv\beta'(t/r)|^qdt\right)^{1/q}\\
&\le C_2\|f\|^2_{\ce}\left(r\inv\int_1^2|\beta''(u)|\,du+
r^{1-q}\int_1^2|\beta'(u)|^qdu\right)\\
&\text{$\to0$ as $r\to\infty$},
\end{align*}
where $C_1,C_2>0$ are constants and
$\frac1p+\frac1q=1$. Hence $\ci_\Phi f=0$.
\end{proof}

\begin{lemma}\label{ind-r}
$\ci^*\ci\col \llw qw2(X)\to\lw qw(X)$ is Fredholm of index $-r_0$, for
$1<q<\infty$.
\end{lemma}

\begin{proof} Because $\ci^*\ci$ is elliptic the operator in the
Lemma is Fredholm if the operator
\begin{equation}\label{Yjlapl}
-\prtl_1^2+\Delta_{Y_j}+|\psi_j|^2\col \llw q{\si_j}2\to
\lw q{\si_j},
\end{equation}
acting on functions on $\ry_j$, is Fredholm for each $j$.
The proof of \cite[Proposition~3.21]{D5} (see also \cite{LMcO})
can be generalized to show that \Ref{Yjlapl} is Fredholm if $\si_j^2$ is
not an eigenvalue of $\Delta_{Y_j}+|\psi_j|^2$. Since we are taking 
$\si_j\ge0$ small, and $\si_j>0$ if $\psi_j=0$, this establishes the
Fredholm property in the Lemma.

We will now compute the index. Set
\[\text{ind}^\pm=
\text{index}\{\ci^*\ci\col \llw q{\pm w}2(X)\to\lw q{\pm w}(X)\}.\]
Expressing functions on $Y_j$ in terms of eigenvectors of 
$\Delta_{Y_j}+|\psi_j|^2$ as in \cite{APS1,D5}
one finds that the kernel of
$\ci^*\ci$ in $\llw {{q'}}{\pm w}{{k'}}$
is the same for all $q'>1$ and integers $k'\ge0$. 
Combining this with the
fact that $\ci^*\ci$ is formally self-adjoint we see that
\[\text{ind}^+=-\text{ind}^-.\]
Now choose smooth functions
$w_j\col \R\to\R$ such that $w_j(t)=\si_j|t|$ for $|t|\ge1$. We will apply the 
addition formula for the index, see \cite[Proposition~3.9]{D5}.
This is proved only
for first order operators in \cite{D5} but holds for higher order operators
as well, with essentially the same proof.
The addition formula gives
\begin{align*}
\ind^-&=\ind^++\sum_j\text{index}\{\ci^*_{\psi_j}\ci_{\psi_j}\col
\llw q{-w_j}2(\ry_j)\to\lw q{-w_j}(\ry_j)\}\\
&=\ind^++2\sum_j\dim\ker(\Delta_{Y_j}+|\psi_j|^2)\\
&=\ind^++2r_0.
\end{align*}
Therefore, $\ind^+=-r_0$ as claimed.
\end{proof}

\begin{prop}\label{laplace-iso}
For any $\Phi\in\Phi\rr+\llw pw1(X;\bs)$ the following hold:
\begin{enumerate}
\item[\rm(i)]The operator
\begin{equation}\label{cici}
\ci_\Phi^*\ci_\Phi\col L\cg\to\lw pw
\end{equation}
is Fredholm of index~$0$, and it has the same kernel as
$\ci_\Phi\col L\cg\to\llw pw1$ and the same image as
$\ci_\Phi^*\col \llw pw1\to\lw pw$.
\item[\rm(ii)]$\ci_\Phi(L\cg)$ is closed in $\llw pw1$ and
\begin{equation}\label{vplus}
\llw pw1(i\La^1\oplus\bs)=\ci_\Phi(L\cg)\oplus\ker(\ci_\Phi^*).
\end{equation}
\end{enumerate}
\end{prop}

\begin{proof} It is easy to deduce Part~(ii) from Part~(i).
We will now prove Part~(i). Since
\[\ci^*_\Phi\ci_\Phi-\ci^*\ci=|\Phi|^2-|\Phi\rr|^2\col \llw pw2\to\lw pw\]
is a compact operator,
$\ci^*_\Phi\ci_\Phi$ and
$\ci^*\ci$ have the same index as operators
between these Banach spaces. It then follows from Lemma~\ref{ind-r} and
Proposition~\ref{banach-alg}~(iv) that the operator
\Ref{cici} is Fredholm of
index~$0$. The statement about the kernels is the same as
Lemma~\ref{laplace-inj}. To prove the statement about the images, we may
as well assume $X$ is connected. If $\Phi\neq0$ then the operator~\Ref{cici}
is surjective and there is nothing left to prove.
Now suppose $\Phi=0$. Then all the
weights $\si_j$ are positive, and the kernel of $\ci_\Phi$ in $\ce$
consists of the constant functions. Hence the image of \Ref{cici}
has codimension~$1$. But $\int_Xd^*a=0$ for every
$1$--form $a\in\llw pw1$, so $d^*\col \llw pw1\to\lw pw$ is not surjective.
\end{proof}

In the course of the proof of (i) we obtained:

\begin{prop}\label{dd-image}
If $X$ is connected and $\si_j>0$ for all $j$ then 
\[d^*d(\ce)=\left\{g\in L^{p,w}(X;\co)\st\int_Xg=0\right\}.\]
\end{prop}

We conclude this subsection with a result that will be
needed in the proofs of Proposition~\ref{prop:loc-comp} and Lemma~\ref{econv}
below.
Let $1<q<\infty$ and for any $L^q_{1,\loc}$ function $f\col X\to\R$ set
\[\del_jf=\int_{\{0\}\times Y_j}\prtl_1f,\quad j=1,\dots,r.\]
The integral is well-defined because if $n$ is any positive integer
then there is a bounded restriction map
$L^q_1(\R^n)\to L^q(\{0\}\times\R^{n-1})$.

Choose a point $x_0\in X$.

\begin{prop}\label{laplace-iso2}
If $X$ is connected, $1<q<\infty$, $r\ge1$, and if
$\si_j>0$ is sufficiently small for each $j$ then the operator
\begin{align*}
\beta\col L^{q,-w}_2(X;\R)&\to L^{q,-w}(X;\R)\oplus\R^r,\\
f&\mapsto(\Delta f,(\del_1f,\dots,\del_{r-1}f,f(x_0)))
\end{align*}
is an isomorphism.
\end{prop}

\begin{proof} By the proof of Lemma~\ref{ind-r},
$\Delta\col L^{q,-w}_2\to L^{q,-w}$ has index 
$r$, hence $\ind(\beta)=0$. We will show $\beta$ is injective.
First observe that $\sum_{j=1}^r\del_jf=0$ whenever
$\Del f=0$, so if $\beta f=0$ then $\del_jf=0$ for all $j$.

Suppose $\beta f=0$. 
To simplify notation we will now assume $Y$ is connected.
Over $\rpy$ we have $\Del=-\prtl_1^2+\Del_Y$.
Let $\{h_\nu\}_{\nu=0,1,\dots}$ be a maximal orthonormal set of
eigenvectors of $\Del_Y$,
with corresponding eigenvalues $\lla^2_\nu$, where
$0=\lla_0<\lla_1\le\lla_2\le\cdots$. Then 
\[f(t,y)=a+bt+g(t,y),\]
where $a,b\in\R$, and $g$ has the form
\[g(t,y)=\sum_{\nu\ge1} c_\nu e^{-\lla_\nu t}h_\nu(y)\]
for some real constants $c_\nu$.
Elliptic estimates show that $g$ decays exponentially,
or more precisely,
\[|(\nabla^jf)_{(t,y)}|\le d_je^{-\lla_1t}\]
for $(t,y)\in\rpy$ and $j\ge0$, where $d_j>0$ is a constant.
Now
\[\prtl_1f(t,y)=b-\sum_{\nu\ge1} c_\nu\lla_\nu e^{-\lla_\nu t}h_\nu(y).\]
Since $\Del_Y$ is formally self-adjoint we have
$\int_Yh_\nu=0$ if $\lla_\nu\neq0$,
hence
\[b\,\text{Vol}(Y)=\int_{\{\tau\}\times Y}\prtl_1f=0.\]
It follows that $f$ is bounded and $df$ decays exponentially over $\rpy$,
so
\[0=\int_Xf\Del f=\int_X|df|^2,\]
hence $f$ is constant. Since $f(x_0)=0$ we have $f=0$.
\end{proof}

\subsection{Local slices}\label{local-slices}

Fix a finite subset $\hx\subset X$. 

\begin{defn}
{\rm Set
\begin{align*}
&\cg_\hx=\{u\in1+\ce\st\text{$u$ maps into $\U1$ and $u|_\hx\equiv1$}\}\\
&L\cg_\hx=\{f\in\ce\st\text{$f$ maps into $i\R$ and $f|_\hx\equiv0$}\}
\end{align*}
and let $\cg_\hx$ and $L\cg_\hx$ have the subspace
topologies inherited from $1+\ce\approx\ce$ and $\ce$, respectively.}
\end{defn}
By $1+\ce$ we mean the set of functions on $X$ of the form $1+f$ where
$f\in\ce$.
If $\hx$ is empty then we write $\cg$ instead of $\cg_\hx$, and similarly 
for $L\cg$.

\begin{prop}$\phantom{9}$
\begin{enumerate}
\item[\rm(i)]$\cg_\hx$ is a smooth submanifold of $1+\ce$ and a
Banach Lie group with Lie algebra $L\cg_\hx$.
\item[\rm(ii)]The natural action $\cg_\hx\times\cc\to\cc$ is smooth.
\item[\rm(iii)]If $S\in\cc$, $u\in L^p_{2,\loc}(X;\U1)$
and $u(S)\in\cc$ then $u\in\cg$.
\end{enumerate}
\end{prop}

\begin{proof} (i)\qua If $r_0<r$ then $1\not\in\ce$, but in any case,
\[f\mapsto\sum_{k=1}^\infty\frac1{k!}f^k=\exp(f)-1\]
defines a smooth map $\ce\to\ce$, by Proposition~\ref{banach-alg}~(v).
Therefore, the exponential map provides the local parametrization around $1$
required for $\cg_\hx$ to be a submanifold of $1+\ce$.
The verification of (ii) and (iii) is left to the reader.
\end{proof}

Let $\cb_\hx=\cc/\cg_\hx$ have the quotient topology. This topology is
Hausdorff because it is stronger than the topology defined by the
$L^{2p}$ metric on $\cb_\hx$ (see \cite{DK}). 
The image in $\cb_\hx$ of a configuration $S\in\cc$ will be denoted $[S]$,
and we say $S$ is a {\em representative} of $[S]$.

Let $\cc^*_\hx$ be the set of all elements of $\cc$ which
have trivial stabilizer in $\cg_\hx$. In other words, $\cc^*_\hx$ consists
of those $(A,\Phi)\in\cc$ such that $\hx$ contains at least one point
 from every component of $X$ where $\Phi$ vanishes almost everywhere.
Let $\cb^*_\hx$ be the image of $\cc^*_\hx\to\cb_\hx$. It is clear that
$\cb^*_\hx$ is an open subset of $\cb_\hx$.

If $\hx$ is empty then $\cc^*\subset\cc$ and $\cb^*\subset\cb$
are the subspaces of
{\em irreducible} configurations. As usual, a configuration that is
not irreducible is called reducible.

We will now
give $\cb^*_\hx$ the structure of a smooth Banach manifold by specifying an
atlas of local parametrizations. Let $S=(A,\Phi)\in\cc^*_\hx$ and set
\[V=\ci^*_\Phi(\llw pw1),\quad W=\ci^*_\Phi\ci_\Phi(L\cg_\hx).\]
By Proposition~\ref{laplace-iso} we have
\[\dim(V/W)=|\hx|-\ell\]
where $\ell$ is the number of components of $X$ where $\Phi$ vanishes
almost everywhere
Choose a bounded linear map $\rho\col V\to W$ such that $\rho|_W=I$, and set
\[\ci^\#_\Phi=\rho\ci^*_\Phi.\]
Then
\[\llw pw1(i\La^1\oplus\bs)=\ci_\Phi(L\cg_\hx)\oplus\ker(\ci^\#_\Phi)\]
by Proposition~\ref{laplace-iso}. Consider the smooth map
\[\Pi\col L\cg_\hx\times\ker(\ci^\#_\Phi)\to\cc,\quad(f,s)\mapsto\exp(f)(S+s).\]
The derivative of this map at $(0,0)$ is
\[D\Pi(0,0)(f,s)=\ci_\Phi f+s,\]
which is an isomorphism by the above remarks. The inverse function theorem
then says that $\Pi$ is a local diffeomorphism at $(0,0)$. 

\begin{prop}\label{local-slice}
In the situation above there is an open
neighbourhood $U$ of $0\in\ker(\ci^\#_\Phi)$ such that
the projection $\cc\to\cb_\hx$ restricts to a topological embedding
of $S+U$ onto an open subset of $\cb^*_\hx$.
\end{prop}

It is clear that the collection of such local parametrizations
$U\to\cb^*_\hx$ is a smooth atlas for $\cb^*_\hx$.

\begin{proof}
It only remains to prove that $S+U\to\cb_\hx$ is injective when
$U$ is sufficiently small. So suppose $(a_k,\phi_k)$, $(b_k,\psi_k)$ are
two sequences in $\ker(\ci^\#_\Phi)$ which both converge to $0$
as $k\to\infty$, and such that
\[u_k(A+a_k,\Phi+\phi_k)=(A+b_k,\Phi+\psi_k)\]
for some $u_k\in\cg_\hx$. We will show that $\|u_k-1\|_{\ce}\to0$.
Since $\Pi$ is a local diffeomorphism at $(0,0)$, this will imply that
$u_k=1$ for $k\gg0$. 

Written out, the assumption on $u_k$ is that
\begin{align*}
u_k\inv du_k&=a_k-b_k,\\
(u_k-1)\Phi&=\psi_k-u_k\phi_k.
\end{align*}
By \Ref{duu} we have $\|du_k\|_{\llw pw1}\to0$, which in turn gives
$\|u_k\phi_k\|_{\llw pw1}\to0$,  hence
\begin{equation}\label{unphi}
\|(u_k-1)\Phi\|_{\llw pw1}\to0.
\end{equation}
Because $u_k$ is bounded and $du_k$ converges to $0$ in $L^p_1$ over
compact subsets, we can find a subsequence $\{k_j\}$ such that
$u_{k_j}$ converges in $L^p_2$ over compact subsets to a locally constant
function $u$. Then $u|_\hx=1$ and $u\Phi=\Phi$, hence $u=1$.
Set $f_j=u_{k_j}-1$ and $\phi=\Phi-\Phi\rr\in\llw pw1$. Then
$\|df_j\otimes\phi\|_{\lw pw}\to0$. Furthermore, given $\eps>0$ we can find
$t>0$ such that 
\[\int_{[t,\infty)\times Y}|e^w\phi|^p<\frac\eps4\]
and $N$ such that
\[\int_{X\endt t}|e^wf_j\phi|^p<\frac\eps2\]
for $j>N$.
Then $\int_X|e^wf_j\phi|^p<\eps$ for $j>N$. Thus $\|f_j\phi\|_{\lw pw}\to0$,
and similarly $\|f_j\nabla\phi\|_{\lw pw}\to0$. Altogether this shows that
$\|f_j\phi\|_{\llw pw1}\to0$. Combined  with \Ref{unphi} this yields
\[\|(u_{k_j}-1)\Phi_0\|_{\llw pw1}\to0,\]
hence $\|u_{k_j}-1\|_\ce\to0$. But we can run the above argument starting with
any subsequence of $\{u_k\}$, so $\|u_k-1\|_\ce\to0$.
\end{proof}

\subsection{Manifolds with boundary}
\label{boundary-slice}

Let $Z$ be a compact, connected, oriented Riemannian $n$--manifold,
perhaps with boundary,
and $\hx\subset Z$ a finite subset.
Let $\bs\to Z$ be a Hermitian vector bundle and $L\to Z$ a principal
$\U1$--bundle. Fix $p>n/2$ and let $\cc$ denote the space of
$L^p_1$ configurations $(A,\Phi)$ in $(L,\bs)$.
Let $\cg_\hx$ be the group of those
$L^p_2$ gauge transformations $Z\to\U1$
that restrict to $1$ on $\hx$, and $\cc^*_\hx$ the set of all
elements of $\cc$ that have trivial stabilizer in $\cg_\hx$.
Then
\[\cb^*_\hx=\cc^*_\hx/\cg_\hx\]
is again a (Hausdorff) smooth Banach manifold. As for orbit spaces of
connections (see \cite[p\,192]{DK}) the main ingredient here is the
solution to the Neumann problem over $Z$, according to which the operator
\begin{align*}
T_\Phi\col L^p_2(Z)&\to L^p(Z)\oplus\prtl L^p_1(\prtl Z),\\
f&\mapsto(\Delta f+|\Phi|^2f,\prtl_\nu f)
\end{align*}
is a Fredholm operator of index $0$ (see \cite[Section~5.7]{Taylor1} and
\cite[pages 85--86]{Hamilton1}). Here $\nu$ is the inward-pointing unit normal
along $\prtl Z$, and $\prtl L^p_1(\prtl Z)$ is the space of
boundary values of $L^p_1$ functions on $Z$. Henceforth we work with 
imaginary-valued functions, and on $\prtl Z$ we identify $3$--forms
with functions by means of the Hodge $*$--operator.
Then $T_\Phi=J_\Phi\ci_\Phi$, where
\[J_\Phi(a,\phi)=(\ci^*_\Phi(a,\phi),(*a)|_{\prtl Z}).\]
Choose a bounded linear map
\[\rho\col L^p(Z)\oplus\prtl L^p_1(\prtl Z)\to W:=T_\Phi(L\cg_\hx)\]
which restricts to the identity on $W$, and set $J^\#_\Phi=\rho J_\Phi$.
An application of Stokes' theorem
shows that
\[\ker(T_\Phi)\subset\ker(\ci_\Phi)\quad\text{in $L^p_2(Z)$,}\]
hence
\[T_\Phi=J^\#_\Phi\ci_\Phi\col L\cg_\hx\to W\]
is an isomorphism. In general, if $V_1\oset{T_1}\to V_2\oset{T_2}\to V_3$
are linear maps between vector spaces such that $T_2T_1$ is an isomorphism,
then $V_2=\im(T_1)\oplus\ker(T_2)$. Therefore, for any $(A,\Phi)\in\cc^*_\hx$
we have
\[L^p_1(Z;i\La^1\oplus\bs)=\ci_\Phi(L\cg_\hx)\oplus\ker(J^\#_\Phi),\]
where both summands are closed subspaces. Thus we obtain the analogue
of Proposition~\ref{local-slice} with local slices of the form
$(A,\Phi)+U$, where $U$ is a small neighbourhood of $0\in\ker(J^\#_\Phi)$.

\section{Moduli spaces}\label{sec:moduli}

\subsection{$\Spc$ structures}\label{subsec:spc-str}

It will be convenient to have a definition of $\spc$ structure that does not
refer to Riemannian metrics. 
So let $X$ be an oriented $n$--dimensional manifold and $\pglp$
its bundle of positive linear frames. Let $\tglp n$ denote the $2$--fold
universal covering group of the identity component $\glp n$ of $\GL(n,\R)$,
and denote by $-1$ the non-trivial element of the
kernel of $\tglp n\to\glp n$. Set
\[\glc n=\tglp n\uset{\pm(1,1)}\times\U1.\]
Then there is a short exact sequence 
\[0\to\z/2\to\glc n\to\glp n\times\U1\to1,\]
and $\Spinc n$ is canonically isomorphic to the preimage of $\SO n$ by the
projection $\glc n\to\glp n$. 

\begin{defn}
{\rm By a $\spc$ structure $\s$ on $X$ we mean a principal $\glc n$--bundle
$\pglc\to X$ together with a $\glc n$ equivariant map
$\pglc\to\pglp$ which covers the identity on $X$. If $\s'$ is another
$\spc$ structure on $X$ given by $\pglc'\to\pglp$ then $\s$ and
$\s'$ are called {\em isomorphic} if there is a $\U1$ equivariant map
$\pglc'\to\pglc$ which covers the identity on $\pglp$.}
\end{defn}

The natural $\U1$--bundle associated to $\pglc$ is denoted $\cll$, and
the Chern class $c_1(\cll)$ is called the
{\em canonical class} of the $\spc$ structure.

Now suppose $X$ is equipped with a Riemannian metric, and let $\pso$ be
its bundle of positive orthonormal frames, which is a principal $\SO n$--bundle.
Then the preimage $\pspc$ of $\pso$ by the projection $\pglc\to\pglp$ is
a principal $\Spinc n$--bundle over $X$, 
ie a $\spc$ structure of $X$ in 
the sense of \cite{LM}. Conversely, $\pglc$ is isomorphic to 
$\pspc\uset{\Spc n}\times\glc n$. Thus there is a natural
1--1 correspondence between
(isomorphism classes of) $\spc$ structures of the smooth oriented
manifold $X$ as defined above, and $\spc$ structures of the oriented
Riemannian manifold $X$ in the sense of \cite{LM}.

By a {\em spin connection} in $\pspc$ 
we shall mean a connection in $\pspc$ that maps
to the Levi--Civita connection in $\pso$.
If $A$ is a spin connection in $\pspc$ 
then $\hatf_A$ will denote the $i\R$ component of the curvature of $A$ with
respect to the isomorphism of Lie algebras
\[\text{spin}(n)\oplus i\R\oset\approx\to\text{spin}^c(n)\]
defined by the double cover $\Spin n\times\U1\to\Spinc n$.
In terms of the induced connection $\check A$ in $\cll$ one has
\[\hatf_A=\frac12 F_{\check A}.\]
If $A,A'$ are spin connections in $\pspc$ then we regard $A-A'$ as an
element of $i\Om^1_X$.

The results of Section~\ref{sec:orbit} carry over to spaces of configurations
$(A,\Phi)$ where $A$ is a spin connection in $\pspc$ and $\Phi$ a section of
some complex vector bundle $\bs\to X$.

When the $\spc$ structure on $X$ is understood then we will say
"spin connection over $X$" instead of "spin connection in $\pspc$".

If $n$ is even then the complex Clifford algebra
$\cocl(n)$ has up to equivalence exactly one
irreducible complex representation. Let $\bs$ denote the associated spin
bundle over $X$. Then the eigenspaces of the
complex volume element $\om_{\co}$ in $\cocl(n)$ defines a
splitting $\bs=\bs^+\oplus\bs^-$ (see \cite{LM}).

If $n$ is odd then $\cocl(n)$
has up to equivalence two irreducible complex representations $\rho_1,\rho_2$.
These restrict to equivalent representations of $\Spc(n)$, so one gets a
well-defined spin bundle $\bs$ for any $\spc$ structure on $X$ \cite{LM}.
If $\al$ is the unique automorphism of $\cocl(n)$ whose restriction to $\R^n$
is multiplication by $-1$ then $\rho_1\approx\rho_2\circ\al$. Hence if $A$
is any spin connection over $X$ then the sign of the Dirac operator $D_A$ 
depends on the choice of $\rho_j$.
To remove this ambiguity we decree that Clifford multiplication of $TX$ on
$\bs$ is to be defined using the representation $\rho_j$ satisfying
$\rho_j(\om_{\co})=1$.

In the case of a Riemannian product $\R\times X$ there is a natural
1--1 correspondence $\Spc(\R\times X)=\Spc(X)$, and we can identify
\[\cll_{\R\times X}=\R\times\cll_X.\]
If $A$ is a spin connection over $\R\times X$ then $A|_{\{t\}\times X}$
will denote the spin connection $B$ over $X$ satisfying
$\check A|_{\{t\}\times X}=\check B$.

When $n$ is odd then we can also identify
\begin{equation}\label{splrs}
\bs^+_{\R\times X}=\R\times\bs_X.
\end{equation}
If $e$ is a tangent vector on $X$ then Clifford multiplication
with $e$ on $\bs_X$ corresponds to multiplication with $e_0e$ on
$\bs^+_{\R\times X}$, where $e_0$ is the positively oriented unit tangent
vector on $\R$. Therefore, reversing the orientation of $X$ changes the
sign of the Dirac operator on $X$.

 From now on, to avoid confusion we will use $\prtl_B$ to denote
the Dirac operator over a $3$--manifold with spin connection $B$,
while the notation $D_A$ will be reserved for
Dirac operators over $4$--manifolds.

By a {\em configuration} over a $\spc$ $3$--manifold $Y$ we shall mean
a pair $(B,\Psi)$ where $B$ is a spin connection over $Y$ and $\Psi$ a 
section of the spin bundle $\bs_Y$. 
By a configuration over a $\spc$ $4$--manifold $X$ we mean
a pair $(A,\Phi)$ where $A$ is a spin connection over $X$ and $\Phi$ a 
section of the positive spin bundle $\bs_X^+$.

\subsection{The Chern--Simons--Dirac functional}
\label{subsec:csd}

Let $Y$ be a closed Riemannian $\spc$ $3$--manifold and $\eta$ 
a closed $2$--form on $Y$ of class $C^1$.
Fix a smooth reference spin
connection $B\rr$ over $Y$ and for any configuration
$(B,\Psi)$ over $Y$ define the Chern--Simons--Dirac functional $\cd=\cd_\eta$ by
\[\cd(B,\Psi)=-\frac12\int_Y(\hatf_B+\hatf_{B\rr}+2i\eta)\wedge(B-B\rr)
-\frac12\int_Y\la\prtl_B\Psi,\Psi\ra.\]
Here and elsewhere $\la\cdot,\cdot\ra$ denotes Euclidean inner products, while
$\la\cdot,\cdot\ra_{\co}$ denotes Hermitian inner products.
Note that reversing the orientation of $Y$ changes the sign of $\cd$.
Let $\cc=\cc_Y$ denote the space of $L^2_1$ configurations $(B,\Psi)$. Then
$\cd$ defines a smooth map $\cc_Y\to\R$ which has an
$L^2$ gradient
\[\nabla\cd_{(B,\Psi)}
=(*(\hatf_B+i\eta)-\frac12\si(\Psi,\Psi),-\prtl_B\Psi).\]
If $\{a_j\}$ is a local orthonormal basis of imaginary-valued $1$--forms
on $Y$ then
\[\si(\phi,\psi)=\sum_{j=1}^3\la a_j\phi,\psi\ra a_j.\]
Here and elsewhere the inner products are Euclidean unless otherwise
specified.
Since $\nabla\cd$ is independent of $B\rr$, $\cd$ is independent
of $B\rr$ up to additive constants. 
If $u\col Y\to\U1$ then
\[\cd(u(S))-\cd(S)=\int_Y(\hatf_B+i\eta)\wedge u\inv du
=2\pi\int_Y\ti\eta\wedge[u],\]
where $[u]\in H^1(Y)$ is the pull-back by $u$ of the fundamental class of
$\U1$, and $\ti\eta$ is as in \Ref{ti-eta}.

The invariance of $\cd$ under null-homotopic gauge transformations imply
\begin{equation}\label{ipncd}
\ci^*_\Psi\nabla\cd_{(B,\Psi)}=0.
\end{equation}
Let $H_{(B,\Psi)}\col L^2_1\to L^2$ be the derivative of $\nabla\cd\col \cc\to L^2$
at $(B,\Psi)$, ie
\[H_{(B,\Psi)}(b,\psi)=(*db-\si(\Psi,\psi),-b\Psi-\prtl_B\psi).\]
Note that $H_{(B,\Psi)}$ is formally self-adjoint, and
$H_{(B,\Psi)}\ci_\Psi=0$ if $\prtl_B\Psi=0$. As in \cite{Fr1},
a critical point $(B,\Psi)$ of $\cd$ is called {\em non-degenerate} if
the kernel of $\ci^*_\Psi+H_{(B,\Psi)}$ in $L^2_1$ is zero,
or equivalently, if $\ci_\Psi+H_{(B,\Psi)}\col L^2_1\to L^2$ is surjective.
Note that if $\eta$
is smooth then any critical point of $\cd_\eta$ has a smooth representative.

Let $\cg$ be the Hilbert Lie group of $L^2_2$
maps $Y\to\U1$, and $\cg_0\subset\cg$
the subgroup of null-homotopic maps. Set
\[\cb=\cc/\cg,\quad\ti\cb=\cc/\cg_0.\]
Then $\cd$ descends to a continuous map $\ti\cb\to\R$ which we also denote
by $\cd$. If Condition~(O1) holds (which we always assume when no statement
to the contrary is made) then there is a real number $q$ such that
\[\cd(\cg S)=\cd(S)+q\z\]
for all configurations $S$.
If (O1) does not hold then $\cd(\cg S)$ is a dense subset
of $\R$. 

If $S$ is any smooth configuration over a band $(a,b)\times Y$, with
$a<b$, let $\nabla\cd_S$ be the section of the bundle
$\pi_2^*(\bs_Y\oplus i\La^1_Y)$ over $(a,b)\times Y$ such that
$\nabla\cd_S|_{\{t\}\times Y}=\nabla\cd_{S_t}$. Here $\pi_2\col \ry\to Y$
is the projection. Note that $S\mapsto\nabla\cd_S$ extends to a smooth
map $L^2_1\to L^2$.

Although we will normally work with $L^2_1$ configurations over $Y$, 
the following lemma is sometimes useful.

\begin{lemma}\label{cd-ext}
$\cd$ extends to a smooth function on the
space of $L^2_{1/2}$ configurations over $Y$.
\end{lemma}

\begin{proof}
The solution to the Dirichlet problem provides
bounded operators
\[E\col L^2_{1/2}(Y)\to L^2_1(\rpy)\]
such that, for any $f\in L^2_{1/2}(Y)$, the function $Ef$ restricts to 
$f$ on $\{0\}\times Y$ and vanishes on $(1,\infty)\times Y$, and 
$Ef$ is smooth whenever $f$ is smooth. (see \cite[p\,307]{Taylor1}).
Similar extension maps can clearly be defined for
configurations over $Y$. The lemma now follows from the observation that
if $S$ is any smooth configuration over $[0,1]\times Y$ then
\[\cd(S_1)-\cd(S_0)=\int_{[0,1]\times Y}
\left\la\nabla\cd_S,\frac{\prtl S}{\prtl t}\right\ra,\]
and the right hand side extends to a smooth function on the space of
$L^2_1$ configurations $S$ over $[0,1]\times Y$.
\end{proof}

We will now relate the Chern--Simons--Dirac functional to the $4$--dimensional
monopole equations, cf \cite{KM4,MST}. Let $X$ be a $\spc$ Riemannian
$4$--manifold.
Given a parameter $\mu\in\Om^2(X)$ there are the
following Seiberg--Witten equations
for a configuration $(A,\Phi)$ over $X$:
\begin{equation}\label{swmu}
\begin{aligned}
(\hatf_A+i\mu)^+&=\quadr(\Phi)\\
D_A\Phi&=0,
\end{aligned}
\end{equation}
where
\[Q(\Phi)=\frac14\sum_{j=1}^3\la\al_j\Phi,\Phi\ra\al_j\]
for any local orthonormal basis $\{\al_j\}$ of imaginary-valued self-dual
$2$--forms on $X$. If $\Psi$ is another section of $\bs^+_X$
then one easily shows that
\[Q(\Phi)\Psi=\la\Psi,\Phi\ra_{\co}\Phi-\frac12|\Phi|^2\Psi.\]
Now let $X=\ry$ and for present and later use recall the standard bundle
isomorphisms
\begin{equation}\label{stand-isom}
\begin{aligned}
\rho^1\col \pi_2^*(\La^0(Y)\oplus\La^1(Y))&\to\La^1(\ry),\quad
(f,a)\mapsto f\,dt+a,\\
\rho^+\col \pi_2^*(\La^1(Y))&\to\La^+(\ry),\quad
a\mapsto\frac12(dt\wedge a+*_Ya).
\end{aligned}
\end{equation}
Here $*_Y$ is the Hodge $*$--operator on $Y$. 
Let $\mu$ be the pull-back of a $2$--form $\eta$ on $Y$.
Set $\cd=\csd\eta$. Let $S=(A,\Phi)$
be any smooth configuration over $\ry$ such that $A$ is in temporal gauge.
Under the identification $\bs^+_{\ry}=\pi_2^*(\bs_Y)$ we have
\[\rho^+(\si(\Phi,\Phi))=2Q(\Phi).\]
Let $\nabla_1\cd$, $\nabla_2\cd$ denote the $1$--form and spinor parts of 
$\nabla\cd$, respectively. Then
\begin{equation}\label{pert1}
\begin{aligned}
\rho^+\left(\frac{\prtl A}{\prtl t}+\nabla_1\cd_S\right)
&=(\hatf_A+i\pi^*_2\eta)^+-\quadr(\Phi)\\
\frac{\prtl\Phi}{\prtl t}+\nabla_2\cd_S&=-dt\cdot D_A\Phi,
\end{aligned}
\end{equation}
Thus, the downward gradient flow equation
\[\frac{\prtl S}{\prtl t}+\nabla\cd_S=0\]
is equivalent to the Seiberg--Witten equations \Ref{swmu}.

\subsection{Perturbations}\label{subsec:pert}

For transversality reasons we will, as in \cite{Fr1}, add further small 
perturbations to the Seiberg--Witten equation over $\ry$. 
The precise shape of these perturbations will depend on the situation
considered.
At this point we will merely describe
a set of properties of these perturbations which will suffice for the Fredholm,
compactness, and gluing theory. 

To any $L^2_1$
configuration $S$ over the band $(-\frac12,\frac12)\times Y$
there will be associated
an element $h(S)\in\R^N$, where $N\ge1$ will depend on the situation
considered. If $S$ is an $L^2_1$
configuration over $\BB^+=(a-\frac12,b+\frac12)\times Y$ where
$-\infty<a<b<\infty$ then the corresponding function
\[h_S\col [a,b]\to\R^N\]
given by $h_S(t)=h(S|_{(t-1/2,t+1/2)\times Y})$ will be smooth. These functions
$h_S$ will have the following properties.
Let $S\rr$ be a smooth reference configuration over $\BB^+$.
\begin{enumerate}
\item[(P1)] For $0\le k<\infty$ the assignment $s\mapsto h_{S\rr+s}$
defines a smooth map
$L^2_1\to C^k$ whose image is a bounded set.
\item[(P2)] If $S_n\to S$ weakly in $L^2_1$ then $\|h_{S_n}-h_S\|_{C^k}\to0$
for every $k\ge0$.
\item[(P3)] $h_S$ is gauge invariant, ie $h_S=h_{u(S)}$ for any smooth 
gauge transformation $u$.
\end{enumerate}
We will also choose a compact codimension~$0$ submanifold
$\Xi\subset\R^N$ which does not
contain $h(\ual)$ for any critical point $\al$, where $\ual$ is the
translationary invariant configuration over $\ry$ determined by $\al$.
Let $\tPrt=\tPrt_Y$ denote the space of all (smooth) $2$--forms on
$\R^N\times Y$ supported in $\Xi\times Y$.
For any $S$ as above and any $\prt\in\tPrt$
let $h_{S,\prt}\in\Om^2([a,b]\times Y)$
denote the pull-back of $\prt$ by the map $h_S\times\id$. It is clear that
$h_{S,\prt}(t,y)=0$ if $h_S(t)\not\in\Xi$.
Moreover,
\begin{equation}\label{es-est}
\|h_{S,\prt}\|_{C^k}\le\ga_k\|\prt\|_{C^k}
\end{equation}
where the constant $\ga_k$ is independent of $S,\prt$.

Now let $-\infty\le a<b\le\infty$ and $\BB=(a,b)\times Y$.
If $\fq\col \BB\to\R$ is a smooth function
then by a {\em $(\prt,\fq)$--monopole} over $\BB$ we shall mean a configuration
$S=(A,\Phi)$ over $\BB^+=(a-\frac12,b+\frac12)\times Y$
(smooth, unless otherwise stated)
which satisfies the equations
\begin{equation}\label{pert2}
\begin{aligned}
(\hatf_A+i\pi^*_2\eta+i\fq h_{S,\prt})^+&=\quadr(\Phi)\\
D_A\Phi&=0,
\end{aligned}
\end{equation}
over $\BB$, where $\eta$ is as before.
If $A$ is in temporal gauge then these equations can also be expressed as
\begin{equation}\label{pst}
\frac{\prtl S_t}{\prtl t}=-\nabla\cd_{S_t}+E_S(t),
\end{equation}
where the perturbation term $E_S(t)$ depends only on the restriction of
$S$ to $(t-\frac12,t+\frac12)\times Y$.

To reduce the number
of constants later, we will always assume that $\fq$ and its differential
$d\fq$ are pointwise bounded (in norm) by $1$ everywhere.
Note that if $\fq$ is constant then the
equations~\Ref{pert2} are translationary invariant. A $(\prt,\fq)$--monopole
with $\fq h_{S,\prt}=0$ is called a {\em genuine} monopole.
In expressions like
$\|F_A\|_2$ and $\|\Phi\|_\infty$ the norms will usually be taken over $\BB$.

For the transversality theory in Section~\ref{section:transv} we will need to
choose a suitable Banach space $\Prt=\Prt_Y$ of forms $\prt$ as above
(of some given regularity). It will be essential that
\begin{equation}\label{sprt}
\cc(\BB^+)\times\Prt\to L^p(\BB,\La^2),\quad(S,\prt)\mapsto h_{S,\prt}
\end{equation}
be a smooth map when $a,b$ are finite (here $p>2$ is the exponent used
in defining the configuration space $\cc(\BB^+)$).
Now, one cannot expect $h_{S,\prt}$ 
to be smooth in $S$
unless $\prt$ is smooth in the $\R^N$ direction (this point was overlooked
in \cite{Fr1}). It seems natural then to
look for a suitable space $\Prt$ consisting of {\em smooth} forms 
$\prt$. Such a $\Prt$ will be provided by Lemma~\ref{banach-construct}.
The topology on $\Prt$ will be stronger
than the $C^\infty$ topology, ie stronger than the $C^k$ topology for
every $k$. The smoothness of the map
\Ref{sprt} is then an easy consequence of property~(P1)
above and the next lemma. 

\begin{lemma}\label{smooth-comp}
Let $A$ be a topological space, $U$ a Banach space, and $K\subset\R^n$
a compact subset. Then the composition map
\[C_B(A,\R^n)\times C^k(\R^n,U)_K\to C_B(A,U)\]
is of class $C^{k-1}$ for any natural number $k$. Here $C_B(A,\sdot)$ denotes
the sup\-re\-mum-normed space
of bounded continuous maps from $A$ into the indicated space, 
and $C^k(\R^n,U)_K$ is the space of $C^k$ maps $\R^n\to U$ with support in $K$.
\end{lemma}

\begin{proof}
This is a formal exercise in the differential calculus.
\end{proof}

\subsection{Moduli spaces}
\label{subsec:modsp}

Consider the situation of Subsection~\ref{intro:comp-single}. (We do not
assume here that (A) holds.) We will define the moduli space
$M(X;\vec\al)$. In addition to the parameter $\mu$ this will depend
on a choice of perturbation forms
$\prt_j\in\tPrt_{Y_j}$
and a smooth function $\fq\col X\to[0,1]$ such that $\|d\fq\|_\infty\le1$,
$\fq\inv(0)=X\endt{\frac32}$, and $\fq=1$ on $[3,\infty)\times Y$.

Choose a smooth reference configuration $S_o=(A_o,\Phi_o)$ over $X$ which is 
translationary invariant and in temporal gauge over the ends, and such that
$S_o|_{\{t\}\times Y_j}$ represents $\al_j\in\ti\fl_{Y_j}$. 
Let $p>4$ and choose $w$ as
in Subsection~\ref{conf-gt}. Let $\hx$ be a finite subset of $X$ and
define $\cc,\cg_\hx,\cb_\hx$ as in Subsections~\ref{conf-gt} and
\ref{local-slices}. For clarity we will sometimes write
$\cc(X;\vec\al)$ etc.
Set $\vec\prt=(\prt_1,\dots,\prt_r)$ and let
\[M_\hx(X;\vec\al)=M_\hx(X;\vec\al;\mu;\vec\prt)\subset\cb_\hx\]
be the subset of
gauge equivalence classes of solutions $S=(A,\Phi)$ (which we simply refer to
as monopoles) to the equations
\begin{equation}\label{pert3}
\begin{gathered}
\left(\hatf_A+i\mu+i\fq\sum_{j=1}^rh_{S,\prt_j}\right)^+-\quadr(\Phi)=0\\
D_A\Phi=0.
\end{gathered}
\end{equation}
It is clear that $\fq\sum_jh_{S,\prt_j}$ vanishes outside a
compact set in $X$. If it
vanishes everywhere then $S$ is called a {\em genuine} monopole. If $\hx$
is empty then we write $M=M_\hx$.

Note that different choices of $S_o$ give canonically homeomorphic
moduli spaces $M_\hx(X;\vec\al)$ (and similarly for $\cb_\hx(X;\vec\al)$).

Unless otherwise stated the forms $\mu$ and $\prt_j$ will be smooth. 
In that case every element of $M_\hx(X;\vec\al;\mu;\vec\prt)$ has a smooth
representative, and in notation like $[S]\in M$ we will often implicitly
assume that $S$ is smooth.

We define the moduli spaces $M(\al,\beta)=M(\al,\beta;\prt)$
of Subsection~\ref{intro-csd}
similarly, except that we here use the equations~\Ref{pert2} with $\fq\equiv1$.

The following estimate will be crucial in compactness arguments later.

\begin{prop}\label{phi-bound}
For any element $[A,\Phi]\in M(X;\vec\al)$ one has that either
\begin{equation*}
\Phi=0\quad\text{or}\quad
\|\Phi\|_\infty^2\le-\frac12\inf_{x\in X}\scal(x)
+4\|\mu\|_\infty+4\ga_0\max_j\|\prt_j\|_\infty,
\end{equation*}
where $\scal$ is the scalar curvature of $X$ and the constant
$\ga_0$ is as in \Ref{es-est}.
\end{prop}

\begin{proof}
Let $\psi_j$ denote the
spinor field of $A\rr|_{\{t\}\times Y_j}$. If $|\Phi|$ has a global
maximum then the conclusion of the proposition holds by the proof of
\cite[Lemma~2]{KM4}. Otherwise one must have
$\|\Phi\|_\infty=\max_j\|\psi_j\|_\infty$ because of the
Sobolev embedding $L^p_1\subset C^0$ on compact $4$--manifolds.
But the argument in \cite{KM4} applied to $\ry$ yields
\[\psi_j=0\quad\text{or}\quad
\|\psi_j\|^2_\infty\le-\frac12\inf_{x\in X}\scal(x)+4\|\mu\|_\infty\]
for each $j$, and the proposition follows.
\end{proof}

The left hand side of \Ref{pert3} can be regarded as a section
$\sw(S)=\tsw(S,\mu,\vec\prt)$ of the bundle
$\La^+\oplus\bs^-$ over $X$. It is clear that $\sw$ defines a smooth map
\[\sw\col \cc\to\lw pw,\]
which we call the {\em monopole map}. Let $D\sw$ denote
the derivative of $\sw$. We claim that
\begin{equation}\label{ci-star-dth}
\ci_\Phi^*+D\sw(S)\col \llw pw1\to\lw pw
\end{equation}
is a Fredholm operator for every $S=(A,\Phi)\in\cc$.
Note that the $\prt_j$--per\-tur\-ba\-tions in \Ref{pert3}
only contribute a compact operator, so we can take $\prt_j=0$ for each $j$. 
We first consider the case $X=\ry$, with $\mu=\pi_2^*\eta$ as before.
By means of the isomorphisms \Ref{stand-isom}, \Ref{splrs} and the isomorphism
$\bs^+\to\bs^-$, $\phi\mapsto dt\cdot\phi$ we can think of $D\Theta_S$
as acting on sections of $\pi_2^*(\La^0_Y\oplus\La^1_Y\oplus\bs_Y)$.
If $A$ is in temporal gauge then a simple computation yields
\begin{equation}\label{dths}
D\Theta_S=\frac d{dt}+P_{S_t},
\end{equation}
where
\[P_{(B,\Psi)}=
\left(\begin{array}{cc}
0 & \ci^*_\Psi \\
\ci_\Psi & H_{(B,\Psi)}
\end{array}\right)\]
for any configuration $(B,\Psi)$ over $Y$. Note that $P_{(B,\Psi)}$ is
elliptic and formally self-adjoint, and if $(B,\Psi)$ is a non-degenerate
critical point of $\cd_\eta$ then $\ker P_{(B,\Psi)}=\ker\ci_\Psi$. Thus,
the structure of the linearized equations over a cylinder is analogous to
that of the instanton equations studied in \cite{D5}, and the results of
\cite{D5} carry over to show that \Ref{ci-star-dth} is a Fredholm operator.

The index of \Ref{ci-star-dth} is
independent of $S$ and is called the {\em expected dimension} of
$M(X;\vec\al)$. If $S\in\cc$ is a 
monopole and $D\sw(S)\col \llw pw1\to\lw pw$ is surjective
then $[S]$ is called a {\em regular} point of $M_\hx(X;\vec\al)$. If in addition
$S\in\cc^*_\hx$ then $[S]$ has
an open neighbourhood in $M_\hx(X;\vec\al)$ which is a smooth submanifold of
$\cb^*_\hx$ of dimension
\[\dim\,M_\hx(X;\vec\al)=\index(\ci_\Phi^*+D\sw(S))+|\hx|.\]

\section{Local compactness I}\label{sec:local-comp1}

This section provides the local compactness results needed for the proof of
Theorem~\ref{thm:neck-comp} assuming (B1).

\subsection{Compactness under curvature bounds}

For the moment let $B$ be an arbitrary compact, oriented Riemannian manifold
with boundary, and $v$ the outward unit normal vector field along $\prtl B$.
Then 
\begin{equation}\label{ell-bs}
\Om^*(B)\to\Om^*(B)\oplus\Om^*(\prtl B),
\quad\phi\mapsto((d+d^*)\phi,\iota(v)\phi)
\end{equation}
is an elliptic boundary system in the sense of \cite{horm3,A3}.
Here $\iota(v)$ is contraction with $v$.
By \cite[Theorems~20.1.2, 20.1.8]{horm3} we then have:

\begin{prop}\label{fred}
For $k\ge1$ the map \Ref{ell-bs} extends to a Fredholm operator
\[L^2_k(B,\La^*_B)\to L^2_{k-1}(B,\La^*_B)\oplus
L^2_{k-\frac12}(\prtl B,\La^*_{\prtl B})\]
whose kernel consists of $C^\infty$ forms.
\end{prop}

\begin{lemma}\label{lemma:loc-comp}
Let $X$ be a $\spc$ Riemannian $4$--manifold and $V_1\subset V_2\subset\cdots$
precompact open subsets of $X$
such that $X=\cup_jV_j$.
For $n=1,2,\dots$ let $\mu_n$ be a $2$--form on $V_n$,
and $S_n=(A_n,\Phi_n)$ a smooth solution to the
Seiberg--Witten equations~\Ref{swmu} over $V_n$
with $\mu=\mu_n$. Let $q>4$. Then there exist a subsequence $\{n_j\}$
and for each $j$ a smooth $u_j\col V_j\to\U1$ with the following 
significance. If $k$ is any non-negative integer such that
\begin{equation}\label{pfm}
\sup_{n\ge j}\left(\|\Phi_n\|_{L^q(V_j)}+\|\hat F(A_n)\|_{L^2(V_j)}+
\|\mu_n\|_{C^k(V_j)}\right)<\infty
\end{equation}
for every positive integer $j$ then for every $p\ge1$ one has that
$u_j(S_{n_j})$ converges weakly in
$L^p_{k+1}$ and strongly in $L^p_k$ over compact subsets of $X$
as $j\to\infty$.
\end{lemma}

Before giving the proof, note that the curvature term in \Ref{pfm}
cannot be omitted.
For if $\om$ is any non-zero, closed, anti-self-dual $2$--form over the
$4$--ball
$B$ then there is a sequence $A_n$ of $\U1$ connections over $B$ 
such that $F(A_n)=in\om$. If $S_n=(A_n,0)$ then there are clearly no
gauge transformations $u_n$ such that $u_n(S_n)$ converges (in any
reasonable sense) over compact subsets of $B$.

\begin{proof}
Let $B\subset X$ be a compact $4$--ball.
After trivializing $\cll$ over
$B$ we can write $A_n|_B=\ti d+a_n$, where $\ti d$ is the spin connection 
over $B$ corresponding to the product connection in $\cll|_B$.
By the solution of the Neumann problem
(see \cite{Taylor1}) there is a smooth $\xi_n\col B\to i\R$ such that
$b_n=a_n-d\xi$ satisfies
\[d^*b_n=0;\quad *b_n|_{\prtl B}=0.\]
Using the fact that $H^1(B)=0$ one easily proves that the map \Ref{ell-bs}
is injective on $\Om^1(B)$. Hence there is a constant $C$ such that
\[\|b\|_{L^2_1(B)}\le C(\|(d+d^*)b\|_{L^2(B)}+
\|{*b}|_{\prtl B}\|_{L^2_{1/2}(\prtl B)})\]
for all $b\in\Om^1(B)$. This gives
\[\|b_n\|_{L^2_1(B)}\le C\|db_n\|_{L^2(B)}=C\|\hat F(A_n)\|_{L^2(B)}.\]
Set $v_n=\exp(\xi_n)$.
It is now an exercise in bootstrapping, using
the Seiberg--Witten equations for $S_n$ and interior elliptic estimates,
to show that, for every $k\ge0$ for which \Ref{pfm} holds and
for every $p\ge1$, the sequence $v_n(S_n)=(d+b_n,v_n\Phi_n)$ is
bounded in $L^p_{k+1}$ over compact subsets
of $\text{int}(B)$. 

To complete the proof, choose a countable collection of such balls such that
the corresponding  balls of half the size cover $X$, and apply
Lemma~\ref{lemma:patch}. 
\end{proof}

\subsection{Small perturbations}
\label{subsec:curv-bounds}

If $S$ is any smooth configuration over a band $(a,b)\times Y$ with
$a<b$, the {\em energy} of $S$ is by definition
\[\text{energy}(S)=\int_{[a,b]\times Y}|\nabla\cd_S|^2.\]
If $S$ is a genuine monopole then
$\prtl_t\cd(S_t)=-\int_Y|\nabla\cd_{S_t}|^2$, and so the energy equals
$\cd(S_a)-\cd(S_b)$. If $S$ is a $(\prt,\fq)$--monopole then one no longer
expects these identities to hold, because the equation~\Ref{pst} is not of
gradient flow type. The main object of this subsection is to show that if
$\|\prt\|_{C^1}$ is sufficiently small then, under suitable assumptions, the
variation of $\cd(S_t)$ still controls the energy locally
(Proposition~\ref{grad-cd-est}), and there is a monotonicity result for
$\cd(S_t)$ (Proposition~\ref{Dcd}).

It may be worth mentioning that the somewhat technical Lemma~\ref{cd-alt}
and Proposition~\ref{grad-cd-est} are not needed in the second approach to
compactness which is the subject of Section~\ref{sec:local-comp2}.

In this subsection $\fq\col \ry\to\R$ may be any smooth function 
satisfying $\|\fq\|_\infty,\|d\fq\|_\infty\le1$. Constants will be
independent of $\fq$. The perturbation forms $\prt$ may be arbitrary 
elements of $\tPrt$.

\begin{lemma}\label{ptcurv}
There is a constant $C_0>0$ such that if $-\infty<a<b<\infty$ and
$S=(A,\Phi)$ is any $(\prt,\fq)$--monopole
over $(a,b)\times Y$ then there is a pointwise bound
\[|\hat F(A)|\le 2|\nabla\cd_S|+|\eta|+C_0|\Phi|^2+\ga_0\|\prt\|_\infty.\]
\end{lemma} 

\begin{proof}
Note that both sides of the inequality are gauge invariant, and if
$A$ is in temporal gauge then
\[F(A)=dt\wedge\frac{\prtl A_t}{\prtl t} + F_Y(A_t),\]
where $F_Y$ stands for the curvature of a connection over $Y$.
Now use inequalities~\Ref{pst} and \Ref{es-est}.
\end{proof}

\begin{lemma}\label{abx}
There exists a constant $C_1>0$ such that for any $\tau>0$
and any $(\prt,\fq)$--monopole $S$ over $(0,\tau)\times Y$ one has
\[\int_{[0,\tau]\times Y}|\nabla\cd_S|^2
\le 2(\cd(S_0)-\cd(S_\tau))+C_1^2\tau\|\prt\|^2_\infty.\]
\end{lemma}

Recall that by convention a $(\prt,\fq)$--monopole over $(0,\tau)\times Y$
is actually a configuration over $(-\frac12,\tau+\frac12)\times Y$, so the
lemma makes sense.

\begin{proof}
We may assume $S$ is in temporal gauge. Then
\begin{align*}
\cd(S_{\tau})-\cd(S_{0})&=\int_{0}^{\tau}\prtl_t\cd(S_t)\,dt\\
&=\int_{[0,\tau]\times Y}\la\nabla\cd_S,-\nabla\cd_S+E_S\ra\,dt\\
&\le\|\nabla\cd_S\|_2 (\|E_S\|_2 - \|\nabla\cd_S\|_2),
\end{align*}
where the norms on the last line are taken over
$[0,\tau]\times Y$. If $a,b,x$ are real numbers satisfying $x^2-bx-a\le0$
then
\[x^2\le2x^2-2bx+b^2\le2a+b^2.\]
Putting this together we obtain
\[\|\nabla\cd_S\|^2_2
\le2(\cd(S_0)-\cd(S_\tau))+\|E_S\|^2_2,\]
and the lemma follows from the estimate \Ref{es-est}.
\end{proof}

\begin{lemma}\label{cd-alt}
For all $C>0$ there exists an $\eps>0$ with the following
significance.
Let $\tau\ge4$, $\prt\in\tPrt$ with
$\|\prt\|_\infty\tau^{1/2}\le\eps$, and let $S=(A,\Phi)$ be a
$(\prt,\fq)$--monopole over $(0,\tau)\times Y$ satisfying
$\|\Phi\|_\infty\le C$. Then at least one of the following two statements must
hold:
\begin{enumerate}
\item[\rm(i)]$\prtl_t\cd(S_t)\le0$ for $2\le t\le \tau-2$,
\item[\rm(ii)]$\cd(S_{t_2})<\cd(S_{t_1})$ for $0\le t_1\le1$,
$\tau-1\le t_2\le \tau$.
\end{enumerate}
\end{lemma}

\begin{proof}
Given $C>0$, suppose that for $n=1,2,\dots$
there exist $\tau_n\ge4$, $\prt_n\in\tPrt$ with
$\|\prt_n\|_\infty \tau_n^{1/2}\le1/n$,
and a $(\prt_n,\fq_n)$--monopole $S_n=(A_n,\Phi_n)$ over $(0,\tau_n)\times Y$ 
satisfying $\|\Phi_n\|_\infty\le C$ such that
(i) is violated at some point
$t=t_n$ and (ii) also does not hold. By Lemma~\ref{abx}
the last assumption implies
\[\|\nabla\cd_{S_n}\|_{L^2([1,\tau_n-1]\times Y)}
\le C_1\|\prt_n\|_\infty\tau_n^{1/2}\le C_1/n.\]
For $s\in\R$ let $\ct_s\col \ry\to\ry$ be translation by $s$:
\begin{equation*}
\ct_s(t,y)=(t+s,y).
\end{equation*}
Given $p>2$ then by Lemmas~\ref{ptcurv} and \ref{lemma:loc-comp}
we can find $u_n\col (-1,1)\times Y\to\U1$ in $L^p_{2,\loc}$ such that
a subsequence of $u_n(\ct_{t_n}^*(S_n))$ converges weakly in $L^p_1$
over $(-\frac12,\frac12)\times Y$ to an $L^p_1$ 
solution $S'$ to the equations~\Ref{swmu} with $\mu=\pi_2^*\eta$. Then
$\nabla\cd_{S'}=0$. After modifying the
gauge transformations we can even arrange that $S'$ is smooth and
in temporal gauge, in which case there is a critical point $\al$ of
$\cd$ such that $S'(t)\equiv\al$.
After relabelling the subsequence above consecutively we then have
\[h_{S_n}(t_n)\to h_{S'}(0)\not\in\Xi.\]
Since $\Xi$ is closed, $h_{S_n}(t_n)\not\in\Xi$ for $n$ sufficiently large.
Hence,
$\prtl_t|_{t_n}\cd(S_n(t))=-\|\nabla\cd_{S_n(t_n)}\|^2\le0$, 
which is a contradiction.
\end{proof}

\begin{prop}\label{grad-cd-est}
For any constant $C>0$ there exist $C',\del>0$ such that if $S=(A,\Phi)$
is any $(\prt,\fq)$--monopole over $(-2,T+4)\times Y$ where $T\ge2$,
$\|\prt\|_\infty\le\del$, and $\|\Phi\|_\infty\le C$, then for $1\le t\le T-1$
one has
\[\int_{[t-1,t+1]\times Y}|\nabla\cd_S|^2
\le2\left(\sup_{0\le r\le1}\cd(S_{-r})
-\inf_{0\le r\le4}\cd(S_{T+r})\right)
+C'\|\prt\|^2_\infty.\]
\end{prop}

\begin{proof}
Choose $\eps>0$ such that the conclusion of Lemma~\ref{cd-alt} holds
(with this constant $C$), and set $\del=\eps/\sqrt6$.
We construct a sequence $t_0,\dots,t_m$ of
real numbers, for some $m\ge1$, with the following properties:
\begin{enumerate}
\item[(i)]$-1\le t_0\le0$ and $T\le t_m\le T+4$,
\item[(ii)]For $i=1,\dots,m$ one has $1\le t_i-t_{i-1}\le5$ and
$\cd(S_{t_i})\le\cd(S_{t_{i-1}})$.
\end{enumerate}
The lemma will then follow from Lemma~\ref{abx}.
The $t_i$'s will be constructed inductively, and this will involve an 
auxiliary sequence $t'_0,\dots,t'_{m+1}$. Set $t_{-1}=t'_0=0$.

Now suppose $t_{i-1},t'_i$ have been constructed for $0\le i\le j$. If
$t'_j\ge T$ then we set $t_j=t'_j$ and
$m=j$, and the construction is finished. If
$t'_j<T$ then we define $t_j,t'_{j+1}$ as follows:

If $\prtl_t\cd(S_t)\le0$ for all $t\in[t'_j,t'_j+2]$ set
$t_j=t'_j$ and $t'_{j+1}=t'_j+2$; otherwise set $t_j=t'_j-1$ and
$t'_{j+1}=t'_j+4$.

Then (i) and (ii) are satisfied, by Lemma~\ref{cd-alt}.
\end{proof}

\begin{prop}\label{Dcd}
For all $C>0$ there exists a $\del>0$ such that if $S=(A,\Phi)$ is any 
$(\prt,\fq)$--monopole in temporal gauge over $(-1,1)\times Y$ such that
$\|\prt\|_{C^1}\le\del$,
$\|\Phi\|_\infty\le C$, and $\|\nabla\cd_S\|_2\le C$ then the
following holds:
Either $\prtl_t|_0\cd(S_t)<0$, or there is a
critical point $\al$ such that $S_t=\al$ for $|t|\le\frac12$.
\end{prop}

\begin{proof}
First observe that if $S$ is any $C^1$ configuration over $\ry$ then
\[\cd(S_{t_2})-\cd(S_{t_1})=
\int_{t_1}^{t_2}\int_Y\la\nabla\cd_{S_t},\prtl_tS_t\ra dy\,dt,\]
hence $\cd(S_t)$ is a $C^1$ function of $t$ whose derivative can be expressed
in terms of the $L^2$ gradient of $\cd$ as usual.

Now suppose there is a $C>0$ and for $n=1,2,\dots$ a $\prt_n\in\tPrt$
and a $(\prt_n,\fq_n)$--monopole $S_n=(A_n,\Phi_n)$ over
$(-1,1)\times Y$ such that
$\|\prt_n\|_{C^1}\le\frac1n$,
$\|\Phi_n\|_\infty\le C$, $\|\nabla\cd_{S_n}\|_2\le C$, and
$\prtl_t|_0\cd(S_n(t))\ge0$. Let $p>4$ and $0<\eps<\frac12$.
After passing to a subsequence and relabelling consecutively
we can find $u_n\col (-1,1)\times Y\to\U1$ in $L^p_{3,\text{loc}}$ such that
$\ti S_n=u_n(S_n)$ converges weakly in $L^p_2$, and strongly in $C^1$, over
$(-\frac12-\eps,\frac12+\eps)\times Y$
to a smooth solution $S'$ of \Ref{swmu} with $\mu=\pi_2^*\eta$.
We may arrange that $S'$ is in temporal gauge. Then
\[0\le\prtl_t|_0\cd(\ti S_n(t))=
\int_Y\la\nabla\cd_{\ti S_n(0)},\prtl_t|_0\ti S_n(t)\ra
\to\prtl_t|_0\cd(S'_t).\]
But $S'$ is a genuine monopole,
so $\prtl_t\cd(S'_t)=-\|\nabla\cd_{S'_t}\|^2_2$. It also follows that
$\nabla\cd_{S'(0)}=0$, hence $\nabla\cd_{S'}=0$ in
$(-\frac12-\eps,\frac12+\eps)\times Y$ by unique
continuation as in \cite[Appendix]{Fr1}. 
Since $h_{S_n}\to h_{S'}$ uniformly in $[-\eps,\eps]$, and
$h_{S'}\equiv\text{const}\not\in\Xi$, 
the function $h_{S_n}$ maps $[-\eps,\eps]$ into
the complement of $\Xi$ when $n$ is sufficiently large. In that case,
$S_n$ restricts to a genuine monopole on $[-\eps,\eps]\times Y$, and the
assumption $\prtl_t|_0\cd(S_n(t))\ge0$ implies that $\nabla\cd_{S_n}=0$
on $[-\eps,\eps]\times Y$. Since this holds
for any $\eps\in(0,\frac12)$, the proposition follows.
\end{proof}

We say a $(\prt,\fq)$--monopole $S$ over $\rpy$ has {\em finite energy} if
$\inf_{t>0}\cd(S_t)>-\infty$. A monopole over a $4$--manifold with 
tubular ends is said to have finite energy if it has finite energy over
each end.

\begin{prop}\label{crit-lim}
Let $C,\del$ be given such that the conclusion of  Proposition~\ref{Dcd}
holds. If $S=(A,\Phi)$ is any finite energy
$(\prt,\fq)$--monopole over $\rpy$ with $\|\prt\|_{C^1}\le\del$, $\fq\equiv1$,
$\|\Phi\|_\infty\le C$, and
\[\sup_{t\ge1}\|\nabla\cd_S\|_{L^2((t-1,t+1)\times Y)}\le C\]
then the following hold:
\begin{enumerate}
\item[\rm(i)]There is a $t>0$ such that $S$ restricts to a genuine monopole on
\linebreak $(t,\infty)\times Y$,
\item[\rm(ii)]$[S_t]$ converges in $\cb_Y$ to some critical point as $t\to\infty$.
\end{enumerate}
\end{prop}

\begin{proof}
Let $p>4$. If $\{t_n\}$ is any sequence with
$t_n\to\infty$ as $n\to\infty$ then
by Lemmas~\ref{lemma:loc-comp} and \ref{ptcurv} there exist
$u_n\in L^p_{3,\loc}(\ry;\U1)$
such that a subsequence of $u_n(\ct_{t_n}^*S)$ converges weakly in $L^p_2$
(hence strongly in $C^1$) over compact subsets of $\ry$ to
a smooth $(\prt,\fq)$--monopole $S'$ in temporal gauge.
Proposition~\ref{Dcd} guarantees that $\prtl_t\cd(S_t)\le0$ for $t\ge1$, so the
finite energy assumption implies that $\cd(S'_t)$ is constant.
By Proposition~\ref{Dcd} there is a critical point $\al$ such that
$S'_t=\al$ for all $t$. This implies (i) by choice of the set $\Xi$
(see Subsection~\ref{subsec:pert}). Part~(ii) follows by a continuity 
argument from the facts that $\cb_Y$ contains only finitely many 
critical points, and the topology on $\cb_Y$ defined by the $L^2$--metric
is weaker than the usual topology.
\end{proof}

The following corollary of Lemma~\ref{cd-ext} 
shows that elements of the moduli spaces defined in 
Subsection~\ref{subsec:modsp} have finite energy. 

\begin{lemma}\label{cdlim}
Let $S$ be a configuration over $\orp\times Y$ and $\al$ a critical
point of $\cd$ such that $S-\ual\in L^p_1$ for some $p\ge2$. Then
\[\cd(S_t)\to\cd(\al)~\text{as $t\to\infty$.}\]
\end{lemma}

\subsection{Neck-stretching I}
\label{subsec:neckstr1}

This subsection contains the crucial step in the proof of
Theorem~\ref{thm:neck-comp} assuming (B1), namely what should be thought of
as a global energy bound.

\begin{lemma}
\label{lemma:bockner}
Let $X$ be as in Subsection~\ref{intro:comp-single} and set $Z=X\endt 1$.
We identify $Y=\prtl Z$.
Let $\mu_1,\mu_2\in\Om^2(Z)$, where $d\mu_1=0$. Set $\eta=\mu_1|_Y$
and $\mu=\mu_1+\mu_2$. Let $A\rr$ be a spin connection over $Z$, and let
the Chern--Simons--Dirac functional $\cd_\eta$ over $Y$
be defined in terms of the reference connection $B\rr=A\rr|_Y$. Then 
for all configurations $S=(A,\Phi)$ over $Z$ which satisfy the
monopole equations~\Ref{swmu} one has
\begin{multline*}
\left|2\cd_\eta(S|_Y)+\int_Z\left(|\nabla_A\Phi|^2+|\hatf_A+i\mu_1|^2\right)
\right|\\
\le C\text{Vol}(Z)\left(1+\|\Phi\|^2_\infty+\|F_{A\rr}\|_\infty
+\|\mu_1\|_\infty+\|\mu_2\|_\infty+\|\scal\|_\infty\right)^2,
\end{multline*}
for some universal constant $C$, where $\scal$ is the scalar curvature of $Z$. 
\end{lemma}
The upper bound given here is not optimal but suffices for our purposes.

\proof
Set $F'_A=\hatf_A+i\mu_1$ and define $F'_{A\rr}$ similarly. Set $B=A|_Y$.
Without the assumption $d\mu_1=0$ we have
\begin{align*}
\int_Z|F'_A|^2
&=\int_Z(2|(F'_A)^+|^2+F'_A\wedge F'_A)\\
&=\int_Z\left(2|\quadr(\Phi)-i\mu_2^+|^2
+F'_{A\rr}\wedge F'_{A\rr}-2i d\mu_1\wedge(A-A\rr)\right)\\
&+\int_Y(\hatf_B+\hatf_{B\rr}+2i\eta)\wedge(B-B\rr).
\end{align*}
Without loss of generality we may assume $A$ is in temporal gauge
over the collar $\iota([0,1]\times Y)$. By the Weitzenb\"ock formula we have
\[0=D_A^2\Phi=\nabla_A^*\nabla_A\Phi+\hatf_A^++\frac\scal4.\]
This gives
\begin{align*}
\int_Z|\nabla_A\Phi|^2&=\int_Z\la\nabla_A^*\nabla_A\Phi,\Phi\ra
+\int_Y\la\prtl_t\Phi,\Phi\ra\\
&=\int_Z\left(-\frac12|\Phi|^4-\frac\scal4|\Phi|^2+\la i\mu^+
\Phi,\Phi\ra\right)
+\int_Y\la\prtl_B\Phi,\Phi\ra.\qquad\qquad\quad\qed
\end{align*}

Consider now the situation of Subsection~\ref{intro:comp-neck}. If (B1) holds
then we can find a closed $2$--form $\mu_1$ on $X$ whose restriction to
$\R_+\times(\pm Y_j)$ is the pull-back of $\eta_j$, and whose restriction to
$\rpy'_j$ is the pull-back of $\eta'_j$. 
From Lemma~\ref{lemma:bockner} we deduce:

\begin{prop}\label{prop:facd-bound}
For every constant $C_1<\infty$ there
exists a constant $C_2<\infty$ with the following significance.
Suppose we are given 
\begin{itemize}
\item $\tau,C_0<\infty$ and an $r$--tuple $T$ such
that $\tau\le T_j$ for each $j$,
\item real numbers $\tau^\pm_j$, $1\le j\le r$ and $\tau'_j$, $1\le j\le r'$
satisfying $0\le T_j-\tau^\pm_j\le\tau$ and $0\le\tau'_j\le\tau$.
\end{itemize}
Let $Z$ be the result of deleting from $\xt$
all the necks $(-\tau^-_j,\tau^+_j)\times Y_j$, $1\le j\le r$ and
all the ends $(\tau'_j,\infty)\times Y'_j$, $1\le j\le r'$.
Then for any configuration $S=(A,\Phi)$ representing an element of a moduli
space $M(\xt;\vec\al';\mu;\vec\prt,\vec\prt')$ where
$\sum_{j=1}^{r'}\cd(\al'_j)>-C_0$ and $\prt_j,\prt'_j,\mu$ all have $L^\infty$
norm $<C_1$ one has that
\begin{align*}
\int_Z(|\nabla_A\Phi|^2+|\hatf_A+i\mu_1|^2)
&+2\sum_{j=1}^r\left(\cd(S|_{\{-\tau^-_j\}\times Y_j})
-\cd(S|_{\{\tau^+_j\}\times Y_j})\right)\\
&+2\sum_{j=1}^{r'}\left(\cd(S|_{\{\tau'_j\}\times Y'_j})
-\cd(\al'_j)\right) < C_2(1+\tau)+2C_0.
\end{align*}
\end{prop}

Thus, if each $\prt_j$, $\prt'_j$ has sufficiently
small $L^\infty$ norm then
Lemma~\ref{cdlim} and
Proposition~\ref{grad-cd-est} provides local energy bounds over necks and ends
for such monopoles. (To apply Proposition~\ref{grad-cd-est} one can take
$\tau^\pm_j$ to be the point $t$ in a suitable interval where
$\pm\cd(S|_{\{t\}\times Y_j})$ attains its maximum, and similarly for 
$\tau'_j$.)
Moreover, if these perturbation forms
have sufficiently small $C^1$ norms then we can apply
Proposition~\ref{Dcd} over necks and ends.
(How small the $C^1$ norms have to be depends on $C_0$.)

\section{Local compactness II}\label{sec:local-comp2}

This section, which is logically independent from
Section~\ref{sec:local-comp1},
provides the local compactness results needed for the proof of
Theorem~\ref{thm:neck-comp} assuming (B2).

While the main result of this section, Proposition~\ref{prop:xe-conv},
is essentially concerned with local convergence of monopoles, the
arguments will, in contrast to those of Section~\ref{sec:local-comp1},
be of a global nature. In particular, function spaces over manifolds
with tubular ends will play a central role.

\subsection{Hodge theory for the operator $-d^*+d^+$}\label{subsec:hodge}

In this subsection we will study the kernel (in certain function spaces)
of the elliptic operator
\begin{equation}\label{cDdef}
\cD=-d^*+d^+\col \Om^1(X)\to\Om^0(X)\oplus\Om^+(X),
\end{equation}
where $X$ is an oriented Riemannian $4$--manifold with
tubular ends. 
The notation $\ker(\cD)$ will refer to the
kernel of $\cD$ in the space $\Om^1(X)$ of all smooth $1$--forms,
where $X$ will be understood from the context.
The results of this subsection complement those of \cite{D5}.

We begin with the
case of a half-infinite cylinder $X=\rpy$, where $Y$ is any closed, oriented,
connected Riemannian $3$--manifold.
Under the isomorphisms \Ref{stand-isom}
there is the identification $\cD=\pp t+P$ over $\rpy$,
where $P$ is the self-adjoint elliptic operator
\[P=\left(\begin{array}{cc}
0 & -d^* \\
-d & *d
\end{array}\right)\]
acting on sections of $\La^0(Y)\oplus\La^1(Y)$ (cf \Ref{dths}).
Since $P^2$ is the Hodge Laplacian,
\[\ker(P)=H^0(Y)\oplus H^1(Y).\]
Let $\{h_\nu\}$ be a maximal orthonormal set of eigenvectors of $P$,
say $Ph_\nu=\lla_\nu h_\nu$.

Given a smooth $1$--form $a$ over $\rpy$ we can express it as
$a=\sum_\nu f_\nu h_\nu$, where $f_\nu\col \R_+\to\R$.
If $\cD a=0$ then $f_\nu(t)=c_\nu e^{-\lla_\nu t}$ for some constant $c_\nu$.
If in addition $a\in L^p$ for some $p\ge1$ then $f_\nu\in L^p$
for all $\nu$, hence
$f_\nu\equiv0$ when $\lla_\nu\le0$.
Elliptic inequalities for $\cD$ then show that
$a$ decays exponentially, or more precisely,
\[|(\nabla^ja)_{(t,y)}|\le\beta_je^{-\del t}\]
for $(t,y)\in\rpy$ and $j\ge0$, where $\beta_j$ is a
constant and $\del$ the smallest positive eigenvalue of $P$.

Now let $\si>0$ be a small constant and
$a\in\ker(\cD)\cap L^{p,-\si}$. Arguing as above we find that
\begin{equation}\label{a-decomp}
a=b+c\,dt+\pi^*\psi,
\end{equation}
where $b$ is an exponentially decaying form,
$c$ a constant, $\pi\col \rpy\to Y$, and $\psi\in\Om^1(Y)$ harmonic.

We now turn to the case when $X$ is an oriented, connected
Riemannian $4$--manifold with tubular end $\orpy$ (so $X\setminus\rpy$ is
compact). Let $Y_1,\dots,Y_r$ be the connected components of $Y$ and set
\[Y'=\bigcup_{j=1}^sY_j,\qquad Y''=Y\setminus Y',\]
where $0\le s\le r$.
Let $\si>0$ be a small constant and $\ka\col X\to\R$ a smooth function such that
\[\ka=\begin{cases}
-\si t & \text{on $\rpy'$}\\
\si t  & \text{on $\rpy''$},
\end{cases}\]
where $t$ is the $\R_+$ coordinate. Our main
goal in this subsection is to describe $\ker(\cD)\cap\lw p\ka$.

We claim that all elements $a\in\ker(\cD)\cap\lw p\ka$ are closed.
To see this, note first that the decomposition~\Ref{a-decomp} shows that
$a$ is pointwise bounded, and $da$ decays exponentially over the ends.
Applying the proof of \cite[Proposition~1.1.19]{DK} to
$X\endt T=X\setminus(T,\infty)\times Y$ we get
\begin{align*}
\int_{X\endt T}(|d^+a|^2-|d^-a|^2)&=\int_{X\endt T}da\wedge da
=\int_{X\endt T}d(a\wedge da)\\
&=\int_{\prtl X\endt T}a\wedge da\to0
\quad\text{as $T\to\infty$.}
\end{align*}
Since $d^+a=0$, we conclude that $da=0$. 

Fix $\tau\ge0$ and for any $a\in\Om^1(X)$ and $j=1,\dots,r$ set
\begin{equation}\label{R-def}
R_ja=\int_{\{\tau\}\times Y_j}*a.
\end{equation}
Recall that $d^*=-*d*$ on $1$--forms, so if $d^*a=0$ then $R_ja$ is independent
of $\tau$. Therefore, if $a\in\ker(\cD)\cap\lw p\ka$ then $R_ja=0$ for $j>s$,
hence
\[\sum_{j=1}^sR_ja=\int_{\prtl X\endt \tau}*a=\int_{X\endt \tau}d{*}a=0.\]
Set
\[\Xi=\{(z_1,\dots,z_s)\in\R^s\st\sum_jz_j=0\}.\]

\begin{prop}\label{kerD}
In the situation above the map
\begin{align*}
\al\col \ker(\cD)\cap L^{p,\ka}&\to\ker(H^1(X)\to H^1(Y''))\oplus\Xi,\\
a&\mapsto([a],(R_1a,\dots,R_sa))
\end{align*}
is an isomorphism.
\end{prop}

\begin{proof}
We first prove $\al$ is injective. Suppose $\al(a)=0$. Then $a=df$
for some function $f$ on $X$. From the
decomposition \Ref{a-decomp} we see that $a$ decays exponentially over the
ends. Hence $f$ is bounded,
in which case
\[0=\int_Xf\,d^*a=\int_X|a|^2.\]
This shows $\al$ is injective.

Next we prove $\al$ is surjective. Suppose $b\in\Om^1(X)$, $db=0$,
$[b|_{Y''}]=0$, and
$(z_1,\dots,z_s)\in\Xi$. Let $\psi\in\Om^1(Y)$ be the harmonic
form representing $[b|_Y]\in H^1(Y)$. Then
\[b|_{\rpy}=\pi^*\psi+df,\]
for some $f\col \rpy\to\R$. Choose a smooth function
$\rho\col X\to\R$ which vanishes in a neighbourhood of $X\endt 0$ and satisfies
$\rho\equiv1$ on $[\tau,\infty)\times Y$. Set $z_j=0$ for $j>s$ and let
$z$ be the function on $Y$ with $z|_{Y_j}\equiv\text{Vol}(Y_j)\inv z_j$.
Define
\[\ti b=b+d(\rho(tz-f)).\]
Then over $[\tau,\infty)\times Y$ we have $\ti b=\pi^*\psi+z\,dt$, so
$d^*\ti b=0$ in this region, and
\[\int_Xd^*\ti b=-\int_{\prtl X\endt \tau}*\ti b=-\int_Yz=0.\]
Let $\bar\ka\col X\to\R$ be a smooth function which agrees with $|\ka|$ outside
a compact set.
By Proposition~\ref{dd-image} we can find a smooth $\xi\col X\to\R$ such that
$d\xi\in L^{p,\bar\ka}_1$ and
\[d^*(\ti b+d\xi)=0.\]
Set $a=\ti b+d\xi$. Then $(d+d^*)a=0$ and $\al(a)=([b],(z_1,\dots,z_s))$.
\end{proof}

The following proposition is essentially \cite[Proposition~3.14]{D5}
and is included
here only for completeness.

\begin{prop}\label{cd-index}
If $b_1(Y)=0$ and $s=0$ then the operator
\begin{equation}\label{cdlp1}
\cD\col \llpk\to\lpk
\end{equation}
has index $-b_0(X)+b_1(X)-b^+_2(X)$.
\end{prop}

\begin{proof}
By Proposition~\ref{kerD} the dimension of the kernel of \Ref{cdlp1}
is $b_1(X)$. From Proposition~\ref{laplace-iso}~(ii) with $\bs=0$ we see
that the image of \Ref{cdlp1} is the sum of $d^*\llpk$ and $d^+\llpk$.
The codimensions of these spaces in $\lpk$ are $b_0(X)$ and $b^+_2(X)$, 
respectively.
\end{proof}

\subsection{The case of a single moduli space}

Consider the situation of Subsection~\ref{intro:comp-single}.
Initially we do not assume Condition~(A).

\begin{prop}\label{prop:loc-comp}
Fix $1<q<\infty$.
Let $\si>0$ be a small constant and $\ka\col X\to\R$ a smooth function such that
$\ka(t,y)=-\si t$ for all $(t,y)\in\rpy$.
Let $A\rr$ be a spin connection over $X$ which is translationary invariant
over the ends of $X$.
For $n=1,2,\dots$ let $S_n=(A\rr+a_n,\Phi_n)$ be a smooth
configuration over $X$ which satisfies the monopole equations
\Ref{pert3} with $\mu=\mu_n$, $\vec\prt=\vec\prt_n$. Suppose $a_n\in\llw q\ka1$
for every $n$, and $\sup_n\|\Phi_n\|_\infty<\infty$. Then there exist smooth
$u_n\col X\to\U1$ such that if $k$ is any non-negative integer with
\begin{equation}\label{mupk}
\sup_{j,n}\left(\|\mu_n\|_{C^k}+\|\prt\nj\|_{C^k}\right)<\infty
\end{equation}
then the sequence $u_n(S_n)$ is
bounded in $L^{p'}_{k+1}$ over compact subsets of $X$ for every $p'\ge1$.
\end{prop}
Before giving the proof we record the following two elementary lemmas:

\begin{lemma}\label{gen-ineq}
Let $E,F,G$ be Banach spaces, and $E\oset S\to F$ and $E\oset T\to G$
bounded linear maps. Set
\[S+T\col E\to F\oplus G,\quad x\mapsto(Sx,Tx).\]
Suppose $S$ has finite-dimensional kernel and closed range, and that
$S+T$ is injective. Then $S+T$ has closed range, hence
there is a constant $C>0$ such that
\[\|x\|\le C(\|Sx\|+\|Tx\|)\]
for all $x\in E$.
\end{lemma}

\begin{proof}
Exercise.
\end{proof}

\begin{lemma}\label{g-char}
Let $X$ be a smooth, connected manifold and $x_0\in X$. Let
$\text{Map}_0(X,\U1)$ denote the set of smooth maps $u\col X\to\U1$ such that
$u(x_0)=1$, and let $V$ denote the set of all closed $1$--forms $\phi$
on $X$ such that $[\phi]\in H^1(X;\z)$. Then
\[\text{Map}_0(X,\U1)\to V,\quad u\mapsto\frac1{2\pi i}u\inv du\]
is an isomorphism of Abelian groups.
\end{lemma}

\begin{proof}
If $\phi\in V$ define 
\[u(x)=\exp(2\pi i\int_{x_0}^x\phi),\]
where $\int_{x_0}^x\phi$ denotes the integral of $\phi$ along any path
 from $x_0$ to $x$. Then $\frac1{2\pi i}u\inv du=\phi$.
The details are left to the reader.
\end{proof}

{\bf Proof of Proposition~\ref{prop:loc-comp}}\qua
We may assume $X$ is connected and that \Ref{mupk} holds at least for $k=0$.
Choose closed
$3$--forms $\om_1,\dots,\om_{b_1(X)}$ which are supported in the interior of
$X\endt0$ and represent a basis for $H^3_c(X)$.
For any $a\in\Om^1(X)$ define the coordinates of $Ja\in\R^{b_1(X)}$ by
\[(Ja)_k=\int_Xa\wedge\om_k.\]
Then $J$ induces an isomorphism $H^1(X)\to\R^{b_1(X)}$,
by Poincar\'e duality. By Lemma~\ref{g-char}
we can find smooth $v_n\col X\to\U1$ such that $J(a_n-v_n\inv dv_n)$ is bounded
as $n\to\infty$. We can arrange that $v_n(t,y)$ is independent of $t\ge0$
for every $y\in Y$. 
Then there are
$\xi_n\in\llw q\ka2(X;i\R)$ such that $b_n=a_n-v_n\inv dv_n-d\xi_n$
satisfies
\[d^*b_n=0;\qquad R_jb_n=0,\quad j=1,\dots,r-1\]
where $R_j$ is as in \Ref{R-def}. If $r\ge1$ this follows from
Proposition~\ref{laplace-iso2}, while if $r=0$ (ie if $X$ is closed) 
it follows from Proposition~\ref{dd-image}. (By Stokes' theorem we have
$R_rb_n=0$ as well, but we don't need this.)
Note that $\xi_n$ must be smooth,
by elliptic regularity for the Laplacian $d^*d$. Set 
$u_n=\exp(\xi_n)v_n$. Then $u_n(A\rr+a_n)=A\rr+b_n$.
By Proposition~\ref{kerD} and Lemma~\ref{gen-ineq} there is a $C_1>0$
such that
\[\|b\|_{L^{q,\ka}_1}\le C_1\left(\|(d^*+d^+)b\|_{L^{q,\ka}}+
\sum_{j=1}^{r-1}|R_jb|+\|Jb\|\right)\]
for all $b\in L^{q,\ka}_1$. From inequality~\Ref{es-est} and
the curvature part of the Seiberg--Witten equations we find that 
$\sup_n\|d^+b_n\|_\infty<\infty$, hence
\[\|b_n\|_{L^{q,\ka}_1}\le C_1(\|d^+b_n\|_{L^{q,\ka}}+\|Jb_n\|)\le C_2\]
for some constant $C_2$. We can now complete the proof by bootstrapping
over compact subsets 
of $X$, using alternately the Dirac and curvature parts of the Seiberg--Witten
equation.
\endproof

Combining Proposition~\ref{prop:loc-comp} (with $k\ge1$) and
Proposition~\ref{phi-bound} we obtain, for fixed closed $2$--forms $\eta_j$
on $Y_j$:

\begin{cor}
\label{cor:energy-bound}
If (A) holds
then for every constant $C_0<\infty$ there exists a constant $C_1<\infty$
with the following
significance. Suppose $\|\mu\|_{C^1},\|\prt_j\|_{C^1}\le C_0$ for each $j$.
Then for any $\vec\al=(\al_1,\dots,\al_r)$ with $\al_j\in\ti\fl_{Y_j}$, and
any $[S]\in M(X;\vec\al;\mu;\vec\prt)$ and
$t_1,\dots,t_r\in[0,C_0]$ one has
\[\left|\sum_{j=1}^r\lla_j\cd(S|_{\{t_j\}\times Y_j})\right|
\le C_1.\]
\end{cor}
Note that if $\sum_j\lla_j\cd(\al_j)\ge-C_0$ then this gives
\begin{equation}\label{energy-bound}
\sum_{j=1}^r\lla_j\left(\cd(S|_{\{t_j\}\times Y_j})-
\cd(\al_j)\right)\le C_0+C_1.
\end{equation}

\subsection{Condition~(C)}
\label{subsec:onC}

Consider the situation in Subsection~\ref{intro:comp-neck} and suppose
$\ga$ is simply-connected and equipped with an orientation $\ortn$.
Throughout this subsection (and the next)
(co)homology groups will have real coefficients, unless otherwise indicated.

We associate to $(\ga,\ortn)$ a ``height function'', namely 
the unique integer valued function $h$ on the set of nodes
of $\ga$ whose minimum value is $0$ and which satisfies $h(e')=h(e)+1$
whenever there is an oriented edge from $e$ to $e'$. 

Let $Z\he k$ and $Z\hg k$ denote the union
of all subspaces $Z_e\subset X\gl$ where $e$ has height $k$ and $\ge k$,
respectively. 
Set
\[\prtl^-Z\he k=\bigcup_{h(e)=k}\prtl^-Z_e.\] 
For each node $e$ of $\ga$ choose a subspace $G_e\subset H_1(Z_e)$ such that
\[H_1(Z_e)=G_e\oplus\im\left(H_1(\prtl^-Z_e)\to H_1(Z_e)\right).\]
Then the natural map $F_e\to G^*_e$ is an isomorphism, where $G^*_e$ is
the dual of the vector space $G_e$.

\begin{lemma}\label{gk-inj}
The natural map $H^1(X\gl)\to\oplus_eG_e^*$ is injective.
Therefore, this map is an isomorphism if and only if $\Si(X,\ga,\ortn)=0$.
\end{lemma}

\begin{proof}
Let $N$ be the maximum value of $h$ and suppose $z\in H^1(X\gl)$
lies in the kernel of the map in the lemma. It is easy
to show, by induction on $k=0,\dots,N$, that $H^1(X\gl)\to H^1(Z\he k)$
maps $z$ to zero for each $k$. We now invoke the Mayer--Vietoris sequence
for the pair of subspaces $(Z\he{k-1},Z\hg k)$ of $X\gl$:
\begin{equation*}
H^0(Z\he{k-1})\oplus H^0(Z\hg k)\oset a\to H^0(\prtl^-Z\he k)\to
H^1(Z\hg{k-1})\oset b\to H^1(Z\he{k-1})\oplus H^1(Z\hg k).
\end{equation*}
Using the fact that $\ga$ is simply-connected it is not hard to see
that $a$ is surjective, hence $b$ is injective.
Arguing by induction on $k=N,N-1,\dots,0$ we then find that 
$H^1(X\gl)\to H^1(Z\hg k)$ maps $z$ to zero for each $k$.
\end{proof}

We will now formulate a condition on $(X,\ga)$ which is stronger than (C)
and perhaps simpler to verify.
A connected, oriented graph is called a {\em tree} if
it has a unique node (the {\em root node}) with no incoming edge, and
any other node has a unique incoming edge. 

\begin{prop}\label{CCpl}
Suppose there is an orientation $\ortn$ of $\ga$ such that $(\ga,\ortn)$
is a tree and
\[H^1(Z\hg k)\to H^1(Z\he k)\]
is surjective for all $k$. Then Condition~(C) holds.
\end{prop}

\proof
It suffices to verify that $\Si(X,\ga,\ortn)=0$.  Set
\begin{align*}
F\hg k&=\ker(H^1(Z\hg k)\to H^1(\prtl^-Z\he k)),\\
F\he k&=\ker(H^1(Z\he k)\to H^1(\prtl^-Z\he k)).
\end{align*}
The Mayer--Vietoris sequence for $(Z\he k,Z\hg{k+1})$ yields
an exact sequence
\[0\to F\hg k\to F\he k\oplus H^1(Z\hg{k+1})\to H^1(\prtl^-Z\hg{k+1}).\]
If $H^1(Z\hg k)\to H^1(Z\he k)$ is surjective then so is $F\hg k\to F\he k$,
hence $\ker(F\hg k\to F\he k)\to F\hg {k+1}$ is an isomorphism, in which case
\[\dim\,F\hg k=\dim\,F\hg {k+1}+\dim\,F\he k.\]
Therefore,
$$\Si(X,\ga,\ortn)=\dim\,H^1(X\gl)-\sum_k\dim\,F\he k
=\dim\,H^1(X\gl)-\dim\,F\hg 0=0.\eqno{\qed}$$

\subsection{Neck-stretching II}
\label{subsec:neckstr2}

Consider again the situation in Subsection~\ref{intro:comp-neck}.
The following set-up will be used in the next two lemmas. We assume that
$\ga$ is simply-connected and that $\ortn$ is an
orientation of $\ga$ with $\Si(X,\ga,\ortn)=0$. Let $1<p<\infty$.

An end of $X$ that corresponds to an edge of $\ga$ is either {\em incoming}
or {\em outgoing} depending on the orientation $o$. (These are the ends
$\R_+\times(\pm Y_j)$, but the sign here is unrelated to $o$.)
All other ends (ie $\rpy'_j$, $1\le j\le r'$)
are called {\em neutral}.

Choose subspaces $G_e$ of $H_1(Z_e)\approx H_1(X_e)$
as in the previous subsection, and set $g_e=\dim\,G_e$.
For each component $X_e$ of $X$ let $\{q_{ek}\}$ be a collection of 
closed $3$--forms on $X_e$ supported in the interior of
$Z'_e=(X_e)\endt0$ which represents a basis
for the image of $G_e$ in $H^3_c(X_e)$ under the Poincar\'e duality
isomorphism. For any $a\in\Om^1(Z'_e)$ define $J_ea\in\R^{g_e}$ by
\[(J_ea)_k=\int_{X_e}a\wedge q_{ek}.\]
For each $e$ let $\rpy_{em}$, $m=1,\dots,h_e$ be the outgoing ends of $X_e$.
For any $a\in\Om^1(Z'_e)$ define $R_ea\in\R^{h_e}$ by
\[(R_ea)_m=\int_{\{0\}\times Y_{em}}*a.\]
Set $n_e=g_e+h_e$ and
\[L_ea=(J_ea,R_ea)\in\R^{n_e}.\]
For any $a\in\Om^1(\xt)$ let $La\in V=\oplus_e\R^{n_e}$ be the element
with components $L_ea$.

For any tubular end $\R_+\times P$ of $X$ let $t\col \R_+\times P\to\R_+$
be the projection. Choose a small $\si>0$ and for each $e$ a smooth function
$\ka_e\col X_e\to\R$ such that
\[\ka_e=
\begin{cases}
\si t&\text{on incoming ends,}\\
-\si t&\text{on outgoing and neutral ends.}
\end{cases}\]
Let $\pxt\subset X$ be as in Subsection~\ref{intro:comp-neck} and
let $\ka=\ka_T\col \xt\to\R$ be a
smooth function such that $\ka_T-\ka_e$ is constant on $X_e\cap\pxt$ for
each $e$. (Such a function exists because $\ga$ is simply-connected.)
This determines $\ka_T$ up to an additive constant. 

Fix a point $x_e\in X_e$ and define a norm $\|\cdot\|_T$ on $V$
by
\[\|v\|_T=\sum_e\exp(\ka_T(x_e))\|v_e\|,\]
where $\|\cdot\|$ is the Euclidean norm on $R^{n_e}$
and $\{v_e\}$ the components of $v$.

Let $\cD$ denote the operator $-d^*+d^+$ on $\xt$.

\begin{lemma}\label{neck-inequal}
There is a constant $C$ such that for every $r$--tuple $T$ with $\min_jT_j$
sufficiently large and every
$L^{p,\ka}_1$ $1$--form $a$ on $\xt$ we have
\[\|a\|_\llpk\le C(\|\cD a\|_\lpk+\|La\|_T).\]
\end{lemma}
Note that adding a constant to $\ka$ rescales all norms in the above
inequality by the same factor.

\proof Let $\tau$ be a function on $X$ which is equal to $2T_j$ on the ends
$\R_+\times(\pm Y_j)$ for each $j$.
Choose smooth functions $f_1,f_2\col \R\to\R$
such that $(f_1(t))^2+(f_2(1-t))^2=1$ for all $t$, and $f_k(t)=1$ for
$t\le\frac13$, $k=1,2$.
For each $e$ define $\beta_e\col X_e\to\R$ by
\[\beta_e=\begin{cases}
f_1(t/\tau)&\text{on outgoing ends,}\\
f_2(t/\tau)&\text{on incoming ends,}\\
1&\text{elsewhere.}
\end{cases}\]
Let $\bbe$ denote the smooth function on $\xt$ which agrees with
$\beta_e$ on $X_e\cap\pxt$ and is zero elsewhere.

In the following, $C,C_1,C_2,\dots$ will be constants that are
independent of $T$. Set $\tmin=\min_jT_j$. Assume $\tmin\ge1$.

Note that $|\nabla\beta_e|\le C_1\tmin\inv$ everywhere, and similarly for $\bbe$.
Therefore
\[\|\beta_ea\|_\llpke\le C_2\|a\|_\llpke\]
for $1$--forms $a$ on $X_e$.

Let $\cD_e$ denote the operator $-d^*+d^+$ on $X_e$.
By Proposition~\ref{kerD} the Fredholm operator
\[\cD_e\oplus L_e\col \llpke\to\lpkj\oplus\R^{n_e}\]
is injective, hence it has a bounded left inverse $P_e$,
\[P_e(\cD_e\oplus L_e)=\id.\]
If $a$ is a $1$--form on $\xt$ and $v\in V$ set
\[\bbe(a,v)=(\bbe a,v_e).\]
Here we regard $\bbe a$ as a $1$--form on $X_e$. Define
\[P=\sum_e\beta_eP_e\bbe\col \lpk\oplus V\to\llpk.\]
If we use the norm $\|\cdot\|_T$ on $V$ then $\|P\|\le C_3$. Now
\begin{align*}
P(\cD\oplus L)a&=\sum_e\beta_eP_e(\bbe\cD a,L_ea)\\
&=\sum_e\beta_eP_e(\cD_e\bbe a+[\bbe,\cD]a,L_e\bbe a)\\
&=\sum_e(\beta_e\bbe a+\beta_eP_e([\bbe,\cD]a,0))\\
&=a+Ea,
\end{align*}
where
\[\|Ea\|_\llpk\le C_4\tmin\inv\|a\|_\lpk.\]
Therefore,
\[\|P(\cD\oplus L)-I\|\le C_4\tmin\inv,\]
so if $\tmin>C_4$ then $z=P(\cD\oplus L)$ will
be invertible, with
\[\|z\inv\|\le(1-\|z-I\|)\inv.\]
In that case we can define a left inverse of $\cD\oplus L$ by
\[Q=(P(\cD\oplus L))\inv P.\]
If $\tmin\ge2C_4$ say, then $\|Q\|\le2C_3$, whence for any $a\in\llpk$ we have
\[\|a\|_\llpk=\|Q(\cD a,La)\|_\llpk\le C(\|\cD
a\|_\lpk+\|La\|_T).\eqno{\qed}\]

\begin{lemma}\label{econv}
Let $e$ be a node of $\ga$ and for $n=1,2,\dots$ let $T(n)$ be an $r$--tuple
and $a_n$ an $\llpkn$ $1$--form on $\xtn$, where $\katn=\ka_{T(n)}$. Suppose
\begin{enumerate}
\item[\rm(i)]$\Si(X,\ga,\ortn)=0$,
\item[\rm(ii)]$\min_jT_j(n)\to\infty$ as $n\to\infty$,
\item[\rm(iii)]There is a constant $C'<\infty$ such that 
\[\katn(x_{e'})\le\katn(x_e)+C'\]
for all nodes $e'$ of $\ga$ and all $n$,
\item[\rm(iv)]$\sup_n\|d^+a_n\|_\infty<\infty$.
\end{enumerate}
Then there are smooth $u_n\col \xtn\to\U1$ such that the sequence
$b_n=a_n-u_n\inv du_n$ is bounded in $L^p_1$ over compact subsets of $X_e$,
and $b_n\in\llpkn$ and $d^*b_n=0$ for every $n$.
\end{lemma}
Note that (iii) implies that $e$ must be a source of $(\ga,\ortn)$.

\begin{proof}
Without loss of generality we may assume that $\katn(x_e)=1$ for all
$n$, in which case
\[\sup_n\|1\|_{\lpkn}<\infty.\]
By Lemmas~\ref{g-char} and \ref{gk-inj} we can find smooth
$v_n\col \xtn\to\U1$ such that
\[\sup_n\|J_{e'}(a_n-v_n\inv dv_n)\|<\infty\]
for every node $e'$, where $\|\cdot\|$ is the Euclidean norm.
(Compare the proof of Proposition~\ref{prop:loc-comp}.) Moreover, we can take
$v_n$ translationary invariant over each end of $\xtn$.
Proposition~\ref{laplace-iso2} then provides smooth
$\xi_n\in\llw p{\ka_n}2(X;i\R)$ such that
\[b_n=a_n-v_n\inv dv_n-d\xi_n\in\llpkn\]
satisfies
\[d^*b_n=0,\qquad \int_{\{0\}\times Y'_j}*b_n=0\]
for $j=1,\dots,r'-1$. Stokes' theorem shows that the integral vanishes
for $j=r'$ as well, and since $\ga$ is simply-connected we obtain,
for $j=1,\dots,r$,
\[\int_{\{t\}\times Y_j}*b_n=0\qquad\text{for $|t|\le T_j(n)$.}\]
In particular, $R_{e'}b_n=0$ for all nodes $e'$ of $\ga$. 

Set $u_n=v_n\exp(\xi_n)$, so that $b_n=a_n-u_n\inv du_n$.
By Lemma~\ref{neck-inequal} we have
\begin{align*}
\|b_n\|_{\llpkn}&\le C\left(\|d^+b_n\|_{\lpkn}
+\sum_{e'}\exp(\katn(x_{e'}))\|J_{e'}(b_n)\|\right)\\
&\le C\left(\|d^+a_n\|_{\lpkn}
+\exp(1+C')\sum_{e'}\|J_{e'}(a_n-v_n\inv dv_n)\|\right),
\end{align*}
which is bounded as $n\to\infty$.
\end{proof}

\begin{prop}\label{prop:xe-conv}
Suppose $\ga$ is simply-connected and that Condition~(C) holds for $(X,\ga)$.
For $n=1,2,\dots$ let
$[S_n]\in M(\xtn;\vec\al'_n;\mu_n;\vec\prt_n;\vec\prt'_n)$,
where $\min_jT_j(n)\to\infty$. Then there exist smooth maps $w_n\col X\to\U1$
such that if $k$ is any positive integer with
\[\sup_{j,j',n}(\|\mu_n\|_{C^k}+
\|\prt_{j,n}\|_{C^k}+\|\prt'_{j',n}\|_{C^k})<\infty\]
then the sequence $w_n(S_n)$ is
bounded in $L^{p'}_{k+1}$ over compact subsets of $X$ for every $p'\ge1$.
\end{prop}

\begin{proof}
Consider the set-up in the beginning of this subsection where now
$p>4$ is the exponent used in defining configuration spaces, and
$\ortn$ is an orientation of $\ga$ for which (C) is fulfilled.
By passing to a subsequence we can arrange that
$\katn(x_e)-\katn(x_{e'})$ converges to a point $\ell(e,e')\in[-\infty,\infty]$
for each pair of nodes $e,e'$ of $\ga$. Define an equivalence relation $\sim$
on the set $\mathcal N$ of nodes of $\ga$ by declaring that
$e\sim e'$ if and only if $\ell(e,e')$ is finite. Then we have
a linear ordering on $\mathcal N/\sim$ such that
$[e]\ge[e']$ if and only if $\ell(e,e')>-\infty$. Here $[e]$ denotes
the equivalence class of $e$. 

Choose $e$ such that $[e]$ is the maximum with respect to this linear ordering.
Let $S_n=(A\rr+a_n,\Phi_n)$.
Then all the hypotheses of Lemma~\ref{econv} are satisfied. If $u_n$ is as in
that lemma then,
as in the proof of Proposition~\ref{prop:loc-comp},
$u_n(S_n)=(A\rr+b_n,u_n\Phi_n)$ will be bounded
in $L^{p'}_{k+1}$ over compact subsets of $X_e$ for every $p'\ge1$.

For any $r$--tuple $T$ let $\wt$ be the result of gluing ends of
$X\setminus X_e$ according to the graph $\ga\setminus e$ and (the relevant 
part of) the vector $T$. To simplify notation let us assume that the outgoing
ends of $X_e$ are $\R_+\times(-Y_j)$, $j=1,\dots,r_1$. Then $\rpy_j$ is an
end of $\wtn$ for $j=1,\dots,r_1$. Let $b'_n$ be the $1$--form on $\wtn$
which away from the ends $\rpy_j$, $1\le j\le r_1$ agrees with $b_n$, and
on each of these ends is defined by cutting off $b_n$:
\[b'_n(t,y)=\begin{cases}
f_1(t-2T_j(n)+1)\cdot b_n(t,y),& 0\le t\le 2T_j(n),\\
0,& t\ge2T_j(n).
\end{cases}\]
Here $f_1$ is as in the proof of Lemma~\ref{neck-inequal}.
Then $\sup_n\|d^+b'_n\|_\infty<\infty$. After choosing
an orientation of $\ga\setminus e$ for which (C) holds we can apply 
Lemma~\ref{econv} to each component of $\ga\setminus e$, with $b'_n$
in place of $a_n$. Repeating this process proves the proposition.
\end{proof}

\begin{cor}\label{lemma:xtn-cd-bound}
If (B2) holds then for every constant $C_0<\infty$ there exists a constant
$C_1<\infty$ such that for
any element $[S]$ of a moduli space
$M(\xt;\vec\al';\mu;\vec\prt;\vec\prt')$ where $\min_jT_j>C_1$
and $\mu,\prt_j,\prt'_j$ all have $C^1$--norm $<C_0$
one has
\[\left|
\sum_{j=1}^r\lla_j
\left(\cd(S|_{\{-T_j\}\times Y_j})
-\cd(S|_{\{T_j\}\times Y_j})\right)
+\sum_{j=1}^{r'}\lla'_j
\cd(S|_{\{0\}\times Y'_j})\right|
< C_1.\]
\end{cor}

The next proposition, which is essentially a corollary of
Proposition~\ref{prop:xe-conv}, exploits the fact that Condition~(C) is
preserved under certain natural extensions of $(X,\ga)$.

\begin{prop}\label{prop:local-energy-bound2} 
Suppose $\ga$ is simply-connected and (C) holds for $(X,\ga)$.
Then for every constant $C_0<\infty$ there is a constant
$C_1<\infty$ such that if $S=(A,\Phi)$ 
represents an element of
a moduli space $M(\xt;\vec\al';\mu;\vec\prt;\vec\prt')$ where $\min_jT_j>C_1$
and $\mu,\prt_j,\prt'_j$ all have $L^\infty$ norm $<C_0$ then
\begin{align*}
&\|\nabla\cd_S\|_{L^2((t-1,t+1)\times Y_j)}
<C_1\quad\text{for $|t|\le T_j-1$},\\
&\|\nabla\cd_S\|_{L^2((t-1,t+1)\times Y'_j)}<C_1\quad\text{for $t\ge1$}.
\end{align*}
\end{prop}

\begin{proof}
Given an edge $v$ of $\ga$ corresponding to a pair of ends
$\R_+\times(\pm Y_j)$ of $X$, we can form a new pair $(X\sj,\ga\sj)$ where
$X\sj=X\coprod(\ry_j)$, and $\ga\sj$ is obtained from $\ga$ by splitting $v$
into two edges with a common end-point representing the component $\ry_j$
of $X\sj$.

Similarly, if $e$ is a node of $\ga$ and $\rpy'_j$ an end of $X_e$
then we can form a new pair $(X\Sj,\ga\Sj)$ where
$X\Sj=X\coprod(\ry'_j)$, and $\ga\Sj$ is obtained from $\ga$ by adding one
node $e'_j$ representing the component $\rpy'_j$ of $X\Sj$ and one edge
joining $e$ and $e'_j$.

One easily shows, by induction on the number of nodes of $\ga$, that if (C)
holds for $(X,\ga)$ then (C) also holds for each of the new pairs
$(X\sj,\ga\sj)$ and $(X\Sj,\ga\Sj)$. 
Given this observation, the proposition
is a simple consequence of Proposition~\ref{prop:xe-conv}.
\end{proof}

\section{Exponential decay}\label{sec:exp-decay}

In this section we will prove exponential decay results for genuine
monopoles over
half-cylinders $\rpy$ and long bands $[-T,T]\times Y$. The overall scheme of 
proof will be the same as that for instantons in \cite{D5}, and
Subsections~\ref{subsec:diff-ineq} and \ref{subsec:decay-mon}
follow \cite{D5} quite closely. On the other hand, the proof of the
main result of Subsection~\ref{subsec:est-band}, Proposition~\ref{gauge3},
is special to monopoles (and is new, as far as we know).

Throughout this section $Y$ will be a closed, connected Riemannian $\spc$
$3$--manifold, and $\eta\in\Om^2(Y)$ closed. We will study exponential
decay towards a non-degenerate critical point $\al$ of $\cd=\cd_\eta$.
We make no non-degeneracy assumptions on any other (gauge equivalence classes
of) critical points, and we do not assume that (O1) holds,
except implicitly in Proposition~\ref{half-cyl-conv}.
All monopoles will be genuine (ie $\prt=0$).

Previously, Nicolaescu~\cite{nico1} has proved a slightly weaker
exponential decay result in the case $\eta=0$, using different techniques.
Other results (with less complete proofs) were obtained
in \cite{MST} and \cite{Marcolli-Wang1}.

\subsection{A differential inequality}\label{subsec:diff-ineq}

We begin by presenting an argument from \cite{D5} in a more abstract setting,
so that it applies equally well to the Chern--Simons and the
Chern--Simons--Dirac
functionals.

Let $E$ be a real Banach space with norm $\|\cdot\|$ and $E'$ a real
Hilbert space with inner product $\la\cdot,\cdot\ra$.
Let $E\to E'$ be an injective,
bounded operator with dense image. We will identify $E$ as a vector space
with its image in $E'$. Set $|x|=\la x,x\ra^{1/2}$ for $x\in E'$.

Let $U\subset E$ be an open set containing $0$ and
\[f\col U\to\R,\qquad g\col U\to E'\]
smooth maps satisfying $f(0)=0$, $g(0)=0$, and
\[Df(x)y=\la g(x),y\ra\]
for all $x\in U$, $y\in E$. Here $Df(x)\col E\to\R$ is the derivative of $f$
at $x$. Suppose $H=Dg(0)\col E\to E'$ is an isomorphism
(of topological vector spaces). Note that $H$ can
be thought of as a symmetric operator in $E'$. Suppose $E$ contains a 
countable set $\{e_j\}$ of eigenvectors for $H$ which forms an orthonormal
basis for $E'$. Suppose $\si,\lla$ are real numbers satisfying
$0\le\lla<\si$ and such that $H$ has no positive eigenvalue less than $\si$.

\begin{lemma}\label{grad-estimate}
In the above situation there is a constant $C>0$ such that for every
$x\in U$ with $\|x\|\le C\inv$ one has
\begin{align*}
2\si f(x)&\le|g(x)|^2+C|g(x)|^3,\\
2\lla f(x)&\le|g(x)|^2.
\end{align*}
\end{lemma}

\begin{proof}
It clearly suffices to establish the first inequality for some $C$.
By Taylor's formula \cite[8.14.3]{Dieud} there is a $C_1>0$ such that
for all $x\in U$ with $\|x\|\le C_1\inv$ one has
\begin{align*}
|f(x)-\frac12\la Hx,x\ra|&\le C_1\|x\|^3,\\
|g(x)-Hx|&\le C_1\|x\|^2.
\end{align*}
Let $He_j=\lla_je_j$, and set $x_j=\la x,e_j\ra e_j$. Then
\[\si\la Hx,x\ra=\si\sum_j\lla_j|x_j|^2\le\sum_j\lla_j^2|x_j|^2=|Hx|^2.\]
By assumption, there is a $C_2>0$ such that
\[\|x\|\le C_2|Hx|\]
for all $x\in E$. Putting the above inequalities together we get, for
$r=\|x\|\le C_1\inv$,
\begin{align*}
|Hx|&\le|g(x)|+C_1r^2\\
&\le|g(x)|+rC_1C_2|Hx|.
\end{align*}
If $r<(C_1C_2)\inv$ this gives
\[|Hx|\le(1-rC_1C_2)\inv|g(x)|,\]
hence
\begin{align*}
2\si f(x)&\le\si\la Hx,x\ra+2\si C_1r^3\\
&\le|Hx|^2+2\si rC_1C_2^2|Hx|^2\\
&\le\frac{1+2\si rC_1C_2^2}{(1-rC_1C_2)^2}|g(x)|^2\\
&\le(1+C_3r)|g(x)|^2\\
&\le|g(x)|^2+C_4|g(x)|^3
\end{align*}
for some constants $C_3,C_4$.
\end{proof}

We will now apply this to the Chern--Simons--Dirac functional.
Let $\al=(B_0,\Psi_0)$ be a non-degenerate critical point of $\cd$.
Set $K_\al=\ker(\ci^*_\al)\subset\Ga(i\La^1\oplus\bs)$
and $\ti H_\al=H_\al|_{K_\al}\col K_\al\to K_\al$. Note that any eigenvalue
of $\ti H_\al$ is
also an eigenvalue of the self-adjoint elliptic operator
\[\left(
\begin{array}{cc}
0 & \ci^*_\al \\
\ci_\al & H_\al
\end{array}
\right)\]
over $Y$ acting on sections of $i\La^0\oplus i\La^1\oplus\bs$.
Let $\lla^\pm$ be positive real numbers such that $\ti H_\al$ has no
eigenvalue in the interval $[-\lla^-,\lla^+]$.

In the following lemma, Sobolev norms of sections of the spinor bundle
$\bs_Y$ over $Y$ will be taken with respect to $B_0$ and some fixed connection
in the tangent bundle $TY$. This means that the same constant $\eps$ 
will work if $\al$ is replaced with some monopole gauge equivalent to $\al$.

\begin{lemma}\label{csd-est1}
In the above situation there exists an $\eps>0$ such that if $S$ is
any smooth monopole over the band $(-1,1)\times Y$
satisfying
$\|S_0-\al\|_{L^2_1(Y)}\le\eps$ then
\[\pm2\lla^\pm(\cd(S_0)-\cd(\al))\le-\prtl_t|_0\cd(S_t).\]
\end{lemma}

\begin{proof}
Choose a smooth $u\col (-1,1)\times Y\to\U1$ such that $u(S)$ is in temporal
gauge. Then
\[\prtl_t\cd(S_t)=\prtl_t\cd(u_t(S_t))=-\|\nabla\cd(u_t(S_t))\|_2^2
=-\|\nabla\cd(S_t)\|_2^2.\]
If $\eps>0$ is sufficiently small then by the local slice theorem
we can find a smooth $v\col Y\to\U1$ which is $L^2_2$ close
to $1$ and such that $\ci_\al^*(v(S_0)-\al)=0$. 
We now apply Lemma~\ref{grad-estimate} with $E$ the kernel of $\ci^*_\al$ in
$L^2_1$, $E'$ the kernel of $\ci^*_\al$ in $L^2$, and
$f(x)=\pm(\cd(\al+x)-\cd(\al))$.
The assumption that $\al$ be non-degenerate means
that $H=\ti H_\al\col E\to E'$ is an isomorphism, so the lemma follows.
\end{proof}

\subsection{Estimates over $[0,T]\times Y$}\label{subsec:est-band}

Let $\al$ be a non-degenerate critical point of $\cd$
and $\ual=(B,\Psi)$ the monopole over $\ry$ that $\al$
defines.
Throughout this subsection, the same convention for Sobolev norms of sections
of $\bs_Y$ will apply as in Lemma~\ref{csd-est1}. For Sobolev norms of sections
of the spinor bundles over (open subsets of)
$\ry$ we will use the connection $B$.

Throughout this subsection
$S=(A,\Phi)$ will be a monopole
over a band $\BB=[0,T]\times Y$ where $T\ge1$. Set $s=(a,\phi)=S-\ual$ and
\begin{equation}\label{del-nu}
\begin{gathered}
\del=\|s\|_{L^2_2(\BB)},\\
\nu^2=\|\nabla\cd_S\|^2_{L^2(\BB)}=\cd(S_0)-\cd(S_T).
\end{gathered}
\end{equation}

The main result of this subsection is Proposition~\ref{gauge3}, which asserts
in particular that if $\del$ is sufficiently small then $S$ is gauge 
equivalent to a configuration $\ti S$ which is in Coulomb gauge
with respect to $\ual$ and satisfies $\|\ti S-\ual\|_{L^2_1(\BB)}\le\const\nu$.

We will assume $\del\le1$.
Let $\adt$ denote the contraction of $a$ with the
vector field $\prtl_1=\frac\prtl{\prtl t}$.
Quantities referred to as constants or denoted ``const''
may depend on $Y,\eta,[\al],T$ but not on $S$. Note that
\begin{equation}\label{nubound}
\nu\le\const(\|s\|_{1,2}+\|s\|_{1,2}^2)\le\text{const},
\end{equation}
the last inequality because $\del\le1$.

For real numbers $t$ set
\[i_t\col Y\to\ry,\quad y\mapsto(t,y).\]
If $\om$ is any differential form over $\BB$ set $\om_t=i^*_t\om$,
$0\le t\le T$.
Similar notation will be used for connections and spinors over $\BB$.

\begin{lemma}\label{prtlphi}
There is a constant $C_0>0$ such that
\[\|\prtl_1\phi\|_2\le C_0(\nu+\|\adt\|_3).\]
\end{lemma}

\begin{proof}
We have
\[\prtl_1\phi=\prtl_1\Phi=\nabla^A_1\Phi-\adt\Phi,\]
where $\nabla^A_1$ is the covariant derivative with respect to $A$ in the 
direction of the vector field $\prtl_1=\frac\prtl{\prtl t}$.
Now $|\nabla^A_1\Phi|$ depends only on the gauge equivalence class of
$S=(A,\Phi)$, and if $A$ is in 
temporal gauge (ie if $\adt=0$) then
$(\nabla^A_1\Phi)_t=\prtl_{A_t}\Phi_t$. 
The lemma now follows because $\del\le1$.
\end{proof}

\begin{lemma}\label{gauge1}
There is a constant $C_1>0$ such that if $\del$ is sufficiently
small then the following hold:
\begin{enumerate}
\item[\rm(i)]$\|\phi\|_{1,2}\le C_1(\|\ci^*_\ual s\|_2+\|\adt\|_{1,2}+\nu)$,
\item[\rm(ii)]There is a smooth $\check f\col \BB\to i\R$ such that 
$\check s=(\check a,\check\phi)=\exp(\check f)(S)-\ual$ satisfies
\[\|\check s_t\|_{1,2}\le C_1\|\nabla\cd_{S_t}\|_2,\quad0\le t\le T,\]
\item[\rm(iii)]$\|da\|_2\le C_1\nu$,
\end{enumerate}
where in (i) and (iii) all norms are taken over $\BB$.
\end{lemma}

\begin{proof}
The proof will use an elliptic inequality over
$Y$, the local slice theorem for $Y$,
and the gradient flow description of the Seiberg--Witten equations over $\ry$.

(i)\qua Since $\al$ is non-degenerate we have
\[\|z\|_{1,2}\le\const\|(\ci^*_\al+H_\al)z\|_2\]
for $L^2_1$ sections $z$ of $(i\La\oplus\bs)_Y$. Recall that
\[\nabla\cd_{\al+z}=H_\al z+z\otimes z\]
where $z\otimes z$ represents a pointwise quadratic function of $z$.
Furthermore, $\|z\otimes z\|_2\le\const\|z\|_{1,2}^2$.
If $\|z\|_{1,2}$ is sufficiently small then we can rearrange to get
\begin{equation}\label{z12}
\|z\|_{1,2}\le\const(\|\ci^*_\al z\|_2+\|\nabla\cd_{\al+z}\|_2).
\end{equation}
By the Sobolev embedding theorem we have
\[\|s_t\|_{L^2_1(Y)}\le\const\|s\|_{L^2_2(\BB)},\quad t\in[0,T],\]
for some constant independent of $t$, so we can apply inequality~\Ref{z12}
with $z=s_t$ when $\del$ is sufficiently small. Because
\[(\ci^*_\ual s-\prtl_1\adt)_t=\ci^*_\al s_t\]
we then obtain
\begin{equation*}
\int_0^T\|s_t\|^2_{L^2_1(Y)}dt\le
\const(\|\ci^*_\ual s\|^2_{L^2(\BB)}+\|\prtl_1\adt\|^2_{L^2(\BB)}+\nu^2).
\end{equation*}
This together with Lemma~\ref{prtlphi} establishes (i).

(ii)\qua Choose a base-point $y_0\in Y$.
By the local slice theorem there is a constant $C$ such that
if $\del$ is sufficiently small then for each $t\in[0,T]$ there is 
a unique smooth $\check f_t\col Y\to i\R$ such that
\begin{itemize}
\item $\|\check f_t\|_{2,2}\le C\del$,
\item $\check f_t(y_0)=0$ if $\al$ is reducible,
\item $\check s_t=\exp(\check f_t)(S_t)-\al$ satisfies $\ci^*_\al\check s_t=0$.
\end{itemize}
It is not hard to see that the function
$\check f\col \BB\to i\R$ given by $\check f(t,y)=\check f_t(y)$
is smooth. Moreover,
$\|\check s_t\|_{1,2}\le\const\|s_t\|_{1,2}$. Part~(ii) then follows by
taking $z=\check s_t$ in \Ref{z12}.

(iii)\qua Choose a smooth $u\col \BB\to\U1$ such that $u(S)$ is in
temporal gauge, and set $(\uline a,\uline\phi)=u(S)-\ual$. Then
\[da=d\uline a=dt\wedge\prtl_1\uline a+d_Y\uline a
=-dt\wedge\nabla_1\cd_{S_t}+d_Y\check a_t,\]
where $\nabla_1\cd$ is the first component of $\nabla\cd$.
This yields the desired estimate on $da$.
\end{proof}

\begin{lemma}\label{b-ineq}
Let $\{v_1,\dots,v_{b_1(Y)}\}$ be a family of closed $2$--forms on $Y$ which
represents a basis for $H^2(Y;\R)$. Then there is a constant $C$ such that
\begin{equation}\label{bvest}
\|b\|_{L^2_1(\BB)}\le C\left(
\|(d^*+d)b\|_{L^2(\BB)}+
\|(*b)|_{\prtl\BB}\|_{L^2_{1/2}(\prtl\BB)}+
\sum_j\left|\int_\BB dt\wedge\pi^*v_j\wedge b\right|\right)
\end{equation}
for all $L^2_1$ $1$--forms  $b$ on $\BB$, where $\pi\col \BB\to Y$ is the
projection.
\end{lemma}

\begin{proof}
Let $K$ denote the kernel of the operator
\[\Om^1_{\BB}\to\Om^0_{\BB}\oplus\Om^2_{\BB}\oplus\Om^0_{\prtl\BB},\quad
b\mapsto(d^*b,db,*b|_{\prtl\BB}).\]
Then we have an isomorphism
\[\rho\col K\oset\approx\to H^1(Y;\R),\quad b\mapsto[b_0].\]
For on the one hand, an application of Stokes' theorem shows that $\rho$
is injective. On the other hand, any $c\in H^1(Y;\R)$ can be represented by
an harmonic $1$--form $\om$, and $\pi^*\om$ lies in $K$,
hence $\rho$ is surjective.

It follows that every element of $K$ is of the form $\pi^*(\om)$.
Now apply Proposition~\ref{fred} and Lemma~\ref{gen-ineq}.
\end{proof}

\begin{lemma}\label{gauge2}
There is a smooth map $\hat f\col \BB\to i\R$, unique up to an additive constant,
such that $\hat a=a-d\hat f$ satisfies
\[d^*\hat a=0,\quad(*\hat a)|_{\prtl\BB}=0.\]
Given any such $\hat f$,
if we set $\hat s=(\hat a,\hat \phi)=\exp(\hat f)(S)-\ual$ then
\[\|\hat a\|_{L^2_1(\BB)}\le C_2\nu,\quad\|\hat s\|_{L^2_2(\BB)}\le C_2\del\]
for some constant $C_2>0$.
\end{lemma}

\begin{proof}
The first sentence of the lemma is just the solution to the
Neumann problem.
If we fix $x_0\in\BB$ then there is a unique $\hat f$ as in the lemma
such that $\hat f(x_0)=0$, and we have
$\|\hat f\|_{3,2}\le\const\|a\|_{2,2}$. 
Writing
\[\hat\phi=\exp(\hat f)\Phi-\Psi=(\exp(\hat f)-1)\Phi+\phi\]
and recalling that, for functions on $\BB$, multiplication is a continuous map
$L^2_3\times L^2_k\to L^2_k$ for $0\le k\le3$, we get
\begin{align*}
\|\hat\phi\|_{2,2}&\le C\|\exp(\hat f)-1\|_{3,2}\|\Phi\|_{2,2}+
\|\phi\|_{2,2}\\
&\le C'\|\hat f\|_{3,2}\exp(C''\|\hat f\|_{3,2})+\|\phi\|_{2,2}\\
&\le C'''\|s\|_{2,2}
\end{align*}
for some constants $C,\dots,C'''$, since we assume $\del\le1$.
There is clearly a similar $L^2_2$ bound on $\hat a$, so this establishes
the $L^2_2$ bound on $\hat s$.

We now turn to the $L^2_1$ bound on $\hat a$.
Let $\check a$ be as in Lemma~\ref{gauge1}. 
Since $\hat a-\check a$ is exact we have
\[\left|\int_Yv\wedge\hat a_t\right|=\left|\int_Yv\wedge\check a_t\right|
\le\const\|v\|_2\|\check a_t\|_2\]
for any closed $v\in\Om^2_Y$.
Now take $b=\hat a$ in Lemma~\ref{b-ineq} and use Lemma~\ref{gauge1}, 
remembering that $d\hat a=da$.
\end{proof}

\begin{defn}
{\rm For any smooth $h\col Y\to i\R$ define $\uh,P(h)\col \BB\to i\R$ by
$\uh(t,y)=h(y)$ and
\[P(h)=\Delta\uh+i\la i\Psi,\exp(\uh)\Phi\ra,\]
where $\Delta=d^*d$ is the Laplacian over $\ry$.
Let $P_t(h)$ be the restriction of $P(h)$ to $\{t\}\times Y$.}
\end{defn}

Note that
\[\ci^*_\ual(\exp(\uh)(S)-\ual)=-d^*a+P(h).\]

\begin{lemma}\label{ph}
If $\al$ is irreducible then the following hold:
\begin{enumerate}
\item[\rm(i)]There is a $C_3>0$ such that if $\del$
is sufficiently small then there exists a unique smooth $h\col Y\to i\R$ satisfying
$\|h\|_{3,2}\le C_3\del$ and $P_0(h)=0$.
\item[\rm(ii)]If $h\col Y\to i\R$ is any smooth function satisfying $P_0(h)=0$ 
then
\[\|P(h)\|_{L^2(\BB)}\le\const(\nu+\|a'\|_{L^3(\BB)}).\]
\end{enumerate}
\end{lemma}

\begin{proof}
(i)\qua We will apply Proposition~\ref{quant-inv2} (ie the inverse
function theorem) to the smooth map
\[P_0\col L^2_3\to L^2_1,\quad h\mapsto\Delta_Yh+i\la i\Psi_0,\exp(h)\Phi_0\ra.\]
The first two derivatives of this map are
\begin{align*}
DP_0(h)k&=\Delta_Yk+k\la\Psi_0,\exp(h)\Phi_0\ra,\\
D^2P_0(h)(k,\ell)&=ik\ell\la i\Psi_0,\exp(h)\Phi_0\ra.
\end{align*}
The assumption $\del\le1$ gives
\[\|D^2P_0(h)\|\le\const(1+\|\nabla h\|_3).\]
Set $L=DP_0(0)$. Then
\[(L-\Delta_Y-|\Psi_0|^2)k=k\la\Psi_0,\phi_0\ra,\]
hence
\[\|L-\Delta_Y-|\Psi_0|^2\|\le\const\del.\]
Thus if $\del$ is sufficiently small then $L$ is invertible and
\[\|L\inv\|\le\|(\Delta_Y+|\Psi_0|^2)\inv\|+1.\] 
Furthermore, we have $P_0(0)=i\la i\Psi_0,\phi_0\ra$, so
\[\|P_0(0)\|_{1,2}\le\const\del.\]
By Proposition~\ref{quant-inv2}~(i) there exists a constant $C>0$ such that
if $\del$ is sufficiently small then there is a unique $h\in L^2_3$ such that
$\|h\|_{3,2}\le C$ and $P_0(h)=0$ (which implies that $h$ is smooth).
Proposition~\ref{quant-inv2}~(ii) then yields
\[\|h\|_{3,2}\le\const\|P_0(0)\|_{1,2}\le\const\del.\]

(ii)\qua Setting $Q=P(h)$ we have, for $0\le t\le T$,
\[\int_Y|Q(t,y)|^2dy=\int_Y\left|\int_0^t\prtl_1Q(s,y)\,ds\right|^2dy
\le\const\int_{\BB}|\prtl_1Q|^2.\]
Now, $\prtl_1Q=i\la i\Psi,\exp(\uh)\prtl_1\Phi\ra$, hence
\[\|\prtl_1Q\|_2\le\const\|\prtl_1\Phi\|_2\le\const(\nu+\|a'\|_3)\]
by Lemma~\ref{prtlphi}.
\end{proof}

\begin{prop}\label{gauge3}
There is a constant $C_4$ such that if $\del$ is sufficiently small then
there exists a smooth $\ti f\col \BB\to i\R$ such that $\ti s=(\ti a,\ti\phi)
=\exp(\ti f)(S)-\ual$ satisfies
\[\ci^*_\ual\ti s=0,\quad(*\ti a)|_{\prtl\BB}=0,\quad
\|\ti s\|_{L^2_1(\BB)}\le C_4\nu,\quad\|\ti s\|_{L^2_2(\BB)}\le C_4\del,\]
where $\del,\nu$ are as in \Ref{del-nu}.
\end{prop}

This is analogous to Uhlenbeck's theorem \cite[Theorem~1.3]{U1} (with p=2), 
except that we assume a bound on $\del$ rather than on $\nu$.

\proof
To simplify notation we will write $\ci=\ci_\ual$ in this proof.

{\bf Case 1: $\al$ reducible}\qua In that case the operator $\ci^*$
is given by $\ci^*(b,\psi)=-d^*b$. Let $\ti f$ be the $\hat f$ provided by 
Lemma~\ref{gauge2}. Then apply Lemma~\ref{gauge1}~(ii),
taking the $S$ of that lemma to be the present $\exp(\ti f)(S)$.

{\bf Case 2: $\al$ irreducible}\qua Let $\hat f,\hat S$ etc be as in
Lemma~\ref{gauge2}. Choose $h\col Y\to i\R$ such that the conclusions of
Lemma~\ref{ph}~(i) holds with the $S$ of that lemma taken to be the present
$\hat S$. Set $\acute S=(\acute A,\acute\Phi)=\exp(\uh)(\hat S)$ and
$\acute s=(\acute a,\acute\phi)=\acute S-\ual$.
By Lemmas~\ref{ph} and \ref{gauge1}~(ii) we have
\[\|\ci^*\acute s\|_2\le\const\nu,\quad\|\acute s\|_{2,2}\le\const\del,
\quad\|\acute\phi\|_{1,2}\le\const\nu.\]
Since $-d^*\acute a=\ci^*\acute s-i\la i\Psi,\acute\phi\ra$ we also get
\[\|d^*\acute a\|_2\le\const\nu.\]
Applying Lemma~\ref{b-ineq} as in the proof of Lemma~\ref{gauge2} we see that
\[\|\acute a\|_{1,2}\le\const\nu.\]

It now only remains to make a small perturbation to $\acute S$ so as to
fulfil the Coulomb gauge condition. To this end we invoke
the local slice theorem for $\BB$. This says that
there is a $C>0$ such that if $\del$ is
sufficiently small then there exists a unique smooth $f\col \BB\to i\R$ such that
setting $\ti s=(\ti a,\ti\phi)=\exp(f)(\acute S)-\ual$ one has
\[\|f\|_{3,2}\le C\del,\quad\ci^*\ti s=0,\quad
*\ti a|_{\prtl\BB}=0.\]

We will now estimate first $\|f\|_{2,2}$, then $\|\ti s\|_{1,2}$ in terms of
$\nu$. First note that $*\acute a|\pb=*\hat a|\pb=0$, and
\[\ti a=\acute a-df,\qquad \ti\phi=\exp(f)\acute\Phi-\Psi\]
by definition. Write the imaginary part of $\exp(f)$ as $f+f^3u$.
Then $f$ satisfies the equations $(\prtl_tf)|\pb=0$ and
\begin{align*}
0&=-d^*\ti a+i\la i\Psi,\ti\phi\ra_{\R}\\
&=\Delta f-d^*\acute a+i\la i\Psi,\exp(f)(\acute\phi+\Psi)\ra_{\R}\\
&=\Delta f+|\Psi|^2f-d^*\acute a+i\la i\Psi,\exp(f)\acute\phi+f^3u\Psi\ra_{\R}.
\end{align*}
By the Sobolev embedding theorem we have
\[\|f\|_\infty\le\const\|f\|_{3,2}\le\const\del,\]
and we assume $\del\le1$, so $\|u\|_\infty\le\text{const}$. Therefore,
\begin{align*}
\|f\|_{2,2}&\le\const\|\Delta f+|\Psi|^2f\|_2\\
&\le\nu+\const\|f^3\|_2\\
&\le\nu+\const\|f\|^3_{2,2},
\end{align*}
cf Subsection~\ref{boundary-slice} for the first inequality.
If $\del$ is sufficiently small
then we can rearrange to get $\|f\|_{2,2}\le\const\nu$. Consequently,
$\|\ti a\|_{1,2}\le\const\nu$. To estimate $\|\ti\phi\|_{1,2}$ we write
\[\ti\phi=g\Psi+\exp(f)\acute\phi,\]
where $g=\exp(f)-1$. Then $|dg|=|df|$ and $|g|\le\const|f|$. Now
\begin{align*}
\|\ti\phi\|_2&\le\const\|f\|_2+\|\acute\phi\|_2\le\const\nu,\\
\|\nabla\ti\phi\|_2&\le\const(\|g\|_{1,2}+\|df\otimes\acute\phi\|_2+
\|\nabla\acute\phi\|_2)\\
&\le\const(\nu+\|df\|_4\|\acute\phi\|_4)\\
&\le\const(\nu+\|f\|_{2,2}\|\acute\phi\|_{1,2})\\
&\le\const(\nu+\nu^2)\\
&\le\const\nu
\end{align*}
by \Ref{nubound}. Therefore,
$\|\ti\phi\|_{1,2}\le\const\nu$. Thus, the proposition holds with
\[\ti f=\hat f+\uh+f.\eqno{\qed}\]

\begin{prop}\label{int-est}
Let $k$ be a positive integer and $V\Subset\text{int}(\BB)$ an open subset.
Then there are constants $\eps_k,C_{k,V}$, where $\eps_k$ is independent of
$V$, such that if
\[\ci^*_\ual s=0,\quad\|s\|_{L^2_1(\BB)}\le\eps_k\]
then
\[\|s\|_{L^2_k(V)}\le C_{k,V}\|s\|_{L^2_1(\BB)}.\]
\end{prop}

\begin{proof}
The argument in \cite[pages 62--63]{DK} carries over, if one replaces the
operator $d^*+d^+$ with $\ci^*_\ual+D\sw_\ual$, where
$D\sw_\ual$ is the linearization of the monopole map at $\ual$.
Note that $\ci^*_\ual+D\sw_\ual$ is injective over $S^1\times Y$
because $\al$ is non-degenerate, so if $\ga\col \BB\to\R$ is a smooth function
supported in $\text{int}(\BB)$ then 
\[\|\ga s\|_{k,2}\le C'_k\|(\ci^*_\ual+D\sw_\ual)(\ga s)\|_{k-1,2}\]
for some constant $C'_k$.
\end{proof}

\subsection{Decay of monopoles}\label{subsec:decay-mon}

The two theorems in this subsection are analogues of Propositions~4.3 and
4.4 in \cite{D5}, respectively.

Let $\beta$ be a non-degenerate monopole over $Y$, and $U\subset\cb_Y$ an
$L^2$--closed subset which contains no monopoles except perhaps $[\beta]$.
Choose $\lla^\pm>0$ such that $\ti H_\beta$ has no eigenvalue in the interval
$[-\lla^-,\lla^+]$, and set $\lla=\min(\lla^-,\lla^+)$. Define
\[\BB_t=[t-1,t+1]\times Y.\]

\begin{thm}\label{exp-decay1}
For any $C>0$ there are constants $\eps,C_0,C_1,\dots$ such that the following
holds. Let $S=(A,\Phi)$ be any monopole in temporal
gauge over $(-2,\infty)\times Y$ such that $[S_t]\in U$ for some $t\ge0$. 
Set
\[\bar\nu=\|\nabla\cd_S\|_{L^2((-2,\infty)\times Y)},\quad
\nu(t)=\|\nabla\cd_S\|_{L^2(\BB_t)}.\]
If $\|\Phi\|_\infty\le C$ and $\bar\nu\le\eps$ then there is a smooth monopole
$\al$ over $Y$, gauge equivalent to $\beta$,
such that if $B$ is the connection part of $\ual$ then
for every $t\ge1$ and non-negative integer $k$ one has
\begin{equation}\label{lambda-decay}
\sup_{y\in Y}|\nabla^k_B(S-\ual)|_{(t,y)}\le C_k\sqrt{\nu(0)}e^{-\lla^+t}.
\end{equation}
\end{thm}

\begin{proof}
It follows from the local slice theorem that $\ti\cb_Y\to\cb_Y$ is a 
(topological) principal $H^1(Y;\z)$--bundle. Choose a small open neighbourhood
$V$ of $[\beta]\in\cb_Y$ which is the image of a convex set in $\cc_Y$. We 
define a continuous function $\bar f\col V\to\R$ by
\[\bar f(x)=\cd(\si(x))-\cd(\si([\beta]))\]
where $\si\col V\to\ti\cb_Y$ is any continuous cross-section. It is clear that
$\bar f$ is independent of $\si$.

Given $C>0$, let $S=(A,\Phi)$ be any monopole over $(-2,\infty)\times Y$
such that $\|\Phi\|_\infty\le C$ and $[S_t]\in U$ for some $t\ge0$. If
$\del>0$, and $k$ is any non-negative integer, then provided $\bar\nu$ is
sufficiently small, our local compactness results
(Lemmas~\ref{ptcurv} and \ref{lemma:loc-comp}) imply that for every
$t\ge0$ we can find a smooth $u\col \BB_0\to\U1$ such that
\[\|u(S|_{\BB_t})-\uline\beta\|_{C^k(\BB_0)}<\del.\]
In particular, if $\bar\nu$ is sufficiently small then 
\[f\col \orp\to\R,\quad t\mapsto\bar f([S_t])\]
is a well-defined smooth function. Since $f(t)-\cd(S_t)$ is locally constant,
and $f(t)\to0$ as $t\to\infty$, we have
\[f(t)=\cd(S_t)-L,\]
where $L=\lim_{t\to\infty}\cd(S_t)$. If $\bar\nu$ is sufficiently small then
Lemma~\ref{csd-est1} gives $2\lla^+f\le-f'$, hence
\[0\le f(t)\le e^{-2\lla^+ t}f(0),\quad t\ge0.\]
This yields
\[\nu(t)^2=f(t-1)-f(t+1)\le\const e^{-2\lla^+t}f(0),\quad t\ge1.\]
If $\bar\nu$ is sufficiently small then by Propositions~\ref{gauge3} and
\ref{int-est} we have
\begin{equation}\label{nu-t-bound}
f(t)\le\const\nu(t),\quad
\sup_{y\in Y}|\nabla^k_A(\nabla\cd_S)|_{(t,y)}\le C'_k\nu(t)
\end{equation}
for every $t\ge0$ and non-negative integer $k$, where $C'_k$ is
some constant. Here we are using the simple fact that if
$E,E'$ are Banach spaces, $W\subset E$
an open neighbourhood of $0$, and $h\col W\to E'$ a differentiable map with 
$h(0)=0$ then
$\|h(x)\|\le(\|Dh(0)\|+1)\|x\|$ in some neighbourhood
of $0$. For instance, to deduce
the second inequality in \Ref{nu-t-bound} we can apply this to the map
\[h\col L^2_{k+j+1}\to L^2_j,\quad
s=(a,\phi)\mapsto\nabla^k_{B'+a}(\nabla\cd_{\uline\beta+s})\]
where $j\ge3$, say, and $B'$ is the connection part of $\uline\beta$.

Putting the inequalities above together we get
\[\sup_{y\in Y}|\nabla^k_A(\nabla\cd_S)|_{(t,y)}\le
C''_k\sqrt{\nu(0)}e^{-\lla^+t},\quad t\ge1\]
for some constants $C''_k$.
If $S$ is in temporal gauge we deduce, by taking $k=0$, that
$S_t$ converges uniformly to some continuous configuration $\al$.
One can now prove by induction on $k$ that $\al$ is of class $C^k$ and
that \Ref{lambda-decay} holds.
\end{proof}

\begin{thm}\label{exp-decay2}
For any $C>0$ there are constants $\eps,C_0,C_1,\dots$ such that the
following holds for every $T>1$. Let $S=(A,\Phi)$ be any smooth monopole
in temporal gauge over the band
$[-T-2,T+2]\times Y$, and suppose $[S_t]\in U$ for some
$t\in[-T,T]$. Set
\[\bar\nu=\|\nabla\cd_S\|_{L^2((-T-2,T+2)\times Y)},\quad
\nu(t)=\|\nabla\cd_S\|_{L^2(\BB_t)}.\]
If $\|\Phi\|_\infty\le C$ and $\bar\nu\le\eps$ then there is a smooth monopole
$\al$ over $Y$, gauge equivalent to $\beta$,
such that if $B$ is the connection part of $\ual$ then
for $|t|\le T-1$ and every non-negative integer $k$ one has
\begin{equation*}
\sup_{y\in Y}|\nabla^k_B(S-\ual)|_{(t,y)}
\le C_k(\nu(-T)+\nu(T))^{1/2}e^{-\lla(T-|t|)}.
\end{equation*}
\end{thm}

\begin{proof}
Given $C>0$, let $S=(A,\Phi)$ be any monopole over
$[-T-2,T+2]\times Y$ such that $\|\Phi\|_\infty\le C$
and $[S_t]\in U$ for some $t\in[-T,T]$. If $\bar\nu$ is sufficiently small
then we can define the function $f(t)$ for $|t|\le T$ as in the
proof of Theorem~\ref{exp-decay1}, and \Ref{nu-t-bound} will hold with
$f(t)$ replaced by $|f(t)|$, for
$|t|\le T$. Again, $f(t)=\cd(S_t)-L$ for some constant $L$.
Lemma~\ref{csd-est1} now gives
\[e^{-2\lla_-(T-t)}f(T)\le f(t)\le e^{-2\lla_+(T+t)}f(-T),\quad|t|\le T,\]
which implies
\begin{gather*}
|f(t)|\le(|f(-T)|+|f(T)|)e^{-2\lla(T-|t|)},\quad|t|\le T,\\
\nu(t)^2\le\const(\nu(-T)+\nu(T))e^{-2\lla(T-|t|)},\quad|t|\le T-1.
\end{gather*}
By Propositions~\ref{gauge3} and \ref{int-est} there is a critical point $\al$
gauge equivalent to $\beta$ such that
\[\|\nabla_B^k(S_0-\al)\|_{L^\infty(Y)}\le C'''_k\nu(0)\]
for some constants $C'''_k$.
It is now easy to complete the proof by induction on $k$.
\end{proof}

\subsection{Global convergence}
\label{subsec:global-conv}

The main result of this subsection is Proposition~\ref{half-cyl-conv}, which
relates local and global convergence of monopoles over a half-cylinder.
First some lemmas.

\begin{lemma}\label{du-ineq}
Let $Z$ be a compact Riemannian $n$--manifold (perhaps with boundary),
$m$ a non-negative integer, and $q\ge n/2$. Then there is
a real polynomial $P_{m,q}(x)$ of degree $m+1$ satisfying $P_{m,q}(0)=0$,
such that for any smooth $u\col Z\to\U1$ one has
\[\|du\|_{m,q}\le P_{m,q}(\|u\inv du\|_{m,q}).\]
\end{lemma}

\begin{proof}
Argue by induction on $m$ and use the Sobolev embedding
$L^r_k(Z)\subset L^{2r}_{k-1}(Z)$ for $k\ge1$, $r\ge n/2$.
\end{proof}

\begin{lemma}\label{df-ineq1}
Let $Z$ be a compact, connected
Riemannian $n$--manifold (perhaps with boundary), $z\in Z$,
$m$ a positive integer, and $q\ge 1$. Then there is a $C>0$
such that for any smooth $f\col Z\to\co$ one has
\begin{enumerate}
\item[\rm(i)]$\|f-f\av\|_{m,q}\le C\|df\|_{m-1,q}$,
\item[\rm(ii)]$\|f\|_{m,q}\le C(\|df\|_{m-1,q}+|f(z)|)$,
\end{enumerate}
where $f\av=\text{Vol}(Z)\inv\int_Zf$ is the average of $f$.
\end{lemma}

\begin{proof}
Exercise.
\end{proof}

\begin{lemma}\label{df-ineq2}
Let $Z$ be a compact Riemannian $n$--manifold (perhaps with boundary),
$m$ a positive integer, and $q$ a real number such that $mq>n$.
Let $\Phi$ be a smooth section of some Hermitian vector bundle
$E\to Z$, $\Phi\not\equiv0$. Then there exists a $C>0$ with the following
significance. Let $\phi_1$ be a smooth section of $E$ satisfying
$\|\phi_1\|_q\le C\inv$ and $w\col Z\to\co$ a smooth map. Define another section
$\phi_2$ by
\[w(\Phi+\phi_1)=\Phi+\phi_2.\]
Then
\[\|w-1\|_{m,q}\le C(\|dw\|_{m-1,q}+\|\phi_2-\phi_1\|_q).\]
\end{lemma}

\begin{proof}
The equation
\[(w-1)\Phi=\phi_2-\phi_1-(w-1)\phi_1\]
gives
\begin{align*}
\|w-1\|_{m,q}&\le\const(\|dw\|_{m-1,q}+\|(w-1)\Phi\|_q)\\
&\le\const(\|dw\|_{m-1,q}+\|\phi_2-\phi_1\|_q+\|w-1\|_{m,q}\|\phi_1\|_q).
\end{align*}
Here the first inequality is analogous to Lemma~\ref{df-ineq1}~(ii).
If $\|\phi_1\|_q$ is sufficiently small then we can rearrange to get the
desired estimate on $\|w-1\|_{m,q}$.
\end{proof}

Now let $\al$ be a non-degenerate critical point of $\cd$.
Note that if $S=(A,\Phi)$ is any finite energy
monopole in temporal gauge over
$\rpy$ such that $\|\Phi\|_\infty<\infty$ and
\[\liminf_{t\to\infty}\int_{[t,t+1]\times Y}|S-\ual|^r=0\]
for some $r>1$ then by the results of Section~\ref{sec:local-comp1} we have
$[S_t]\to[\al]$ in $\cb_Y$, hence $S-\ual$ decays exponentially 
by Theorem~\ref{exp-decay1}. In this situation we will simply say that
$S$ is {\em asymptotic to $\al$.}

(Here we used the fact that for any $p>2$ and $1<r\le2p$, say,
the $L^r$ metric on the
$L^p_1$ configuration space $\cb([0,1]\times Y)$ is well-defined.)

\begin{lemma}\label{temp-gauge}
If $S=(A,\Phi)$ is any smooth monopole over $\rpy$ such that
$\|\Phi\|_\infty<\infty$ and $S-\ual\in L^p_1$ for some $p>2$ then
there exists a null-homotopic smooth $u\col \rpy\to\U1$ such that
$u(S)$ is in temporal gauge and asymptotic to $\al$.
\end{lemma}

\begin{proof}
By Theorem~\ref{exp-decay1} there exists a smooth
$u\col \rpy\to\U1$ such that $u(S)$ is in temporal gauge and asymptotic
to $\al$.
Lemma~\ref{du-ineq}, Lemma~\ref{df-ineq1}~(i),
and the assumption $S-\ual\in L^p_1$ then gives
\[\|u-u\av\|_{L^\infty([t,t+1]\times Y)}\to0~\text{as $t\to\infty$},\]
hence $u$ is null-homotopic.
\end{proof}

It follows that all elements of the moduli spaces defined in
Subsection~\ref{subsec:modsp} have smooth representatives that are in temporal
gauge over the ends.

\begin{prop}\label{half-cyl-conv}
Let $\del>0$ and suppose $\cd\col \ti\cb_Y\to\R$ has no critical value
in the half-open
interval $(\cd(\al),\cd(\al)+\del]$ (this implies Condition~(O1)).
For $n=1,2,\dots$ let $S_n=(A_n,\Phi_n)$
be a smooth monopole over $\orp\times Y$ such that
\[S_n-\ual\in L^p_1,\quad \sup_n\|\Phi_n\|_\infty<\infty,\quad
\cd(S_n(0))\le\cd(\al)+\del,\]
for some $p>2$.
Let $v_n\col \orp\times Y\to\U1$ be a smooth map such that the sequence
$v_n(S_n)$ converges in $C^\infty$ over compact subsets of $\orp\times Y$ to a 
configuration $S$ in temporal gauge. Then the following hold:
\begin{enumerate}
\item[\rm(i)]$S$ is asymptotic to a critical point $\al'$ gauge equivalent to
$\al$,
\item[\rm(ii)]If $\al=\al'$ then $v_n$ is null-homotopic for all sufficiently
large $n$, and there exist smooth
$u_n\col \orp\times Y\to\U1$ with the following significance: For every $t\ge0$
one has $u_n=1$ on $[0,t]\times Y$ for all sufficiently large $n$.
Moreover, for any $\si<\lla^+$, $q\ge1$ and non-negative integer $m$ one has
\[\|u_nv_n(S_n)-S\|_{\llw q\si m}\to0~\text{as $n\to\infty$}.\]
\end{enumerate}
\end{prop}

Here $\lla^+$ is as in Subsection~\ref{subsec:diff-ineq}.

\begin{proof}
It clearly suffices to prove the proposition
when $q\ge2$ and $mq>4$, which we assume from now on.

By Lemma~\ref{cdlim} we have
\begin{equation}\label{intcd}
\int_{\rpy}|\nabla\cd_{S_n}|^2=\cd(S_n(0))-\cd(\al)\le\del
\end{equation}
for each $n$, hence $\int_{\rpy}|\nabla\cd_S|^2\le\del$.
Part~(i) of the proposition
is now a consequence of Theorem~\ref{exp-decay1} and the following

\begin{claim}\label{stal}
$[S(t)]$ converges in $\cb_Y$ to $[\al]$ as $t\to\infty$.
\end{claim}

{\bf Proof of claim}\qua For $r>0$ let $B_r\subset\cb_Y$
denote the open $r$--ball around $[\al]$ in the $L^2$ metric, and let
$\oline B_r$ be the corresponding closed ball. Choose $r>0$ such that
$\oline B_{2r}$ contains no monopole other than $[\al]$. Assuming the
claim does not hold then by Lemma~\ref{lemma:loc-comp}
one can find a sequence $t'_j$ such that
$t'_j\to\infty$ as $j\to\infty$ and $[S(t'_j)]\not\in\oline B_{2r}$ for
each $j$. Because of the convergence of $v_n(S_n)$
it follows by a continuity argument that there are sequences $n_j,t_j$
with $t_j,n_j\to\infty$ as $j\to\infty$, such that
\[[S_{n_j}(t_j)]\in\oline B_{2r}\setminus B_r\]
for all $j$.
For $s\in\R$ let $\ct_s\col \ry\to\ry$ be translation by $s$:
\[\ct_s(t,y)=(t+s,y).\]
Again by Lemma~\ref{lemma:loc-comp} there are smooth
$\om_j\col \orp\times Y\to\U1$ such that a subsequence of 
$(\ct_{t_j})^*(\om_j(S_{n_j}))$ converges in $C^\infty$ over compact subsets of
$\ry$ to some finite energy monopole $S'$ whose spinor field is 
pointwise bounded.
Moreover, it is clear that $\cd\circ\om_j(0)-\cd\in\R$ must be bounded as
$j\to\infty$, so by passing to a subsequence and replacing $\om_j$ by
$\om_j\om_{j_0}\inv$ for some fixed $j_0$ we may arrange that
$\cd\circ\om_j(0)=\cd$ for all $n$. 
Then $\ell=\lim_{t\to-\infty}\cd(S'(t))$ must
be a critical value of $\cd$. Since
\[[S'(0)]\in\oline B_{2r}\setminus B_r,\]
$S'(0)$ is not a critical point, whence $\prtl_t|_0\cd(S'(t))<0$. Therefore,
\[\cd(\al)+\del\ge\ell>\cd(S'(0))>\cd(\al),\]
contradicting our assumptions. This proves the claim.
\endproof

We will now prove Part~(ii).
For $\tau\ge0$ let
\[B^-_\tau=[0,\tau]\times Y,\quad B^+_\tau=[\tau,\infty)\times Y,\quad
\cO_\tau=[\tau,\tau+1]\times Y.\]
By Lemma~\ref{temp-gauge}
there is, for every $n$, a null-homotopic, smooth $\ti v_n\col \orp\times Y\to\U1$
such that $S''_n=\ti v_n(S_n)$ is in temporal gauge and asymptotic to $\al$.

Note that 
\[\lim_{t\to\infty}\limsup_{n\to\infty}\cd(S_n(t))=\cd(\al).\]
For otherwise we could find an $\eps>0$ and for every natural number $j$ a pair
$t_j,n_j\ge j$ such that
\[\cd(S_{n_j}(t_j))\ge\cd(\al)+\eps,\]
and we could then argue as in the proof of Claim~\ref{stal} to produce a 
critical value of $\cd$ in the interval $(\al,\al+\del]$.
Since $|\nabla\cd_{S_n}|=|\nabla\cd_{S''_n}|$ it follows from \Ref{intcd} and
Theorem~\ref{exp-decay1} that there exists
a $t_1\ge0$ such that if $\tau\ge t_1$ then
\begin{equation*}
\limsup_{n\to\infty}\|S''_n-\ual\|_{\llw q\si m(B^+_\tau)}
\le\const e^{(\si-\lla^+)\tau}
\end{equation*}
where the constant is independent of $\tau$.
Then we also have
\[\limsup_{n\to\infty}\|S''_n-S\|_{\llw q\si m(B^+_\tau)}\le
\const e^{(\si-\lla^+)\tau}.\]
Set $S'_n=v_n(S_n)$ and $w_n=\ti v_nv_n\inv$. Then we get
\[\limsup_{n\to\infty}\|S''_n-S'_n\|_{L^q_m(\cO_\tau)}\le
\const e^{-\lla^+\tau},\]
which gives
\begin{equation*}
\limsup_{n\to\infty}\|dw_n\|_{L^q_m(\cO_\tau)}\le \const e^{-\lla^+\tau}
\end{equation*}
by Lemma~\ref{du-ineq}. In particular, $w_n$ is
null-homotopic for all sufficiently large $n$.

Fix $y_0\in Y$ and set $x_\tau=(\tau,y_0)$.
Choose a sequence $\tau_n$ such that $\tau_n\to\infty$ as $n\to\infty$ and
\begin{equation}\label{tn-est}
\begin{aligned}
{}&\|S'_n-S\|_{\llw q\si m(B^-_{\tau_n+1})}\to0,\\
&\|S''_n-S\|_{\llw q\si m(B^+_{\tau_n})}\to0
\end{aligned}
\end{equation}
as $n\to\infty$. If $\al$ is reducible then by multiplying each
$\ti v_n$ by a constant and redefining
$w_n,S''_n$ accordingly we may arrange that $w_n(x_{\tau_n})=1$ for all $n$.
(If $\al$ is irreducible we keep $\ti v_n$ as before.)
Then \Ref{tn-est} still holds. Applying
Lemma~\ref{du-ineq} together with Lemma~\ref{df-ineq1}~(ii) (if $\al$
is reducible) or Lemma~\ref{df-ineq2} (if $\al$ is irreducible) we see that
\[e^{\si\tau_n}\|w_n-1\|_{L^q_{m+1}(\cO_{\tau_n})}\to0\]
as $n\to\infty$.

Let $\beta\col \R\to\R$ be a
smooth function such that $\beta(t)=0$ for $t\le\frac13$
and $\beta(t)=1$ for $t\ge\frac23$. Set $\beta_\tau(t)=\beta(t-\tau)$.
Given any function $w\col \cO_\tau\to\co\setminus(-\infty,0]$ define
\[\mathcal U_{w,\tau}=\exp(\beta_\tau\log w)\]
where $\log(\exp(z))=z$ for complex numbers $z$ with $|\Im(z)|<\pi$.
Let $m'$ be any integer such that $m'q>4$.
If $\|w-1\|_{m',q}$ is sufficiently small then 
\begin{equation}\label{fder}
\begin{aligned}
\|\mathcal U_{w,\tau}-1\|_{m',q}&\le\const\|w-1\|_{m',q},\\
\|w\inv dw\|_{m'-1,q}&\le\const\|w-1\|_{m',q}.
\end{aligned}
\end{equation}
To see this recall that for functions on $\rf$, multiplication defines a 
continuous map $L^q_{m'}\times L^q_k\to L^q_k$ for $0\le k\le m'$. Therefore,
if $V$ is the set of all functions in $L^q_{m'}(\cO_\tau,\co)$ that map into
some fixed small open ball about $1\in\co$ then $w\mapsto\mathcal U_{w,\tau}$
defines a $C^\infty$ map $V\to L^q_{m'}$. This yields the first inequality in
\Ref{fder}, and the proof of the second inequality is similar.

Combining \Ref{tn-est} and \Ref{fder} we conclude that Part~(ii) of the
proposition holds with
\[u_n=\begin{cases}
1& \text{in $B^-_{\tau_n}$},\\
\mathcal U_{w_n,\tau_n} & \text{in $\cO_{\tau_n}$},\\
w_n & \text{in $B^+_{\tau_n+1}$}.
\end{cases}\]
This completes the proof of Proposition~\ref{half-cyl-conv}.
\end{proof}

\section{Global compactness}\label{sec:global-comp}

In this section we will prove Theorems~\ref{thm:sing-comp} and
\ref{thm:neck-comp}. Given the results of Sections~\ref{sec:local-comp1}
and \ref{sec:local-comp2}, what remains to be understood is convergence over
ends and necks.
We will use the following terminology:
\[\text{c-convergence = $C^\infty$ convergence over compact subsets.}\]

\subsection{Chain-convergence}\label{subsec:chain-conv}

We first define the notion of chain-convergence. For simplicity we only
consider two model cases: first the case of one end and no necks, then the
case of one neck and no ends. It should be clear how to extend the notion
to the case of multiple ends and/or necks.

\begin{defn}\label{defn:chain-conv}
{\rm
Let $X$ be a $\spc$ Riemannian $4$--manifold with one tubular end $\rpy$,
where $Y$ is connected. Let $\al_1,\al_2,\dots$ and $\beta_0,\dots,\beta_k$
be elements of $\ti\fl_Y$,
where $k\ge0$ and $\cd(\beta_{j-1})>\cd(\beta_j)$ for $j=1,\dots,k$.
Let $\om\in M(X;\beta_0)$ and
$\vec v=(v_1,\dots,v_k)$, where $v_j\in\check M(\beta_{j-1},\beta_j)$. 
We say a sequence $[S_n]\in M(X;\al_n)$ {\em chain-converges}
to $(\om,\vec v)$ if there exist, for each $n$,
\begin{itemize}
\item a smooth map $u_n\col X\to\U1$,
\item for $j=1,\dots,k$ a smooth map $u\nj\col \ry\to\U1$,
\item a sequence $0=t_{n,0}<t_{n,1}<\dots<t_{n,k}$,
\end{itemize}
such that
\begin{enumerate}
\item[(i)]$u_n(S_n)$ c-converges over $X$ to a representative of $\om$
(in the sense of Subsection~\ref{local-slices}),
\item[(ii)]$t_{n,j}-t_{n,j-1}\to\infty$ as $n\to\infty$,
\item[(iii)]$u\nj(\ct^*_{t\nj}S_n)$ c-converges
over $\ry$ to a representative of $v_j$,
\item[(iv)]$\limsup_{n\to\infty}\left[\cd(S_n(t_{n,j-1}+\tau))
-\cd(S_n(t\nj-\tau))\right]\to0$ as $\tau\to\infty$,
\item[(v)]$\limsup_{n\to\infty}\left[\cd(S_n(t_{n,k}+\tau))-\cd(\al_n)\right]
\to0$ as $\tau\to\infty$,
\end{enumerate}
where (ii), (iii) and (iv) should hold for $j=1,\dots,k$.}
\end{defn}

Conditions (iv) and (v) mean, in familiar language, that ``no energy is lost
in the limit''. As before, $\ct_s$ denotes translation by $s$, ie
$\ct_s(t,y)=(t+s,y)$.

We now turn to the case of one neck and no ends.

\begin{defn} 
{\rm In the situation of Subsection~\ref{intro:comp-neck}, suppose $r=1$ and
$r'=0$. Let $\beta_0,\dots,\beta_k\in\ti\fl_Y$,
where $k\ge0$ and $\cd(\beta_{j-1})>\cd(\beta_j)$, $j=1,\dots,k$.
Let $\om\in M(X;\beta_0,\beta_k)$ and $\vec v=(v_1,\dots,v_k)$,
where $v_j\in\check M(\beta_{j-1},\beta_j)$. Let $T(n)\to\infty$ as
$n\to\infty$. We say a sequence $[S_n]\in M(\xtn)$ {\em chain-converges}
to $(\om,\vec v)$ if there exist, for every $n$, 
\begin{itemize}
\item a smooth map $u_n\col \xtn\to\U1$,
\item for $j=1,\dots,k$ a smooth map $u\nj\col \ry\to\U1$,
\item a sequence $-T(n)=t_{n,0}<t_{n,1}<\dots<t_{n,k+1}=T(n)$,
\end{itemize}
such that (i)--(iv) of Definition~\ref{defn:chain-conv} hold for the values
of $j$ for which they are defined (in other words,
(ii) and (iv) should hold for $1\le j\le k+1$
and (iii) for $1\le j\le k$).}
\end{defn}

In the notation of Subsection~\ref{intro-csd}, if $J\subset\R$ is an
interval with non-empty interior then a smooth configuration
$S$ over $J\times Y$ is called {\em regular} (with respect to $\cd$)
if either $\prtl_t\cd(S_t)<0$ for every $t\in J$,
or $S$ is gauge equivalent to the translationary invariant configuration
$\ual$ determined by some critical point $\al$ of $\cd$.
Proposition~\ref{Dcd} guarantees the regularity of certain
$(\prt,\fq)$--monopoles when $\prt$ is sufficiently small. In particular,
genuine monopoles are always regular.

Consider now the situation of Subsection~\ref{intro:comp-neck}
(without assuming (B1) or (B2)), and let the $2$--form $\mu$ on $X$ be fixed.

\begin{defn}\label{defn:suff-small}
{\rm A set of perturbation parameters $\vec\prt,\vec\prt'$ is
{\em admissible} for a vector $\vec\al'$ of critical points if for some
$t_0\ge1$ the following holds. Let $\cm$ be the disjoint union of all
moduli spaces $M(\xt;\vec\al';\vec\prt;\vec\prt')$ with $\min_jT_j\ge t_0$.
Then we require, for all $j,k$, that
\begin{enumerate}
\item[(i)]If $\ti S$ is any configuration over
$[-1,1]\times Y_j$ which is a $C^\infty$ limit of configurations
of the form $S|_{[t-1,t+1]\times Y_j}$ with $S\in\cm$
and $|t|\le T_j-1$, then $\ti S$ is regular.
\item[(ii)]If $\ti S$ is any configuration over
$[-1,1]\times Y'_{k}$ which is a $C^\infty$ limit of configurations
of the form $S|_{[t-1,t+1]\times Y'_{k}}$ with $S\in\cm$
and $t\ge1$, then $\ti S$ is regular.
\end{enumerate}}
\end{defn}
In particular, the zero perturbation parameters are always admissible.

The next two propositions describe some properties of chain-convergence.

\begin{prop}
In the notation of Theorem~\ref{thm:sing-comp}, suppose $\om_n$
chain-converges to $(\om,\vec v_1,\dots,\vec v_r)$,
where $\vec\al_n=\vec\beta$ for all $n$,
each $\vec v_j$ is empty, and $\vec\prt$ is admissible for $\vec\beta$.
Then $\om_n\to\om$ in $M(X;\vec\beta)$ with its
usual topology.
\end{prop}

\begin{proof}
This follows from Proposition~\ref{half-cyl-conv}.
\end{proof}

In other words, if a sequence $\om_n$
in a moduli space $M$ chain-converges to an element $\om\in M$,
then $\om_n\to\om$ in $M$ provided the perturbations are admissible.

\begin{prop}\label{prop:chain-conv}
In the notation of Theorem~\ref{intro:comp-neck}, suppose that the sequence
$\om_n\in M(\xtn;\vec\al'_n)$ chain-converges to
$\bv=(\om,\vec v_1,\dots,\vec v_r,\vec v'_1,\dots,\vec v'_{r'})$, where
$\min_jT_j(n)\to\infty$. Suppose also that the perturbation parameters
$\vec\prt,\vec\prt'$ are admissible for each $\vec\al'_n$. Then the
following hold:
\begin{enumerate}
\item[\rm(i)]For sufficiently large $n$ there is a smooth map $u_n\col \xtn\to\U1$
such that $v\nj=u_n|_{\{0\}\times Y'_j}$ satisfies $v\nj(\al'\nj)=\ga'_j$,
$j=1,\dots,r'$.
\item[\rm(ii)]The chain limit is unique up to gauge equivalence,
ie if $\mathbb V$, $\mathbb V'$ are two chain limits of $\om_n$ then there
exists a smooth $u\col X\gl\to\U1$ which is translationary invariant over the ends
of $X\gl$, and such that $u(\mathbb V)=\mathbb V'$.
\end{enumerate}
\end{prop}

In (i), recall that moduli spaces are labelled by
critical points modulo null-homotopic gauge transformations.
Note that we can arrange that the maps $u_n$  are translationary
invariant over the ends. This allows us to identify the moduli spaces
$M(X;\vec\al'_n)$ and $M(X;\vec\ga')$, so that we obtain a sequence
$u_n(\om_n)$, $n\gg0$ in a {\em fixed} moduli space.

In (ii) we define $u(\mathbb V)$ as follows.
Let $w_j\col \ry_j\to\U1$ and $w'_j\col \ry'_j\to\U1$ be the translationary invariant
maps which agree with $u$ on $\{0\}\times Y_j$ and $\rpy'_j$, respectively.
Let $w\col X\to\U1$ be the map which is translationary invariant over each end
and agrees with $u$ on $X\endt 1$. Then $u(\mathbb V)$ is the result of
applying the appropriate maps $w,w_j,w'_j$ to the various components of
$u(\mathbb V)$.

\begin{proof} (i)\qua For simplicity we only discuss the case of one
end and no necks, ie the situation of Definition~\ref{defn:chain-conv}.
The proof in the general case is similar. 

Using Condition~(v) of Definition~\ref{defn:chain-conv} and a simple
compactness argument it is easy to see that $\al_n$ is gauge equivalent
to $\beta_k$ for all sufficiently large $n$. Moreover,
Conditions~(iv) and (v) of Definition~\ref{defn:chain-conv}
ensure that there exist
$\tau,n'>0$ such that if $n>n'$ then $\om_n$ restricts to a genuine
monopole on $(t_{n,k}+\tau,\infty)\times Y$ and on
$(t_{n,j-1}+\tau,t\nj-\tau)\times Y$ for $j=1,\dots,k$. It then follows from
Proposition~\ref{half-cyl-conv} that $v_n=u_{n,k}|_{\{0\}\times Y}$
satisfies $v_n(\al_n)=\beta_k$ for $n\gg0$. 
(Recall again that
$\al_n,\beta_k\in\ti\fl_Y$ are critical points modulo null-homotopic
gauge transformations, so $v_n(\al_n)$ depends only on the homotopy class
of $v_n$.)
Similarly, it follows from Theorems~\ref{exp-decay1},
\ref{exp-decay2}
that $u_{n,j-1}|_{\{0\}\times Y}$ is homotopic to
$u\nj|_{\{0\}\times Y}$ for $j=1,\dots,r$
and $n\gg0$, where $u_{n,0}=u_n$. Therefore, $v_n$ extends over $X\endt0$.

(ii)\qua This is a simple exercise.
\end{proof}

\subsection{Proof of Theorem~\ref{thm:sing-comp}}
\label{subsec:single}

By Propositions~\ref{phi-bound},\ref{prop:local-energy-bound2}, and \ref{Dcd},
if each $\prt_j$ has sufficiently small $C^1$ norm then $\vec\prt$ will
be admissible for all $\vec\al$. Choose $\vec\prt$ so that this is the case.
Set
\[C_0=-\inf_n\sum_j\lla_j\cd(\al\nj)<\infty.\]

Let $S_n$ be a smooth representative for $\om_n$. 
The energy assumption on the asymptotic
limits of $S_n$ is unaffected if we replace $S_n$ by $u_n(S_n)$ for 
some smooth $u_n\col X\to\U1$ which is translationary invariant on
$(t_n,\infty)\times Y$ for some $t_n>0$. After passing to a subsequence
we can therefore, by Proposition~\ref{prop:loc-comp}, assume that
$S_n$ c-converges over $X$ to some monopole $S'$ which is in temporal
gauge over the ends. Because $\vec\prt$ is admissible we have that
\[\prtl_t\cd(S_n|_{\{t\}\times Y_j})\le0\]
for all $j,n$ and $t\ge0$. From the energy bound \Ref{energy-bound}
we then see that $S'$ must have finite energy.
Let $\ga_j$ denote the asymptotic limit of $S'$ over the end $\rpy_j$
as guaranteed by Proposition~\ref{crit-lim}. Then
\[\limsup_n\cd(\al\nj)\le\cd(\ga_j)\]
for each $j$. Hence there is a constant $C_2<\infty$
such that for $h=1,\dots,r$ and all $n$ one has
\[C_2+\lla_h\cd(\al_{n,h})\ge\sum_j\lla_j\cd(\al\nj)\ge-C_0.\]
Consequently,
\[\sup\nj|\cd(\al\nj)|<\infty.\]
For the remainder of this proof we fix $j$ and focus on one end $\ry_j$.
For simplicity we drop $j$ from notation and write $Y,\al_n$ instead of
$Y_j,\al\nj$ etc.

After passing to a subsequence we may arrange that $\cd(\al_n)$ has the
same value $L$ for all $n$ (here we use Condition~(O1)).
If $\cd(\ga)=L$ then we set $k=0$ and the proof is complete.
Now suppose $\cd(\ga)>L$. Then there is an
$n'$ such that $\prtl_t\cd(S_n(t))<0$ for all $n\ge n'$, $t\ge0$. Set
\[\del=\frac12\min\{|x-y|\st\text{$x,y$ are distinct critical values of
$\cd\col \ti\cb_Y\to\R$}\}.\]
The minimum exists by (O1). For sufficiently large $n$ we define
$t_{n,1}\gg0$ implicitly by
\[\cd(S_n(t_{n,1}))=\cd(\ga)-\del.\]
It is clear that $t_{n,1}\to\infty$ as $n\to\infty$. 
Moreover, Definition~\ref{defn:chain-conv}~(iv) must hold for $j=1$.
For otherwise we can find $\eps>0$ and sequences $\tau_\ell$, $n_\ell$ with
$\tau_\ell,n_\ell\to\infty$ as $\ell\to\infty$, such that
\begin{equation}\label{eqn:tau-ell}
\cd(S_{n_\ell}(\tau_\ell))-\cd(S_{n_\ell}(t_{n_\ell,1}-\tau_\ell))>\eps
\end{equation}
for every $\ell$. As in the proof of Claim~\ref{stal} there are smooth
$\ti u_\ell\col \ry\to\U1$ satisfying $\cd\circ\ti u_\ell(0)=\cd$ such that
a subsequence of $\ti u_\ell(\ct^*_{t_{n_\ell,1}}S_{n_\ell})$ c-converges
over $\ry$ to a finite energy monopole $\ti S$ in temporal gauge.
The asymptotic limit $\ti \ga$ of $\ti S$ at $-\infty$
must satisfy
\[\eps\le\cd(\ga)-\cd(\ti\ga)<\del,\]
where the first inequality follows from \Ref{eqn:tau-ell}.
This contradicts the choice of $\del$. Therefore, 
Definition~\ref{defn:chain-conv}~(iv) holds for $j=1$ as claimed.

After passing to a subsequence we can find
$u_{n,1}\col \ry\to\U1$ such that $u_{n,1}(\ct^*_{t_{n,1}}S_n)$ c-converges
over $\ry$ to some finite energy monopole $S'_1$ in temporal gauge.
Let $\beta_1^\pm$ denote the limit of $S'_1$ at $\pm\infty$.
A simple compactness argument shows that $\ga$ and $\beta_1^-$
are gauge equivalent, so we can arrange that $\ga=\beta_1^-$ by
modifying the $u_{n,1}$ by a fixed gauge transformation $\ry\to\U1$.
As in the proof of Proposition~\ref{prop:chain-conv}~(i) we see that
$u_{n,1}$ must be null-homotopic for all sufficiently large $n$.
Hence $\cd(\beta^+_1)\ge L$. If $\cd(\beta^+_1)=L$ then we set $k=1$ and the
proof is finished. If on the other hand $\cd(\beta^+_1)>L$ 
then we continue the above process.
The process ends when, after passing successively to subsequences and choosing
$u\nj,t\nj,\beta_j^\pm$ for $j=1,\dots,k$
(where $\beta^+_{j-1}=\beta^-_j$, and $u\nj$ is null-homotopic for $n\gg0$)
we have $\cd(\beta^+_k)=L$. This must occur
after finitely many steps; in fact $k\le(2\del)\inv(\cd(\ga)-L)$.
\endproof

\subsection{Proof of Theorem~\ref{thm:neck-comp}}

For simplicity we first consider the case
when there is exactly one neck (ie $r=1$), and we write $Y=Y_1$ etc.
We will make repeated use of the local compactness results proved earlier.

Let $S_n$ be a smooth representative of $\om_n$.
After passing to a subsequence we can find smooth maps
\[u_n\col \xtn\setminus(\{0\}\times Y)\to\U1\]
such that $\ti S_n=u_n(S_n)$ c-converges over $X$ to some
finite energy monopole $S'$ which is in temporal gauge over the ends.
Introduce the temporary notation $S_n(t)=S_n|_{\{t\}\times Y}$, and similarly
for $\ti S_n$ and $u_n$. For $0\le\tau<T(n)$ set
\[\Ttn=\cd(S_n(-T(n)+\tau))-\cd(S_n(T(n)-\tau)).\]
Let $u^\pm_n=u_n(\pm T(n))$ and
\[I^\pm_n=2\pi\int_Y\ti\eta_j\wedge[u^\pm_n],\]
cf Equation~\Ref{cduS}. Since $\Theta_{0,n}$
is bounded as $n\to\infty$,
it follows that $I^+_n-I^-_n$ is bounded as $n\to\infty$.
By Condition~(O1) there is a $q>0$ such that $qI^\pm_n$ is integral for all
$n$. Hence we can arrange, by passing to a subsequence,
that $I^+_n-I^-_n$ is constant. In particular,
\[I^+_n-I^+_1=I^-_n-I^-_1.\]
Choose a smooth map $w\col X\to\U1$ which is translationary invariant over the
ends, and homotopic to $u_1\inv$ over $X\endt0$.
After replacing $u_n$ by $wu_n$ for every $n$ we then
obtain $I^+_n=I^-_n$. Set $I_n=I^\pm_n$.
We now have
\[\Ttn=\cd(\ti S_n(-T(n)))-\cd(\ti S_n(T(n))).\]
Let $\beta_0,\beta'$ denote the asymptotic limits of $S'$ over the ends
$\iota^+(\rpy)$ and $\iota^-(\rpy)$, respectively. 
Set
\[L=\lim_{\tau\to\infty}\lim_{n\to\infty}\Ttn=\cd(\beta_0)-\cd(\beta').\]
Since $\Ttn\ge0$ for $\tau\ge0$ we have $L\ge0$.

Suppose $L=0$. Then a simple compactness argument shows that there is a
smooth $v\col Y\to\U1$ such that $v(\beta_0)=\beta'$. Moreover, there is an $n_0$
such that $v\,u^-_n\sim u^+_n$ for $n\ge n_0$, where $\sim$ means `homotopic'.
Therefore, we can find a smooth $z\col X\to\U1$ which is translationary 
invariant over the ends and homotopic to $u_{n_0}\inv$ over $X\endt0$,
such that after replacing $u_n$ by $zu_n$ for
every $n$ we have that $\beta_0=\beta'$ and $u^+_n\sim u^-_n$.
In that case we can in fact assume that $u_n$ is a smooth map 
$\xtn\to\U1$. The remainder of the proof when $L=0$ (dealing with
convergence over the ends) is now a repetition of the proof of
Theorem~\ref{thm:sing-comp}.

We now turn to the case $L>0$. For large $n$ we must then have
$\prtl_tS_n(t)<0$ for $|t|\le T(n)$. Let $\del$
be as in the proof of Theorem~\ref{thm:sing-comp}.
We define $t_{n,1}\in(-T(n),T(n))$ implicitly for large $n$ by
\[\cd(\beta_0)=\cd(S_n(t_{n,1}))+I_n+\del.\]
Then $|\tno\pm T(n)|\to\infty$ as $n\to\infty$. As in the proof of 
Theorem~\ref{thm:sing-comp} one sees that
\[\limsup_{n\to\infty}\left[\cd(S_n(-T(n)+\tau))
-\cd(S_n(\tno-\tau))\right]\to0\]
as $\tau\to\infty$, and after passing to a subsequence we can find smooth
$u_{n,1}\col \ry\to\U1$ such that $u_{n,1}(\ct^*_{\tno}S_n)$ c-converges
over $\ry$ to a finite energy monopole $S'_1$ in temporal gauge whose
asymptotic limit at $-\infty$ is $\beta_0$. Let $\beta_1$ denote
the asymptotic limit of $S'_1$ at $\infty$. We now repeat the above process.
The process ends when, after passing successively to subsequences and choosing
$u\nj,t\nj,\beta_j$ for $j=1,\dots,k$ one has that
\[\limsup_{n\to\infty}\left[\cd(S_n(t_{n,k}+\tau))
-\cd(S_n(T(n)-\tau))\right]\to0\]
as $\tau\to\infty$. As in the case $L=0$
one sees that $\beta_k,\beta'$ must be gauge equivalent, and after
modifying $u_n,u\nj$ one can arrange that $\beta_k=\beta'$. This establishes
chain-convergence over the neck. As in the case $L=0$ we can in fact
assume that $u_n$ is a smooth map $\xtn\to\U1$, and the rest of the proof
when $L>0$ is again a repetition of the proof of
Theorem~\ref{thm:sing-comp}.

In the case of multiple necks one applies the above argument successively
to each neck. In this case, too, after passing to a subsequence one ends up with
smooth maps $u_n\col \xtn\to\U1$ such that $u_n(S_n)$ c-converges over $X$.
One can then deal with convergence over the ends as before.
\endproof

\section{Transversality}\label{section:transv}

We will address two kinds of transversality problems: non-degeneracy of
critical points of the Chern--Simons--Dirac functional, and regularity of
moduli spaces over $4$--manifolds. 

In this section we do not assume Condition~(O1).

Recall that a subset of a topological space $Z$ is called {\em residual} if
it contains a countable intersection of dense open subsets of $Z$.

\subsection{Non-degeneracy of critical points}

\begin{lemma}\label{irr-crit-nondeg}
Let $Y$ be a closed, connected, Riemannian $\spc$ $3$--manifold and
$\eta$ any closed (smooth) $2$--form on $Y$. Let $G^*$ be the set of all
$\nu\in\Om^1(Y)$ such that all irreducible critical points of
$\cd_{\eta+d\nu}$ are non-degenerate. Then
$G^*\subset\Om^1(Y)$ is residual, hence dense
(with respect to the $C^\infty$ topology).
\end{lemma}

\begin{proof}
The proof is a slight modification of the argument in \cite{Fr1}.
For $2\le k\le\infty$ and $\del>0$
let $\wkd$ be the space of all $1$--forms $\nu$
on $Y$ of class $C^k$ which satisfy $\|d\nu\|_{C^1}<\del$.
Let $\wkd$ have the $C^k$ topology. 
For $1\le k<\infty$ we define a $\cg$--equivariant smooth map
\begin{align*}
\uk\col \cc^*\times L^2_1(Y;i\R)\times\wkd&\to L^2(Y;i\La^1\oplus\bs),\\
(B,\Psi,\xi,\nu)&\mapsto\ci_\Psi\xi + \nabla\cd_{\eta+d\nu}(B,\Psi),
\end{align*}
where $\cg$ acts trivially on forms, and by multiplication on spinors.
Now if $\uk(B,\Psi,\xi,\nu)=0$ then
\[\|\ci_\Psi\xi\|_2^2=-\int_Y\la\nabla\cd_{\eta+d\nu}(B,\Psi),\ci_\Psi\xi\ra=0\]
by \Ref{ipncd}, which implies $\xi=0$ since $\Psi\neq0$.
The derivative of $\uk$ at a point $x=(B,\Psi,0,\nu)$ is
\begin{equation}\label{ups-deriv}
D\uk(x)(b,\psi,f,v)=H_{(B,\Psi)}(b,\psi)+\ci_\Psi f+(i{*}dv,0).
\end{equation}

Let $(B,\Psi)$ be any irreducible critical point of $\cd_{\eta}$. We will show
that $P=D\uk(B,\Psi,0,0)$ is surjective. Note that altering $(B,\Psi)$ by an
$L^2_2$ gauge transformation $u$ has the effect of replacing $P$ by $uPu\inv$.
We may therefore assume that $(B,\Psi)$ is smooth. Since
$P_1=\ci_\Psi+H_{(B,\Psi)}$ has surjective symbol, the image of the induced
operator $L^2_1\to L^2$ is closed and has finite codimension. The
same must then hold for $\im(P)$. Suppose $(b,\psi)\in L^2$ is orthogonal to $\im(P)$,
ie $db=0$ and $P_1^*(b,\psi)=0$. The second equation implies that $b$ and
$\psi$ are smooth, by elliptic regularity. Writing out the equations
we find as in \cite{Fr1} that on the complement of $\Psi\inv(0)$ we have
$-b=idr$ for some smooth function $r\col Y\setminus\Psi\inv(0)\to\R$. We now
invoke a result of B\"ar \cite{Baer1} which says that, because $B$ is smooth
and $\Psi\not\equiv0$, the equation $\prtl_B\Psi=0$ implies that
the zero-set of $\Psi$ is contained in a countable union of smooth 
$1$--dimensional submanifolds of $Y$. 
In particular, any smooth loop in $Y$
can be deformed slightly so that it misses $\Psi\inv(0)$. Hence $b$ is exact.
 From B\"ar's theorem (or unique continuation for $\prtl_B$, which holds
when $B$ is of class $C^1$, see \cite{JKazdan1})
we also deduce that the complement of $\Psi\inv(0)$ is dense and connected.
Therefore, $f$ has a smooth extension to all of $Y$, and as in \cite{Fr1}
this gives $(b,\psi)=0$. Hence $P$ is surjective.

Consider now the vector bundle
\[E=(\cc^*\uset{\cg}\times L^2(Y;i\La^1\oplus\bs))\times L^2_1(Y;i\R)
\to\cb^*\times L^2_1(Y;i\R).\]
For $1\le k<\infty$
the map $\uk$ defines a smooth section $\skd$ of the bundle
\[E\times\wkd\to\cb^*\times L^2_1(Y;i\R)\times\wkd.\]
By the local slice theorem, a zero of $\uk$ is a regular point
of $\uk$ if and only if the
corresponding zero of $\skd$ is regular. Since surjectivity
is an open property for bounded operators between Banach spaces, 
a simple compactness argument shows that the zero-set of $\si_{2,\del}$
is regular when $\del>0$ is sufficiently small. Fix such a $\del$.
Observe that the question of whether the operator \Ref{ups-deriv} 
is surjective for a given $x$ is independent of $k$.
Therefore, the zero-set $\mkd$ of $\skd$ is regular for $2\le k<\infty$.
In the remainder of the proof assume $k\ge2$.

For any $\rho>0$ let $\cbr$ be the set of elements $[B,\Psi]\in\cb$
satisfying
\[\int_Y|\Psi|\ge\rho.\]
Define $\mkdr\subset\mkd$ similarly.
For any given $\nu$, the formula for $\uk$ defines a 
Fredholm section of $E$ which we denote by $\si_\nu$.
Let $\gkdr$ be the set of those $\nu\in\wkd$ such that
$\si_\nu$ has only regular zeros in $\cbr\times\{0\}$. For $k<\infty$ let
\[\pi\col \mkd\to\wkd\]
be the projection, and $\Si\subset\mkd$ the closed subset consisting of all
singular points of $\pi$.
A compactness argument shows that $\pi$ restricts to a closed map on $\mkdr$,
hence 
\[\gkdr=\wkd\setminus\pi(\mkdr\cap \Si)\]
is open in $\wkd$. On the other hand,
applying the Sard--Smale theorem as in \cite[Section~4.3]{DK} we see
that $\gkdr$ is residual (hence dense) in $\wkd$. Because $W_{\infty,\del}$
is dense in $\wkd$, we deduce that $G_{\infty,\del,\rho}$ is open and dense
in $W_{\infty,\del}$. But then
\[\bigcap_{n\in\n}G_{\infty,\del,\frac1n}\]
is residual in $W_{\infty,\del}$, and this is the set of all
$\nu\in W_{\infty,\del}$ such that $\si_\nu$ has only regular zeros. 

An irreducible critical point of $\cd_{\eta+d\nu}$ is
non-degenerate if and only if
the corresponding zero of $\si_\nu$ is regular. Thus we have proved that
among all smooth $1$--forms $\nu$ with $\|d\nu\|_{C^1}<\del$, those
$\nu$ for which all irreducible
critical points of $\cd_{\eta+d\nu}$ are non-degenerate
make up a residual subset in the $C^\infty$ topology.
The same must hold if $\eta$ is replaced with $\eta+d\nu$
for any $\nu\in\Om^1(Y)$, so we conclude that $G^*$ is locally residual in
$\Om^1(Y)$, ie any point in $G^*$ has a neighbourhood $V$ such that $G^*\cap V$
is residual in $V$. Hence $G^*$ is residual in $\Om^1(Y)$. (This last
implication holds if $\Om^1(Y)$ is replaced with any second countable, regular
space.) 
\end{proof}

\begin{prop}\label{frechet}
Let $Y$ be a closed, connected, Riemannian $\spc$ $3$--manifold and
$\eta$ any closed $2$--form on $Y$ such that either $b_1(Y)=0$
or $\ti\eta\neq0$. Let $G$ be the set of all $\nu\in\Om^1(Y)$ such that
all critical points of $\cd_{\eta+d\nu}$ are non-degenerate. Then
$G$ is open and dense in $\Om^1(Y)$ with respect to the $C^\infty$ topology.
\end{prop}

\begin{proof}
A compactness argument shows that $G$ is open.
If $b_1(Y)>0$ then $\cd_{\eta+d\nu}$ has no reducible critical points and
the result follows from Lemma~\ref{irr-crit-nondeg}.

Now suppose $b_1(Y)=0$. Then we may assume $\eta=0$. For $0\le k\le\infty$
let $W_k$ be the space of $1$--forms on $Y$ of class $C^k$, with the
$C^k$ topology.
If $\nu\in W_1$ then $\cd_{d\nu}$ has up to gauge equivalence a unique
reducible critical point, represented by $(B-i\nu,0)$ for any smooth
spin connection $B$ over $Y$ with $\check B$ flat.
This critical point is non-degenerate
precisely when 
\begin{equation}\label{kerdir}
\ker(\prtl_{B-i\nu})=0\quad\text{in $L^2_1$.}
\end{equation}
Let $G_k'$ be the set of all $\nu\in W_k$ such that \Ref{kerdir} holds.
This is clearly an open subset of $W_k$.
The last part of the proof of \cite[Proposition~3]{Fr1} shows that
$G'_k$ is residual (hence dense) in $W_k$ for $2\le k<\infty$.
Hence $G'_\infty$ is open and dense in $\Om^1(Y)=W_\infty$.
Now apply Lemma~\ref{irr-crit-nondeg}.
\end{proof}

Marcolli \cite{Marcolli1} proved a weaker result in the case $b_1(Y)>0$,
allowing $\eta$ to vary freely among the closed $2$--forms.

\subsection{Regularity of moduli spaces}

The following lemma will provide us with suitable Banach spaces of
perturbation forms.

\begin{lemma}\label{banach-construct}
Let $X$ be a smooth $n$--manifold, $K\subset X$ a compact, codimension~$0$
submanifold, and $E\to X$ a vector bundle. Then there exists a Banach space
$W$ consisting of smooth sections of $E$ supported in $K$, such
that the following hold:
\begin{enumerate}
\item[\rm(i)]The natural map $W\to\Ga(E|_K)$ is continuous with respect to the
$C^\infty$ topology on $\Ga(E|_K)$.
\item[\rm(ii)]For every point $x\in\text{int}(K)$ and every $v\in E_x$ there
exists a section $s\in\Ga(E)$ with $s(x)=v$ and a smooth embedding
$g\col \R^n\to X$ with $g(0)=x$ such that for arbitrarily small $\eps>0$
there are elements of $W$ of the form $fs$ where $f\col X\to[0,1]$
is a smooth function which vanishes outside $g(\R^n)$ and satisfies
\[f(g(z))=\begin{cases}
0,&|z|\ge2\eps,\\
1,&|z|\le\eps.
\end{cases}\]
\end{enumerate}
\end{lemma}

\begin{proof}
Fix connections in $E$ and $TX$, and a Euclidean metric on $E$.
For any sequence $a=(a_0,a_1,\dots)$
of positive real numbers and any $s\in\Ga(E)$ set
\[\|s\|_a=\sup_{0\le k<\infty}a_k\|\nabla^ks\|_\infty\]
and
\[W_a=\{s\in\Ga(E)\st\text{supp}(s)\subset K,\;\|s\|_a<\infty\}.\]
Then $W=W_a$, equipped with the norm $\|\sdot\|_a$, clearly satisfies
(i) for any $a$. We claim that one can choose $a$ such that (ii) also holds.
To see this, first observe that
there is a finite dimensional subspace $V\subset\Ga(E)$ such that
\[V\to E_x,\quad s\mapsto s(x)\]
is surjective for every $x\in K$. Fix a smooth function $b\col \R\to[0,1]$ 
satisfying
\[b(t)=\begin{cases}
1,&t\le1,\\
0,&t\ge4.
\end{cases}\]
We use functions $f$
that in local coordinates have the form
\[f_r(z)=b(r|z|^2),\]
where $r\gg0$. Note that for each $k$ there is a bound
$\|f_r\|_{C^k}\le\const r^k$ where the constant is independent of $r\ge1$.
It is now easy to see that a suitable sequence $a$ can be found.
\end{proof}

In the next two propositions, $X,\vec\al,\mu$ will 
be as in Subsection~\ref{intro:comp-single}. Let $K\subset X$ be any 
non-empty compact codimension~$0$ submanifold. Let $W$ be a Banach space
of smooth self-dual $2$--forms on $X$ supported in $K$, as provided by
Lemma~\ref{banach-construct}.
The following proposition will be used in the proof of Theorem~\ref{sw-thm2}.

\begin{prop}\label{transv1}
In the above situation, let $G$ be set of all $\nu\in W$ such that
all irreducible points of the moduli space $M(X;\vec\al;\mu+\nu;0)$
are regular (here $\prt_j=0$ for each $j$).
Then $G\subset W$ is residual, hence dense.
\end{prop}

There is another version of this proposition where
$W$ is replaced with the Fr\'echet space of
all smooth self-dual $2$--forms on $X$ supported in $K$, at least if one
assumes that (O1) holds for each pair $Y_j,\eta_j$ and
that (A) holds for $X,\ti\eta_j,\lla_j$.
The reason for the extra assumptions is
that the proof then seems to require global compactness results
(cf the proof of Lemma~\ref{irr-crit-nondeg}).

\begin{proof}
We may assume $X$ is connected. Let $\tsw$ be
as in Subsection~\ref{subsec:modsp}. Then
\[(S,\nu)\mapsto\tsw(S,\mu+\nu,0)\]
defines a smooth map
\[f\col \cc^*\times W\to\lw pw(X;i\La^+\oplus\bs^-),\]
where $\cc^*=\cc^*(X;\vec\al)$.
We will show that $0$ is a regular value of $f$. Suppose
$f(S,\nu)=0$ and write $S=(A,\Phi)$. We must show that the derivative
$P=Df(S,\nu)$ is surjective. Because of the gauge 
equivariance of $f$ we may assume that
$S$ is smooth. Let $P_1$ denote the derivative of
$f(\,\cdot\,,\nu)$ at $S$. Since the image of $P_1$ in $\lw pw$ is closed 
and has finite codimension, the same holds for the image of $P$. Let
$p'$ be the exponent conjugate to $p$ and suppose
$(z,\psi)\in\lw{p'}{-w}(X;i\La^+\oplus\bs^-)$ is $L^2$ orthogonal to
the image of $P$, ie
\[\int_X\la P(a,\phi,\nu'),(z,\psi)\ra=0\]
for all $(a,\phi)\in\llw pw1$ and $\nu'\in W$. Taking $\nu'=0$
we see that $P_1^*(z,\psi)=0$. Since $P_1^*$ has injective symbol,
$z,\psi$ must be smooth.
On the other hand, taking $a,\phi=0$ and varying $\nu'$ we find that
$z|_K=0$ by choice of $W$.
By assumption, $\Phi$ is not identically zero. Since $D_A\Phi=0$,
the unique continuation theorem in 
\cite{JKazdan1} applied to $D_A^2$ says that $\Phi$
cannot vanish in any non-empty open set.
Hence $\Phi$ must be non-zero at some point $x$ in the interior of $K$.
Varying $a$
alone near $x$ one sees that $\psi$ vanishes in some neighbourhood of $x$.
But $P_1P_1^*$ has the same symbol as
$D_A^2\oplus d^+(d^+)^*$, so another application of the same unique
continuation theorem shows that $(z,\psi)=0$. Hence $P$ is surjective.

Consider now the vector bundle
\[E=\cc^*\times_{\cg}\lw pw(X;i\La^+\oplus\bs^-)\]
over $\cb^*$. The map $f$ defines a smooth section $\si$ of the bundle
\[E\times W\to\cb^*\times W.\]
Because of the local slice theorem
and the gauge equivariance of $f$, the fact that $0$ is a regular value of
$f$ means precisely that $\si$ is transverse to the zero-section.
Since $\si(\,\cdot\,,\nu)$ is a Fredholm section of $E$ for
any $\nu$, the proposition follows by another application of the
Sard--Smale theorem.
\end{proof}

We will now establish transversality results for moduli spaces of the 
form $M(X,\vec\al)$ or $M(\al,\beta)$ involving perturbations of the kind
discussed in Subsection~\ref{subsec:pert}. For the time being we limit 
ourselves to the case where the $3$--manifolds $Y,Y_j$ are all rational
homology spheres. We will use functions $h_S$ that are a small modification of
those in \cite{Fr1}. To define these, let $Y$ be a closed Riemannian
$\spc$ $3$--manifold satisfying $b_1(Y)=0$,
and $\cd$ the Chern--Simons--Dirac functional on $Y$
defined by some closed $2$--form $\eta$. 
Choose a smooth, non-negative function $\chi\col \R\to\R$
which is supported in the interval $(-\frac14,\frac14)$
and satisfies $\int\chi=1$. If $S$ is any $L^2_1$ configuration over
a band $(a-\frac14,b+\frac14)$ where $a\le b$ define the smooth function
$\tcd_S\col [a,b]\to\R$ by
\[\tcd_S(T)=\int_{\R}\chi(T-t)\cd(S_t)\,dt,\]
where we interpret the right hand side as an integral over $\ry$.
A simple exercise, using the Sobolev embedding theorem, shows that if
$S_n\to S$ weakly in $L^2_1$ over $(a-\frac14,b+\frac14)\times Y$ then 
$\tcd_{S_n}\to\tcd_S$ in $C^\infty$ over $[a,b]$.

Choose a smooth function
$\cutoff\col \R\to\R$ with the following properties:
\begin{itemize}
\item $\cutoff'>0$,
\item $\cutoff$ and all its derivatives are bounded,
\item $\cutoff(t)=t$ for all critical values $t$ of $\cd$,
\end{itemize}
where $\cutoff'$ is the derivative of $\cutoff$.
The last condition is added only for convenience.

For any $L^2_1$ configuration $S$ over
$(a-\frac12,b+\frac12)\times Y$ with $a\le b$ define
\[h_S(t)=\int_{\R}\chi(t_1)\cutoff(\tcd_S(t-t_1))dt_1.\]
It is easy to verify that $h_S$ satisfies the properties (P1)--(P3). 

It remains to choose $\Xi$ and $\Prt$.
Choose one compact subinterval (with non-empty interior)
of each bounded connected component of $\R\setminus\text{crit}(\cd)$,
where $\text{crit}(\cd)$ is the set of critical values of $\cd$.
Let $\Xi$ be the union of these compact subintervals. 
Let $\Prt=\Prt_Y$ be a Banach space of $2$--forms on $\R\times Y$
supported in $\Xi\times Y$ as provided by Lemma~\ref{banach-construct}.

We now return to the situation described in the paragraph preceding
Proposition~\ref{transv1}. Let $W'\subset W$ be the open subset consisting
of those elements $\nu$ that satisfy $\|\nu\|_{C^1}<1$. Let $\Pkd$ denote the
set of all $\vec\prt=(\prt_1,\dots,\prt_r)$ where $\prt_j\in\Prt_{Y_j}$
and $\|\prt_j\|_{C^1}<\del$ for each $j$.

\begin{prop}\label{transv2}
Suppose each $Y_j$ is a rational homology sphere and
$K\subset X\endt0$. 
Then there exists a $\del>0$ such that the following holds.
Let $G$ be the set of all
$(\nu,\vec\prt)\in W'\times\Pid$ such that every irreducible point of
the moduli space $M(X;\vec\al;\mu+\nu;\vec\prt)$ is regular. Then
$G\subset W'\times\Pid$ is residual,  hence dense.
\end{prop}
It seems necessary here to let $\vec\prt$ vary as well, because
if any of the $\prt_j$ is non-zero then the linearization of the monopole map
is no longer a differential operator, and it is not clear whether
one can appeal to unique continuation as in the proof of
Proposition~\ref{transv1}. 

\begin{proof}
To simplify notation assume $r=1$ and set $Y=Y_1$, $\al=\al_1$ etc.
(The proof in the general case is similar.)
Note that (A) is trivially satisfied, since each $Y_j$ is a 
rational homology sphere. Therefore, by Propositions~\ref{Dcd} and
\ref{prop:local-energy-bound2}, if $\del>0$ is sufficiently small
then for any
$(\nu,\prt)\in W'\times\Pkd$ and $[S]\in M(X;\al;\mu+\nu;\prt)$ one has
that either
\begin{enumerate}
\item[(i)]$[S_t]=\al$ for $t\ge0$, or
\item[(ii)]$\prtl_t\cd(S_t)<0$ for $t\ge0$.
\end{enumerate}

As in the proof of Proposition~\ref{transv1} it suffices 
to prove that
$0$ is a regular value of the smooth map
\begin{align*}
\ti f\col \cc^*\times W'\times\Pkd & \to\lw pw,\\
(S,\nu,\prt) & \mapsto\tsw(S,\mu+\nu,\prt).
\end{align*}
The smoothness of the perturbation term $g(S,\prt)=\fq h_{S,\prt}$ follows
 from the smoothness of the map \Ref{sprt}, for by (P1) there exist a $t_0$
and a neighbourhood $U\subset\cc$ of $S$ such
that $h_{S'}(t)\not\in\Xi$ for all $t>t_0$ and $S'\in U$.

Now suppose $\ti f(S,\nu,\prt)=0$ and $(z,\psi)\in\lw{p'}{-w}$ is orthogonal
to the image of $D\ti f(S,\nu,\prt)$. We will show that $z$ is orthogonal
to the image of $T=Dg(S,\prt)$, or equivalently, that $(z,\psi)$ is orthogonal
to the image of $Df(S,\nu)$, where $f=\ti f-g$ as before.
The latter implies $(z,\psi)=0$ by the proof of
Proposition~\ref{transv1}.

Let $h_S\col [\frac12,\infty)\to\R$ be defined in terms of the restriction of
$S$ to $\rpy$. If $h_S(J)\subset\R\setminus\Xi$ for some compact interval
$J$ then by (P1) one has that
$h_{S'}(J)\subset\R\setminus\Xi$ for all $S'$ in some neighbourhood of $S$ 
in $\cc$. Therefore, all elements of $\im(T)$ vanish on
$h_S\inv(\R\setminus\Xi)\times Y$. 

We now digress to recall that if $u$ is any locally integrable function
on $\R^n$ then the complement of the Lebesgue set of $u$ has measure zero,
and if $v$ is any continuous function on $\R^n$ then
any Lebesgue point of $u$ is also a Lebesgue point of $uv$.
The notion of Lebesgue set also makes sense for sections $\tau$ of 
a vector bundle of finite rank over a finite dimensional smooth manifold $M$.
In that case a point $x\in M$ is called a Lebesgue point of $\tau$ if it is
a Lebesgue point in the usual sense for some (hence any) choice of local
coordinates and local trivialization of the bundle around $x$.

Returning to our main discussion, there are now two cases:
If (i) above holds then
$h_S(t)=\cd(\al)\not\in\Xi$ for $t\ge\frac12$,
whence $T=0$ and we are done (recall
the overall assumption $\fq\inv(0)=X\endt{\frac32}$ made
in Subsection~\ref{subsec:modsp}).
Otherwise (ii) must hold. In that case we have $\prtl_t\cutoff(\tcd(t))<0$
for $t\ge\frac14$ and $\prtl_th_S(t)<0$ for $t\ge\frac12$. Since $z$
is orthogonal to $\fq h_{S,\prt'}$ for all $\prt'\in\Prt_Y$
we conclude that $z(t,y)=0$ for every Lebesgue point $(t,y)$ of $z$
with $t>\frac32$ and $h_S(t)\in\itr(\Xi)$.
Since $h_S\inv(\prtl\Xi)\cap(\frac32,\infty)$ is a finite set, $z$ must 
vanish almost everywhere in $[h_S\inv(\Xi)\cap(\frac32,\infty)]\times Y$.
Combining this with our earlier result we deduce that
$z$ is orthogonal to $\im(T)$.
\end{proof}

In the next proposition (which is similar to \cite[Proposition~5]{Fr1})
let $\Pid$ be as above with $r=1$, and set $Y=Y_1$.

\begin{prop}
In the situation of Subsection~\ref{intro-csd}, suppose $Y$ is a rational
homology sphere and $\al,\beta\in\fl_Y=\ti\fl_Y$.
Then there exists a $\del>0$ such that the following holds.
Let $G$ be the set of all $\prt\in\Pid$ such that
every point in $M(\al,\beta;\prt)$ is regular. Then
$G\subset\Pid$ is residual,  hence dense.
\end{prop}

\begin{proof}
If $\al=\beta$ then an application of
Proposition~\ref{Dcd} shows that if $\|\prt\|_{C^1}$ is sufficiently small then
$M(\al,\beta;\prt)$ consists of a single point represented by $\ual$,
which is regular because $\al$ is non-degenerate.

If $\al\neq\beta$ and $\|\prt\|_{C^1}$ is sufficiently small then
for any $[S]\in M(\al,\beta;\prt)$ one has $\prtl_t\cd(S_t)<0$ for all $t$.
Moreover, the moduli space contains no reducibles, since $\al,\beta$
cannot both be reducible.
The proof now runs along the same lines as that of Proposition~\ref{transv2}.
Note that the choice of $\Xi$ is now essential: it ensures that
$\im(h_S)=(\cd(\al),\cd(\beta))$ contains interior points of $\Xi$.
\end{proof}

\section{Proof of Theorems~\ref{sw-thm1} and \ref{sw-thm2}}
\label{proofs12}

In these proofs we will only use genuine monopoles.

{\bf Proof of Theorem~\ref{sw-thm1}}\qua
We may assume $Y$ is connected.  Let $\eta$ be a closed non-exact
$2$--form on $Y$ which is the restriction of a closed form on $Z$. Let $Y$
have a metric of positive scalar curvature. 
If $s\neq0$ is a small real number then $\csd{s\eta}$ will have no irreducible
critical points, by the \textit{a priori} estimate on the spinor fields and the
positive scalar curvature assumption. If in addition
$s[\eta]-\pi c_1(\cll_Y)\neq0$ then $\csd{s\eta}$ will have no reducible
critical points either.

Choose a $\spc$ Riemannian $4$--manifold $X$ as in
Subsection~\ref{intro:comp-neck}, with $r=1$, $r'=0$, such that there exists
a diffeomorphism $X\gl\to Z$ which maps $\{0\}\times Y_1$ isometrically onto $Y$.
Let $\eta_1$ be the pull-back of $s\eta$.
Then (B1) is satisfied (but perhaps not (B2)), so it follows from
Theorem~\ref{thm:neck-comp} that $M(\xt)$ is empty for $T\gg0$.
\endproof

We will now define an invariant $h$ for closed $\spc$ $3$--manifolds $Y$ that
satisfy $b_1(Y)=0$ and admit metrics with positive scalar curvature.
Let $g$ be such a metric on $Y$. 
Recall that for the unperturbed Chern--Simons--Dirac functional $\cd$ the
space $\fl_Y$ of critical points modulo gauge equivalence consists of
a single point $\theta$, which is reducible. Let $(B,0)$ be a representative
for $\theta$. Let $Y_1,\dots,Y_r$ be the connected components of $Y$ and
choose a $\spc$ Riemannian $4$--manifold $X$ with tubular
ends $\rpy_j$, $j=1,\dots,r$  (in the sense of Subsection~\ref{intro:comp-single}) and
a smooth spin connection $A$ over $X$ such that the restriction of $\check A$
to $\rpy$ is equal to the pull-back of $\check B$. (The notation here
is explained in Subsection~\ref{subsec:spc-str}.) Define
\begin{align*}
h(Y,g)&=\ind_{\co}(D_A)-\frac18(c_1(\cll_X)^2-\si(X))\\
&=\frac12(\dim\,M(X;\theta)-d(X)+b_0(X)), 
\end{align*}
where $D_A\col L^2_1\to L^2$, `dim' is the expected dimension, and $d(X)$ is the 
quantity defined in Subsection~\ref{intro:vanishing}.
Since $\ind_{\co}(D_A)=\frac18(c_1(\cll_X)^2-\si(X))$
when $X$ is closed, it follows easily from the addition formula for the index
(see \cite[Proposition~3.9]{D5}) that $h(Y,g)$ is independent of $X$ and that
\[h(-Y,g)=-h(Y,g).\]
Clearly,
\[h(Y,g)=\sum_jh(Y_j,g_j),\]
where $g_j$ is the restriction of $g$ to $Y_j$.
To show that $h(Y,g)$ is independent of $g$ we may therefore assume $Y$ is
connected. Suppose
$g'$ is another positive scalar curvature metric on $Y$ and consider the $\spc$
Riemannian  manifold $X=\ry$ where the metric agrees with
$1\times g$ on $(-\infty,-1]\times Y$ and with $1\times g'$ on
$[1,\infty)\times Y$.
Theorem~\ref{thm:sing-comp} and Proposition~\ref{prop:chain-conv}~(ii)
say that $M(X;\theta,\theta)$ is compact. This moduli space contains one
reducible point, and
arguing as in the proof of \cite[Theorem~6]{Fr1} one sees that the moduli
space 
must have non-positive (odd) dimension. Thus,
\[h(Y,g')+h(-Y,g)=\frac12(\dim\,M(X;\theta,\theta)+1)\le0.\]
This shows $h(Y)=h(Y,g)$ is independent of $g$. 

{\bf Proof of Theorem~\ref{sw-thm2}}\qua
Let each $Y_j$ have a positive scalar
curvature metric. Choose a $\spc$ Riemannian $4$--manifold $X$ as in
Subsection~\ref{intro:comp-neck}, with $r'=0$ and with the same $r$, 
such that there exists
a diffeomorphism $f\col X\gl\to Z$ which maps $\{0\}\times Y_j$
isometrically onto $Y_j$.
Then (B1) is satisfied (but perhaps not (B2)). Let $X_0$ be the component of $X$ such
that $W=f((X_0)\endt1)$.
For each $j$ set $\eta_j=0$ and let $\al_j\in\fl_{Y_j}$ be
the unique (reducible) critical point. Choose a reference connection $A\rr$
as in Subsection~\ref{subsec:modsp} and set $A_0=A\rr|_{X_0}$. Since
each $\al_j$ has representatives of the form $(B,0)$ where $\check B$ is flat
it follows that $\hat F(A\rr)$ is compactly supported. In the following,
$\mu$ will denote the (compactly supported) perturbation $2$--form on $X$
and $\mu_0$ its restriction to $X_0$.

Let $\harm^+$ be the space of self-dual closed $L^2$ $2$--forms on $X_0$. Then
$\dim\,\harm^+=b^+_2(X_0)>0$, so $\harm^+$ contains a non-zero element $z$.
By unique continuation for harmonic forms we can find a smooth $2$--form 
$\mu_0$ on $X_0$, supported in any given small ball, such that 
$\hatfpl(A_0)+i\mu^+_0$ is not $L^2$ orthogonal to $z$. (Here $\hatfpl$
is the self-dual part of $\hat F$.) Then
\[\hatfpl(A_0)+i\mu^+_0\not\in\text{im}(d^+\col \llw pw1\to\lw pw),\]
where $w$ is the weight function used in the definition of the configuration space.
Hence $M(X_0;\vec\al)$ contains no reducible monopoles. After perturbing
$\mu_0$ in a small ball we can arrange that $M(X_0;\vec\al)$ is transversally
cut out as well, by Proposition~\ref{transv1}.

To prove (i), recall that
\begin{equation}\label{dwsum}
\dim\,M(X_0;\vec\al)=d(W)-1+2\sum_jh(Y_j),
\end{equation}
so the inequality in (i) simply says that
\[\dim\,M(X_0;\vec\al)<0,\]
hence $M(X_0;\vec\al)$ is empty. Since there are no other moduli spaces over $X_0$,
it follows from Theorem~\ref{thm:neck-comp} that $M(\xt)$ is empty when $\min_jT_j\gg0$.

We will now prove (ii). If $M(\xt)$ has odd or negative dimension then there
is nothing to prove, so suppose this dimension is $2m\ge0$.
Since $M(X_0;\vec\al)$ contains no reducibles we deduce from Theorem~\ref{thm:neck-comp}
that $M(\xt)$ is also free of reducibles when $\min_jT_j$ is sufficiently large.
Let $\BB\subset X_0$ be a compact $4$--ball and $\cb^*(\BB)$
the Banach manifold of irreducible $L^p_1$ configurations over $\BB$ modulo
$L^p_2$ gauge
transformations. Here $p>4$ should be an even integer to ensure the existence
of smooth partitions of unity.  Let
$\bl\to\cb^*(\BB)$ be the natural complex line bundle associated to some
base-point in $\BB$, and $s$ a generic
section of the $m$--fold direct sum $m\bl$. For $\min_jT_j\gg0$ let
\[S\tu\subset M(\xt),\qquad S_0\subset M(X_0;\vec\al)\]
be the subsets consisting of those
elements $\om$ that satisfy $s(\om|_{\BB})=0$. By assumption, $S_0$ is a 
submanifold of codimension $2m$. For any $T$ for which $S\tu$ is transversely cut out
the Seiberg--Witten invariant of $Z$ is equal to the number of points in $S\tu$
counted with sign. Now, the inequality in (ii) is equivalent to
\[d(W)+2\sum_jh(Y_j)<d(Z)=2m+1,\]
which by \Ref{dwsum} gives
\[\dim\,S_0=\dim\,M(X_0;\vec\al)-2m<0.\]
Therefore, $S_0$ is empty. By Theorem~\ref{thm:neck-comp}, 
$S\tu$ is empty too when $\min_jT_j\gg0$, hence $\text{SW}(Z)=0$.
\endproof

\appendix

\section{Patching together local gauge transformations}

In the proof of Lemma~\ref{lemma:loc-comp}
we encounter sequences $S_n$ of configurations
such that for any point $x$ in the base-manifold there is a sequence $v_n$
of gauge transformations defined in a neighbourhood of $x$ such that
$v_n(S_n)$ converges (in some Sobolev norm) in a (perhaps smaller)
neighbourhood of $x$. The problem then is to find a sequence $u_n$ of
global gauge transformations such that $u_n(S_n)$ converges globally.
If $v_n,w_n$ are two such sequences of local gauge transformations then
$v_nw_n\inv$ will be bounded in the appropriate Sobolev norm, so the
problem reduces to the lemma below.

This issue was discussed by Uhlenbeck in \cite[Section~3]{U1}. Our
approach has the advantage that it does not involve any ``limiting bundles''.

\begin{lemma}\label{lemma:patch}
Let $X$ be a Riemannian manifold and $P\to X$ a principal
$G$--bundle, where $G$ is a compact subgroup of some matrix
algebra $M_r(\R)$. Let $M_r(\R)$ be equipped with
an $\text{Ad}_G$--invariant inner product, and fix a connection in the
Euclidean vector bundle $E=P\times_{\text{Ad}_G}M_r(\R)$ (which we use to
define Sobolev norms of automorphisms of $E$).
Let $\{U_i\}_{i=1}^\infty$, $\{V_i\}_{i=1}^\infty$ be open covers
of $X$ such that $U_i\Subset V_i$ for each $i$.
We also assume that each $V_i$ is the interior of a compact codimension~$0$
submanifold of $X$, and that $\prtl V_i$ and $\prtl V_j$ intersect 
transversally for all $i\neq j$.
For each $i$ and $n=1,2,\dots$ let $\vin$ be a continuous
automorphism of $P|_{V_i}$.
Suppose $\vin\vjn\inv$ converges uniformly over $V_i\cap V_j$
for each $i,j$ (as maps into $E$).
Then there exist
\begin{itemize}
\item a sequence of positive integers $n_1\le n_2\le\cdots$,
\item for each positive integer $k$ an open subset $W_k\subset X$ with
\[\bigcup_{i=1}^kU_i\Subset W_k\subset\bigcup_{i=1}^kV_i,\]
\item for each $k$ and $n\ge n_k$ a continuous automorphism $\wkn$ of $P|_{W_k}$,
\end{itemize}
such that
\begin{enumerate}
\item[\rm(i)]If $1\le j\le k$ and $n\ge n_k$ then
$w_{j,n}=\wkn$ on $\bigcup_{i=1}^jU_i$,
\item[\rm(ii)]For each $i,k$ the sequence
$\wkn\vin\inv$ converges uniformly over $W_k\cap V_i$,
\item[\rm(iii)]If $1\le p<\infty$, and $m>\frac np$ is an integer such that $\vin\in L^p_{m,\loc}$
for all $i,n$ then $\wkn\in L^p_{m,\loc}$ for all $k$ and $n\ge n_k$. If in addition
\[\sup_n\|\vin\vjn\inv\|_{L^p_m(V_i\cap V_j)}<\infty\qquad
\text{for all $i,j$}\]
then
\[\sup_{n\ge n_k}\|\wkn\vin\inv\|_{L^p_m(W_k\cap V_i)}<\infty\qquad
\text{for all $k,i$}.\]
\end{enumerate}
\end{lemma}

The transversality condition ensures that the Sobolev embedding theorem
holds for $V_i\cap V_j$ (see \cite{Adams}). Note that this condition can
always be achieved by shrinking the $V_i$'s a little.

\begin{proof}
Let $N'\subset LG$ be a small $\text{Ad}_G$ invariant open neighbourhood
of $0$. Then $\exp\col LG\to G$ maps $N'$ diffeomorphically onto an
open neighbourhood $N$ of $1$. Let $f\col N\to N'$ denote the inverse map.
Let $\aut(P)$ the
bundle of fibre automorphisms of $P$ and $\g_P$ the corresponding bundle of 
Lie algebras. Set $\mathbf{N}=P\times_{\text{Ad}_G}N\subset\g_P$
and let $\exp\inv\col \mathbf{N}\to\aut(P)$
be the map defined by $f$.

Set $w_{1,n}=v_{1,n}$ and $W_1=V_1$. Now suppose $\wkn$, $W_k$
have been chosen for $1\le k<\ell$, where $\ell\ge2$, such that (i)--(iii)
hold for these values of $k$. 
Set $z_n=w_{\ell-1,n}(v_{\ell,n})\inv$ on $W_{\ell-1}\cap V_\ell$. According to
the induction hypothesis the sequence $z_n$ converges uniformly over
$W_{\ell-1}\cap V_\ell$, hence there exists an integer
$n_\ell\ge n_{\ell-1}$ such that $y_n=(z_{n_\ell})\inv z_n$ 
takes values in $\mathbf{N}$ for $n\ge n_\ell$.

Choose an open subset $\cw\subset X$ which is the
interior of a compact codimension~$0$ submanifold of $X$, and which satisfies
\[\bigcup_{i=1}^{\ell-1} U_i\Subset\cw\Subset W_{\ell-1}.\]
We also require that $\prtl\cw$ intersect $\prtl V_i\cap\prtl V_j$
transversally for all $i,j$. (For instance, one can take
$\cw=\al\inv([0,\eps])$ for suitable $\eps$, where $\al\col X\to[0,1]$ is any
smooth function with $\al=0$ on $\cup_{i=1}^{\ell-1} U_i$
and $\al=1$ on $W_{\ell-1}$.) Choose also
a smooth, compactly supported function $\phi\col W_{\ell-1}\to\R$ with
$\phi|_\cw=1$. Set $W_\ell=\cw\cup V_\ell$ and 
for $n\ge n_\ell$ define an automorphism $w_{\ell,n}$ of $P|_{W_\ell}$ by
\[w_{\ell,n}=\begin{cases}
w_{\ell-1,n} & \text{on $\cw$,}\\
z_{n_\ell}\exp(\phi\exp\inv y_n)v_{\ell,n} &
   \text{on $W_{\ell-1}\cap V_\ell$,}\\
z_{n_\ell}v_{\ell,n} & \text{on $V_\ell
\setminus\text{supp}(\phi)$.}
\end{cases}\]
Then (i)--(iii) hold for $k=\ell$ as well. To see that (iii) holds,
note that our transversality assumptions guarantee that the Sobolev
embedding theorem holds for $W_{\ell-1}\cap V_\ell$ and for all
$V_i\cap V_j$. Since $mp>n$, $L^p_m$ is therefore a Banach algebra for
these spaces (see \cite{Adams}). Recalling the proof of this fact, and
the behaviour of $L^p_m$ under composition with smooth maps on the left
(see \cite[p\,184]{McDuff-Salamon1}), one obtains (iii).
\end{proof}

\section{A quantitative inverse function theorem}

In this section $E,E'$ will be Banach spaces. We denote by $\cb(E,E')$
the Banach space of bounded operators from $E$ to $E'$. 
If $T\in\cb(E,E')$ then $\|T\|=\sup_{\|x\|\le1}\|Tx\|$.
If $U\subset E$ is open and $f\col U\to E'$ smooth then $Df(x)\in\cb(E,E')$
is the derivative of $f$ at $x\in U$. The second derivative
$D(Df)(x)\in\cb(E,\cb(E,E'))$ is usually written $D^2f(x)$ and 
can be identified with the symmetric bilinear map $E\times E\to E'$ given by
\[D^2f(x)(y,z)=\left.\frac{\prtl^2}{\prtl s\prtl t}\right|_{(0,0)}f(x+sy+tz).\]
The norm of the second derivative is
\[\|D^2f(x)\|=\sup_{\|y\|,\|z\|\le1}\|D^2f(x)(y,z)\|.\]
For $r>0$ let
\[B_r=\{x\in E\st\|x\|<r\}.\]
\begin{lemma}\label{quant-inv1}
Let $\eps,M>0$ be positive real numbers
such that $\eps M<1$, and suppose $f\col B_\eps\to E$ is a smooth map satisfying
\[f(0)=0;\quad Df(0)=I;\quad \|D^2f(x)\|\le M\quad\text{for $x\in B_\eps$}.\]
Then $f$ restricts to a diffeomorphism $f\inv B_{\frac\eps2}\oset\approx\to
B_{\frac\eps2}$.
\end{lemma}

The conclusion of the lemma holds even when $\eps M=1$, see
Proposition~\ref{quant-inv2} below.

\begin{proof}
The estimate on $D^2f$ gives
\begin{equation}\label{Df-est}
\|Df(x)-I\|=\|Df(x)-Df(0)\|\le M\|x\|.
\end{equation}
Therefore the map
\[h(x)=f(x)-x=\int_0^1(Df(tx)-I)x\,dt\]
satisfies
\begin{align*}
\|h(x)\|&\le\frac M2\|x\|^2,\\
\|h(x_2)-h(x_1)\|&\le\eps M\|x_2-x_1\|
\end{align*}
for all $x,x_1,x_2\in B_\eps$. Hence for every
$y\in B_{\frac\eps2}$ the assignment $x\mapsto y-h(x)$ defines a map
$B_\eps\to B_\eps$ which has a unique fix-point.
In other words, $f$ maps
$f\inv B_{\frac\eps2}$ bijectively onto $B_{\frac\eps2}$. Moreover,
$Df(x)$ is an isomorphism for every $x\in B_\eps$, by \Ref{Df-est}.
Applying the contraction mapping
argument above to $f$ around an arbitrary point in $B_\eps$ then shows
that $f$ is an open map. It is then a simple exercise to prove that
the inverse $g\col B_{\frac\eps2}\to f\inv B_{\frac\eps2}$
is differentiable and $Dg(y)=(Df(g(y)))\inv$ (see \cite[8.2.3]{Dieud}).
Repeated application of the chain rule then shows that $g$ is smooth.
\end{proof}

For $r>0$ let $B_r\subset E$ be as above, and define
$B'_r\subset E'$ similarly.

\begin{prop}\label{quant-inv2}
Let $\eps,M$ be positive real numbers and $f\col B_\eps\to E'$ a smooth map
such that $f(0)=0$, $L=Df(0)$ is invertible, and
\[\|D^2f(x)\|\le M\quad\text{for all $x\in B_\eps$.}\]
Set $\ka=\|L\inv\|\inv-\eps M$ and $\eps'=\eps\|L\inv\|\inv$.
Then the following hold:
\begin{enumerate}
\item[\rm(i)]If $\ka\ge0$ then $f$ is a diffeomorphism onto an open subset of
$E'$ containing $B'_{\eps'/2}$.
\item[\rm(ii)]If $\ka>0$ and $g\col B'_{\eps'/2}\to B_\eps$ is the smooth map
satisfying $f\circ g=I$ then for all $x\in B_\eps$ and
$y\in B'_{\eps'/2}$ one has
\[\|Df(x)\inv\|,\|Dg(y)\|<\ka\inv,\quad \|D^2g(y)\|<M\ka^{-3}.\]
\end{enumerate}
\end{prop}

The reader may wish to look at some simple example (such as a quadratic
polynomial) to understand the various ways in which these results are optimal.

\proof
(i)\qua For every $x\in B_\eps$ we have
\[\|Df(x)L\inv-I\|\le\|Df(x)-L\|\cdot\|L\inv\|<\eps M\|L\inv\|\le1,\]
hence $Df(x)$ is invertible. Thus $f$ is a local diffeomorphism by 
Lemma~\ref{quant-inv1}. Set $h(x)=f(x)-Lx$. If $x_1,x_2\in B_\eps$ and
$x_1\neq x_2$ then
\[\|f(x_2)-f(x_1)\|\ge\|L(x_2-x_1)\|-\|h(x_2)-h(x_1)\|>\ka\|x_2-x_1\|,\]
hence $f$ is injective. By choice of $\eps'$ the map
\[\ti f=f\circ L\inv\col B'_{\eps'/2}\to E'\]
is well defined, and for every $y\in B'_{\eps'/2}$ one has
\[\|D^2\ti f(y)\|\le M\|L\inv\|^2.\]
Because
\[\eps' M\|L\inv\|^2=\eps M\|L\inv\|\le1,\]
Lemma~\ref{quant-inv1} says that the image of $\ti f$ contains every ball
$B'_{\del/2}$ with $0<\del<\eps'$, hence also $B'_{\eps'/2}$.

(ii)\qua Set $c=I-Df(x)L\inv$. Then
\[Df(x)\inv=L\inv\sum_{n=0}^\infty c^n,\]
hence
\[\|Df(x)\inv\|\le\frac{\|L\inv\|}{1-\|c\|}
<\frac{\|L\inv\|}{1-\eps M\|L\inv\|}=\ka\inv.\]
This also gives the desired bound on $Dg(y)=Df(g(y))\inv$.

To estimate $D^2g$, let $\text{Iso}(E,E')\subset\cb(E,E')$ be the open subset
of invertible operators, and let $\iota\col \text{Iso}(E,E')\to\cb(E',E)$ be the
inversion map:
$\iota(a)=a\inv$. Then $\iota$ is smooth, and its derivative is given by
\begin{equation*}
D\iota(a)b=-a\inv b a\inv,
\end{equation*}
see \cite{Dieud}. The chain rule says that
\begin{align*}
Dg&=\iota\circ Df\circ g,\\
D(Dg)(y)&=D\iota(Df(g(y)))\circ D(Df)(g(y))\circ Dg(y).
\end{align*}
This gives
$$\|D(Dg)(y)\|\le\|Df(g(y))\inv\|^2\cdot\|D(Df)(g(y))\|\cdot\|Dg(y)\|
<\ka^{-2}\cdot M\cdot\ka\inv.\eqno{\qed}$$

\end{document}